\documentclass[journal,oneside,web]{template}
\usepackage{generic} 
\usepackage{cite}
\usepackage{amsmath,amssymb,amsfonts}
\usepackage{algorithmic}
\usepackage{graphicx}
\usepackage{textcomp}
\usepackage{latexsym,dsfont}
\usepackage{mathrsfs}
\usepackage{hyperref}
\usepackage{subcaption}
\usepackage{booktabs}
\usepackage{multirow}

\newcommand{\hf}{{\frac 12}}

\DeclareMathOperator*{\argmax}{arg\,max}

\newcommand{\bfa}{\boldsymbol{a}}
\newcommand{\bfb}{\boldsymbol{b}}

\newcommand{\bfp}{\boldsymbol{p}}

\newcommand{\bfs}{\boldsymbol{s}}

\newcommand{\bfu}{\boldsymbol{u}}

\newcommand{\bfw}{\boldsymbol{w}}
\newcommand{\bfx}{\boldsymbol{x}}
\newcommand{\bfy}{\boldsymbol{y}}
\newcommand{\bfz}{\boldsymbol{z}}

\newcommand{\bfA}{\boldsymbol{A}}

\newcommand{\bfI}{\boldsymbol{I}}

\newcommand{\bfK}{\boldsymbol{K}}

\newcommand{\bfX}{\boldsymbol{X}}

\newcommand{\bfmu}{{\boldsymbol \mu}}
\newcommand{\bfxi}{{\boldsymbol \xi}}

\newcommand{\bfSigma}{{\boldsymbol \Sigma}}

\newcommand{\R}{\ensuremath{\mathds{R}}}
\newcommand{\E}{\ensuremath{\mathds{E}}}

\newtheorem{theorem}{Theorem}
\newtheorem{remark}{Remark}

\def\du{\ensuremath{\mathrm{d}}}

\newcommand{\bfth}{\boldsymbol{\theta}}

\graphicspath{{figures/}{./}}

\newtheorem{assumption}{Assumption}

\newenvironment{mytheorem}[1]{%
    \begin{theorem}#1}{%
    \Endofdef\end{theorem}%
}

\newenvironment{myassumption}[1]{%
    \begin{assumption}#1}{%
    \Endofdef\end{assumption}%
}

\newcommand{\xqed}[1]{%
    \leavevmode\unskip\penalty9999 \hbox{}\nobreak\hfill
    \quad\hbox{\ensuremath{#1}}}
\newcommand{\Endofdef}{\xqed{\diamond}}

\makeatletter  
\def\endthebibliography{%
  \def\@noitemerr{\@latex@warning{Empty `thebibliography' environment}}%
  \endlist
}
\makeatother

\begin{document}
\title{A Neural Network Approach for High-Dimensional Optimal Control\\ Applied to Multi-Agent Path Finding}
\author{Derek Onken, Levon Nurbekyan, Xingjian Li, Samy Wu Fung, Stanley Osher, and Lars Ruthotto
\thanks{This work was supported in part by NSF award DMS 1751636, AFOSR Grants FA95550-20-1-0372 \& FA9550-18-1-0167, AFOSR MURI FA9550-18-1-0502,
Binational Science Foundation Grant 2018209,
US DOE Office of Advanced Scientific Computing Research Field Work Proposal 20-023231, ONR Grants No. N00014-18-1-2527 \& N00014-20-1-2093, a gift from UnitedHealth Group R\&D, and a GPU donation by NVIDIA Corporation. (\textit{Corresponding author: Lars Ruthotto})}%
\thanks{D.~Onken is with the Dept of Computer Science,
Emory University, Atlanta, GA, USA (derek@derekonken.com)}%
\thanks{L.~Nurbekyan and S.~Osher are with the Dept of Mathematics,
UCLA, Los Angeles, CA, USA (lnurbek@math.ucla.edu; sjo@math.ucla.edu)}%
\thanks{X.~Li and L.~Ruthotto are with the Dept of Mathematics,
Emory University, Atlanta, GA, USA (xingjian.li@emory.edu; lruthotto@emory.edu)}%
\thanks{S. Wu Fung is with the Dept of Applied Mathematics and Statistics, Colorado School of Mines, Golden, CO, USA (swufung@mines.edu)}%
}

\maketitle
\begin{abstract}
We propose a neural network approach that yields approximate solutions for high-dimensional optimal control problems and demonstrate its effectiveness using examples from multi-agent path finding.
Our approach yields controls in a feedback form, where the policy function is given by a neural network (NN). Specifically, we fuse the Hamilton-Jacobi-Bellman (HJB) and Pontryagin Maximum Principle (PMP) approaches by parameterizing the value function with an NN.
Our approach enables us to obtain approximately optimal controls in real-time without having to solve an optimization problem. 
Once the policy function is trained, generating a control at a given space-time location takes milliseconds; in contrast, efficient nonlinear programming methods typically perform the same task in seconds.
We train the NN offline using the objective function of the control problem and penalty terms that enforce the HJB equations. Therefore, our training algorithm does not involve data generated by another algorithm. By training on a distribution of initial states, we ensure the controls' optimality on a large portion of the state-space. Our grid-free approach scales efficiently to dimensions where grids become impractical or infeasible. We apply our approach to  several multi-agent collision-avoidance problems in up to 150 dimensions. Furthermore, we empirically observe that the number of parameters in our approach scales linearly with the dimension of the control problem, thereby mitigating the curse of dimensionality.
\end{abstract}


\section{Introduction} \label{sec:introduction}

Decision-making for complex systems using optimal control (OC) has become increasingly relevant yet remains challenging, especially when the state dimension is high and decisions are needed in real-time.
Examples include controlling a swarm of quadcopters~\cite{honig2018trajectory} with collision-avoidance and controlling an unmanned aerial vehicle~\cite{kim2020real,elsayed2020uncertainty,florence2018nanomap,yang2021opt} while reacting to possible wind interference during flight.

We consider real-time OC applications that lead to deterministic, finite time-horizon control problems. 
The speed of generating new controls is critical in these real-time problems where unexpected situations occur during deployment, e.g., wind interference~\cite{ross06realtime,tang2018learning,sanchez2018real,mo21fasttrack,bansal2021deepreach}. 
While nonlinear programming (NLP) methods can provide optimal controls for fixed initial states~\cite{betts98survey}, computation may be too slow for real-time applications: seconds vs milliseconds.
We provide controls in a feedback form, where the policy is given by a neural network (NN). Hence, we generate approximately optimal controls in milliseconds (real-time) without having to solve an optimization problem.

Two of the most common frameworks to solve OC problems
are the Pontryagin Maximum Principle (PMP)~\cite{pontryagin62} and
Hamilton-Jacobi-Bellman (HJB) partial differential equation
(PDE)~\cite{bellman57}.
The PMP is often suitable for high-dimensional problems (Sec.~\ref{sec:PMP}). A local solution method, the PMP finds the optimal policy for a single initial state, so deviations of the system from the optimal trajectory require re-computation of the solution. 
In contrast, the HJB approach is a global solution method suitable for real-time applications. 
It is based on solving the HJB PDE to obtain the value function (Sec.~\ref{sec:HJB}). However, state-of-the-art HJB solvers, e.g., ENO/WENO \cite{osher1991high}, are grid-based and can suffer from the curse of dimensionality (CoD)~\cite{bellman57}, i.e., costs increase exponentially with dimension. 
For OC problems with a state-space dimension exceeding four, the CoD renders using grid-based HJB solvers infeasible.

We fuse the principles of the PMP and HJB methods to formulate an NN approach that is semi-global while mitigating CoD. 
In particular, we begin by parameterizing the value function with an NN, which circumvents CoD by approximating the solution to the HJB PDE in the underlying parameter space. Thus, our method is grid-free and suitable for high-dimensional problems.
Using the PMP, we express the control in feedback form.
We train the NN approximation of the value function by minimizing the expected cost on a distribution of initial states. 
As we minimize the cost function directly, our approach does not require generating solutions via an existing algorithm for training---i.e., it is not supervised.
Training the NN on a distribution of initial states ensures the controls' optimality on a large portion of the state-space; hence, our approach is semi-global.
As we demonstrate, the controls are robust to moderate perturbations or shocks to the system, such as wind interference (Sec.~\ref{subsec:shocks}). 
The controls are obtained in a feedback form via prior offline training, so the feedback form can be applied efficiently during deployment. 
Lastly, we improve the NN training by adding residual penalty terms derived from the HJB PDE, similar to \cite{ruthotto2020machine,lin2020apac,onken2020otflow}.

This paper extends a preliminary conference version of the approach~\cite{onken2020neural} with more extensive and thorough experiments. Specifically, we add experiments where agents swap positions with each other and one involving a nonlinear control-affine quadcopter with complicated dynamics. Additionally, we include experiments that investigate the sensitivity of NN hyperparameters, thoroughly compare the semi-global nature of the NN model against thousands of baseline solutions, demonstrate the efficient deployment timings of the NN, and test the influence of CoD on the NN.

Our formulation is applicable to OC problems for which the underlying Hamiltonian can be computed efficiently, e.g., affine controls with convex Lagrangians. 
Real-world applications that fall within this scope arise in centrally controlled multi-agent systems, which are the focus of this work. These also lead to challenging high-dimensional OC problems.
Indeed, for $n$ agents in a $q$-dimensional space, we obtain a $d{=}n\,{\cdot}\,q$-dimensional OC problem. Therefore, even moderate $n,q$ yield problems out of reach for traditional HJB solvers.

In our experiments, we solve a series of multi-agent OC problems whose state-space dimensions range from four ($n{=}2,q{=}2$) to 150 ($n{=}50,q{=}3$).
First, we solve a two-agent corridor problem with a smooth obstacle terrain (Sec.~\ref{sec:softcorridor}). Second, we investigate a two-agent problem where agents swap positions while avoiding hard obstacles and a 12-agent unobstructed version found in~\cite{mylvaganam2017differential} (Sec.~\ref{sec:swap}). Third, we experiment with a 50-agent swarm of three-dimensional agents obstructed by rectangular prisms inspired by~\cite{honig2018trajectory} (Sec.~\ref{sec:swarm}). Finally, we solve a 12-dimensional single-agent quadcopter problem with complicated dynamics from~\cite{lin2018splitting} (Sec.~\ref{sec:quadcopter}). 
Accompanying videos of our NN's solutions to these problems reside at \url{https://imgur.com/a/eWr6sUb}.

Using the corridor problem, we test our model's robustness to external shocks (random additive perturbations) that occur during deployment. We perform an example shock (Fig.~\ref{fig:shock}) and compare the NN's response against a baseline method (Sec.~\ref{sec:baseline}). Furthermore, we compare the solutions from the NN approach and the baseline on thousands of initial points (Fig.~\ref{fig:globalNN}). In this example, the NN reacts approximately optimally to moderate shocks (in terms of solution quality).
For large shocks, the NN learns a suboptimal control but still drives the agents towards the targets. However, the NN (trained offline) demonstrates quick online speed (Table~\ref{tab:stat}).

As one indicator that our approach effectively mitigates the CoD, we demonstrate empirically that increasing the state-space dimension does not lead to an exponential growth in computational costs. Specifically, we obtain approximately optimal controls by increasing the number of NN parameters linearly while keeping all other settings, including the batch size and number of optimization steps, fixed (Fig.~\ref{fig:cod}). We also show that we are able to solve a 150-dimensional problem in less than one hour on a single graphics processing unit (GPU).

\section{Related Work} \label{sec:related_works}

In recent years, many new numerical methods and machine learning approaches have been developed for solving high-dimensional OC problems. We discuss deterministic (Sec.~\ref{subsec:det}) and stochastic (Sec.~\ref{subsec:stoch}) settings separately because they differ considerably.
In Sec.~\ref{subsec:mapf}, we survey the state-of-the-art in the application domain that motivates our experiments.

\subsection{High-Dimensional Deterministic Optimal Control} \label{subsec:det}

A common difficulty in solving high-dimensional OC problems is the CoD. Exceptions are convex OC problems for which high-dimensional solvers can be devised via primal-dual methods and Hopf-Lax representation formulae~\cite{darbonosher16,lin2018splitting,Kirchner18,Kirchner2018APM,Chow2018AlgorithmFO,CHOW2019376,claudel2010lax1,claudel2010lax2}.

Kang and Wilcox \cite{kang2017mitigating} alleviate the CoD by introducing a sparse grid in the state-space and use the method of characteristics to solve boundary value problems over each sparse grid point. 
To approximate the feedback control at arbitrary points, they interpolate the solutions of the grid using high-order polynomials.
The authors solve up to six-dimensional control problems.
Nakamura-Zimmerer \textit{et al.}~\cite{nakamura2019adaptive} also attempt to alleviate CoD by learning a closed-form value function. First, trajectories are generated in a similar manner as in~\cite{kang2017mitigating}.
Using a supervised learning approach, the NN is then trained to match the generated trajectories. 
The trajectories (training data) are generated adaptively using information about the adjoint and by combining progressive batching with an efficient adaptive sampling technique. 

Bansal and Tomlin~\cite{bansal2021deepreach} solve high-dimensional reachability problems by combining the Hamilton-Jacobi-Isaacs (HJI) framework with the Deep Galerkin Method in~\cite{sirignano2018}. More precisely, they approximate the value function with an NN and minimize the empirical average of the HJI residual at randomly drawn space-time points.

Our work stems from the same framework as \cite{kunisch2020semiglobal}, which approximates the feedback control with an NN then optimizes the control cost on a distribution of initial states.
The authors also provide a theoretical analysis of OC solutions via NN approximations.
We extend the framework to finite horizon problems with non-quadratic costs and parameterize the value function instead of the feedback function. This extension enables penalization of the HJB conditions, which empirically improves numerical performance for solving high-dimensional mean-field games, mean-field control, and normalizing flows~\cite{ruthotto2020machine,lin2020apac,onken2020otflow}. We demonstrate similar advantages for OC problems considered in this work, which make similar use of NNs to parameterize the value function.

\subsection{High-Dimensional Stochastic Optimal Control} \label{subsec:stoch}

In the seminal works~\cite{E_2017,Han_2018}, the authors solve high-dimensional semilinear parabolic PDE problems by the method of (stochastic) characteristics. To overcome CoD, they approximate the gradient of the solution at different times by NNs and introduce a loss function that measures the deviation from the correct terminal condition in the characteristic equations. In particular, they solve high-dimensional \textit{stochastic} OC problems by solving the corresponding viscous HJB equation. This method recovers the gradient of the solution as a function of space and time and can be considered a global method. Nevertheless, loss functions employed in \cite{E_2017,Han_2018} consider only one initial point at a time, and the generalization depends on how well the generated random trajectories fill the space. 
The variance of the trajectories increases as time grows.
Finally, in the deterministic limit the method becomes local as there is no diffusion to enforce the trajectories to explore the whole space. Similar techniques are applied in \cite{nusken2020solving, moon2020generalized} based on different loss functions.

In \cite{han2016deep}, the authors solve stochastic OC problems by directly approximating controls and using the control objective as a loss function. As in \cite{E_2017,Han_2018}, the loss function considers a single initial point.

\subsection{Multi-Agent Path-Finding} \label{subsec:mapf}

Multi-Agent Path-Finding (MAPF)~\cite{stern2019multi,jing2019multiagent,zhao2018affine} methods are methods tailored for multi-agent control problems.
These methods tend to focus on collision avoidance rather than optimality. Among these are Conflict-Based Search (CBS) methods~\cite{sharon2015conflict,wagner2011m}, which are two-level algorithms. At the low level, optimal paths are found for individual agents, while at the high-level, a search is performed in a constraint tree whose nodes include constraints on time and location for a single agent. Decoupled optimization approaches~\cite{erdmann1987multiple,honig2018trajectory} first compute independent paths and then try to avoid collision afterwards. Other approaches phrase these as a constrained optimization problem~\cite{richards2002aircraft,blackmore2006optimal,patel2011trajectory,zhang2021opt}. Such methods are often combined with graph-based methods~\cite{standley2011complete}, sub-dimensional expansions~\cite{wagner2015subdimensional}, and CBS approaches~\cite{boyarski2015icbs,Cohen16Improved}.
Another approach phrases the MAPF problem as a differential game~\cite{mylvaganam2017differential}.
Provided certain assumptions, this differential game strategy guarantees that the agents reach their targets while avoiding collisions.
Machine learning approaches for multi-agent control have also been successfully applied in~\cite{riviere2020glas} where supervised learning is used to imitate non-machine-learning solutions generated by~\cite{honig2018trajectory}.
Our approach differs from these methods in that we do not have a data generation and fitting/imitation phases; instead, we directly solve for the control objective. Additionally, localization and interaction modeling techniques such as in \cite{shi2020neuralswarm} can be incorporated in our model in a straightforward manner.

\section{Mathematical Formulation} \label{sec:preliminaries}

We briefly discuss the general OC framework, derive the multi-agent control problems with collision avoidance used in the experiments, and review the theoretical foundations of the NN approach, primarily following~\cite[Chapters I-II]{flemingsoner06}.

\subsection{Optimal Control Formulation}

We consider deterministic, finite time-horizon OC problems. For a fixed time-horizon $[0,T]$, we have system dynamics
\begin{equation}\label{eq:charOC}
    \partial_s \bfz_{t,\bfx}(s) = f(s,\bfz_{t,\bfx}(s),\bfu_{t,\bfx}(s)), \quad \bfz_{t,\bfx}(t)= \bfx,
\end{equation}
for $t \leq  s \leq  T$.
Here, $\bfx \in \R^d$ is the initial state, and $t\in [0,T]$ is the initial time of the system. Next, $\bfz_{t,\bfx}(s) \in \R^d$ is the state of the system at time $s\in[t,T]$ with initial data $(t,\bfx)$, and $\bfu_{t,\bfx}(s) \in U \subset \R^a$ is the control applied at time $s$.
The function $f \colon [0,T]\times \R^d\times U \to \R^d$ models the evolution of the state $\bfz_{t,\bfx} \colon [t,T]\to \R^d$ in response to the control $\bfu_{t,\bfx} \colon [t,T]\to U$.

Next, we suppose that the control $\bfu_{t,\bfx} \colon [t,T]\to U$ and the trajectory $\bfz_{t,\bfx} \colon [t,T] \to \R^d$ satisfying~\eqref{eq:charOC} yield a cost
\begin{equation}\label{eq:Joc}
     \int_{t}^T L\big(s,\bfz_{t,\bfx}(s), \bfu_{t,\bfx}(s)\big) \, \du s  \,\, +  \,\, G\big(\bfz_{t,\bfx}(T)\big),
\end{equation}
where $L \colon [0,T]\times \R^d \times U \to \R$ is the \textit{running cost} or the \textit{Lagrangian}, and $G \colon \R^d \to \R$ is the \textit{terminal cost}. We assume that $f,L,G,U$ are sufficiently regular (see \cite[Sec. I.3, I.8-9]{flemingsoner06} for a list of assumptions). The goal of the OC problem is to find an optimal control $\bfu^*_{t,\bfx}$ that incurs the minimal cost, i.e,
\begin{equation}\label{eq:OC}
	\begin{split}
    \Phi(t,\bfx)=\inf_{\bfu_{t,\bfx}}  &\bigg\{  \int_{t}^T L\big(s,\bfz_{t,\bfx}(s), \bfu_{t,\bfx}(s)\big) \, \du s \\ &+  G\big(\bfz_{t,\bfx}(T)\big)  \bigg\} \,\,\, \text{s.t. }~\eqref{eq:charOC},
    \end{split}
\end{equation}
where $\Phi$ is called the \textit{value function}. A solution $\bfu^*_{t,\bfx}$ of \eqref{eq:OC} is called an \textit{optimal control}. Accordingly, the $\bfz^*_{t,\bfx}$ which corresponds to $\bfu^*_{t,\bfx}$ is called an \textit{optimal trajectory}.

We also define the \textit{Hamiltonian} of the system by
\begin{equation}\label{eq:H}
\begin{split}
    H(t,\bfz,\bfp)=&\sup_{\bfu \in U} \left\{ -\bfp \cdot f(t,\bfz,\bfu)-L(t,\bfz,\bfu) \right\}\\
    =&\sup_{\bfu \in U} \mathcal{H}(t,\bfz,\bfp,\bfu),
\end{split}
\end{equation}
where $\bfp \in \R^d$ is called the \textit{adjoint state}.
The Hamiltonian is a key ingredient in the Pontryagin Maximum Principle~\cite{pontryagin62} (Sec.~\ref{sec:PMP}) and also appears in the Hamilton-Jacobi-Bellman PDE~\cite{bellman57} (Sec.~\ref{sec:HJB}), which together form the foundation of our numerical solution approach.

\subsection{Collision-Avoiding Multi-Agent Control Problems} \label{sec:lagrangian}

While our NN approach is applicable to a broad range of OC problems, our numerical examples are motivated by centrally controlled multi-agent problems with collision avoidance.
Optimal decision-making for this class of problems is complicated due to the high-dimensionality of the control problem and the interactions between the agents.
These difficulties are exacerbated in the presence of random shocks and other forms of uncertainty.
Here, we describe the generic set up of these problems and refer to Section~\ref{sec:numerical_examples} for specific instances.

We seek to control a system of $n$ agents with initial states $x^{(1)},\dots,x^{(n)} \in \R^{q}$.
We denote the initial joint-state of the system by concatenating the agents' initial states, i.e.,
\begin{equation} \label{eq:x_agents}
	\bfx = \left( x^{(1)},x^{(2)},\dots,x^{(n)} \right) \in \R^d.
\end{equation}
Thus, the dimension of the joint-state of the system is $d=q\cdot n$.
Similarly, we denote the joint-state of the system at time $s$ by
\begin{equation}
    \bfz_{t,\bfx}(s) = \Big( z_{t,\bfx}^{(1)}(s), \, z_{t,\bfx}^{(2)}(s),\dots, \,z_{t,\bfx}^{(n)}(s) \Big),
\end{equation}
where, for a fixed $s \in [t,T]$, $z_{t,\bfx}^{(i)}(s) \in \R^{q}$ is the state of the $i$th agent.
Also,  we represent the control as
\begin{equation}
    \bfu_{t,\bfx}(s) = \Big( u_{t,\bfx}^{(1)}(s),  \, u_{t,\bfx}^{(2)}(s),\dots, \, u_{t,\bfx}^{(n)}(s) \Big).
\end{equation}
Hence, both the dimension of the state and the control space are proportional to the number of agents.

In the numerical experiments, the terminal costs depend on the distance between the agents' final positions and their given target states.
We denote the target joint-state of the system by the vector $\bfy \in \R^d$, obtained by concatenating the target states for all the agents as in~\eqref{eq:x_agents},  and consider the terminal cost
\begin{equation} \label{eq:G}
	G\big(\bfz_{t,\bfx}(T)\big) = \frac{\alpha_1}{2}\| \bfz_{t,\bfx}(T) - \bfy \|^2.
\end{equation}
The Lagrangian can be written as
\begin{equation}\label{eq:mult_ag_obj}
	L \big(s,\bfz, \bfu \big) = E \big(\bfu\big) \, + \, \alpha_2 Q \big( \bfz\big)
	 + \alpha_3 W \big( \bfz \big),
\end{equation}
where the scalar weighting parameters $\alpha_1$,$\alpha_2$,$\alpha_3$ are problem-specific and model the trade-off between the individual terms.

The first term in~\eqref{eq:mult_ag_obj}, $E\colon U \to \R$, is the \textit{energy term}, 
which is the total consumption cost comprised of individual ones
\begin{equation}
	E \big( \bfu_{t,\bfx} \big)  = \sum_{i=1}^n  E_i \Big( u_{t,\bfx}^{(i)} \Big).
\end{equation}
In our experiments, we use $E_i(u) {=} \frac12 \|u\|^2 \,{+}\, \kappa$ with a problem-dependent constant $\kappa\in\R$, which simplifies the Hamiltonian computation in~\eqref{eq:H}.
Unlike the other terms, this first term depends explicitly on the control.

The second term in~\eqref{eq:mult_ag_obj}, $Q\colon \R^d \to \R$, models obstacles by penalizing agents at certain spatial locations (e.g., a terrain function) and decouples into
\begin{equation}
	Q \big( \bfz_{t,\bfx} \big)  = \sum_{i=1}^n  Q_i \Big( z_{t,\bfx}^{(i)} \Big),
\end{equation}
where $Q_i \colon \R^q \to \R$ models the $i$th agent's spatial preferences.

The third term in~\eqref{eq:mult_ag_obj}, $W\colon \R^d \to \R$, models interactions among the individual agents. For example, this term can penalize proximity among agents to avoid collisions, i.e,
\begin{equation}
	W(\bfz_{t,\bfx}) = \sum_{j\neq i} w\Big(z_{t,\bfx}^{(i)},z_{t,\bfx}^{(j)}\Big)
\end{equation}
for function $w \colon \R^q \times \R^q \to \R$,
\begin{equation} \label{eq:Wij}
	w \Big(z^{(i)},z^{(j)} \Big) =
	\begin{cases}
     \exp \left( -  \frac{ \| z^{(i)} - z^{(j)} \|_2^2}{2r^2} \right) ,  &\left\| z^{(i)} {-} z^{(j)} \right\|_2 < 2r, \\
     0 , &\text{otherwise}.
	\end{cases}
\end{equation}
Here, $r>0$ is the radius of an agent's safety region or space bubble.
While not guaranteed, this $w$ term can in practice prevent the overlapping of the agents' space bubbles, thus avoiding collisions, when $\alpha_3$ is sufficiently large.
Our approach straightforwardly extends to non-symmetric interaction costs and heterogeneous agents.

We note that the presence of the terrain function $Q$ and the interaction potential $W$ render the objective function non-convex in $\bfz$.
However, in our experiments, the function is strongly convex (in fact, quadratic) in $\bfu$, which simplifies evaluations of the Hamiltonian~\eqref{eq:H} under certain assumptions on $f$.
Our framework can be directly applied to other choices of $G$, $E$, $Q$, and $W$ as long as $H$ can be computed efficiently.

\subsection{The Pontryagin Maximum Principle} \label{sec:PMP}

The Pontryagin Maximum Principle (PMP) provides a set of necessary first-order optimality conditions for the optimal control $\bfu^*_{t,\bfx}(\cdot)$ and trajectory $\bfz_{t,\bfx}^*(\cdot)$ originating from fixed initial data $(t,\bfx)$.
Since a new instance of the problem needs to be solved when the initial data change or the system's state deviates from the optimal curve, the PMP can be considered a \emph{local} solution method.

\begin{mytheorem}[Theorem I.6.3 \cite{flemingsoner06}]\label{thm:PMP} Assume that $(\bfz^*_{t,\bfx},\bfu^*_{t,\bfx})$ is a pair of an optimal trajectory and optimal control that solve~\eqref{eq:charOC}. Furthermore, assume that $\bfp_{t,\bfx} \colon [t,T]\to \R^d$ is the solution of the \emph{adjoint equation}
\begin{equation}\label{eq:adjoint}
    \begin{cases}
    \partial_s \bfp_{t,\bfx}(s)=\nabla_{\bfz} \mathcal{H} \big( s,\bfz^*_{t,\bfx}(s),\bfp_{t,\bfx}(s),\bfu^*_{t,\bfx}(s) \big),\\
    \bfp_{t,\bfx}(T)=\nabla_{\bfz} G \big(\bfz_{t,\bfx}^*(T) \big),
    \end{cases}
\end{equation}
for $t\leq s \leq T$. Then
\begin{equation}\label{eq:PMP}
\begin{split}
    \bfu^*_{t,\bfx}(s) \in \argmax_{\bfu \in U}~ \mathcal{H} \big( s,\bfz_{t,\bfx}^*(s),\bfp_{t,\bfx}(s),\bfu \big)
\end{split}
\end{equation}
for almost all $s\in [t,T]$.
\end{mytheorem}
 Theorem~\ref{thm:PMP} provides necessary conditions, and hence does not guarantee that the computed solutions are optimal.

  In general, finding $\bfu_{t,\bfx}^*, \bfz_{t,\bfx}^*, \bfp_{t,\bfx}$ that satisfy the PMP is difficult.
  Simultaneously solving the initial value problem~\eqref{eq:charOC} and the terminal value problem~\eqref{eq:adjoint} gives the system a particularly challenging forward-backward structure~\cite{kang2017mitigating,kang2019algorithms}.

  As we show below, the PMP can be applied more readily when the value function $\Phi$ is differentiable at $(t,\bfx)$.
  First, in this case, the conditions in Theorem~\ref{thm:PMP} are sufficient~\cite[Theorem 7.3.9]{cannarsa04}. \cite[Theorems 7.3.10, 7.4.20]{cannarsa04} provide a similar result with slightly weaker assumptions.
  Second, as we outline below, the solution of~\eqref{eq:adjoint} can be obtained from $\Phi$. The following is a standing assumption throughout the paper.%
  \begin{myassumption}\label{assum:feedback}
	 Assume that \eqref{eq:PMP} admits a unique continuous closed-form solution
	 \begin{equation}\label{eq:feedback}
	 	\bfu^*(s,\bfz,\bfp) = \argmax_{\bfu \in U}~ \mathcal{H} \big( s,\bfz,\bfp,\bfu \big)
	 \end{equation}
	 for every $s \in [t,T]$ and $\bfz, \bfp \in\R^d$.
  \end{myassumption}
  A closed-form solution for the optimal control exists in a wide variety of OC problems~\cite{sanchez2018real,tang2018learning,mo21fasttrack,bansal2021deepreach,kunisch2020semiglobal}.
  Importantly, the PMP can also be applied efficiently when~\eqref{eq:feedback} does not admit a closed-form solution but can be computed efficiently.

  The next theorem states that the value function $\Phi$ contains complete information about the optimal control and we can easily recover $\bfu_{t,\bfx}^*$ and $\bfp_{t,\bfx}$ from $\Phi$  when Assumption~\ref{assum:feedback} holds.%
 \begin{mytheorem}[Theorem I.6.2 \cite{flemingsoner06}]\label{thm:p=gradPhi}
 Assume that $\bfu^*_{t,\bfx}$ is a right-continuous optimal control  and $\Phi$ is differentiable at $(s,\bfz_{t,\bfx}^*(s))$ for $t \leq s <T$. Then
 \begin{equation}\label{eq:gradPhi}
 \bfp_{t,\bfx}(s)=\nabla_{\bfz} \Phi \big(s,\bfz_{t,\bfx}^*(s) \big)
 \end{equation}
 solves \eqref{eq:adjoint}. Also, \eqref{eq:PMP} simplifies to
 \begin{equation}\label{eq:opt_control_nec}
 \begin{split}
 	\bfu^*_{t,\bfx}(s)=\bfu^*\left(s,\bfz^*_{t,\bfx}(s),\nabla_{\bfz} \Phi \big(s,\bfz_{t,\bfx}^*(s) \big) \right)
 \end{split}
 \end{equation}
 for almost all $s\in [t,T]$.
 \end{mytheorem}
 Note that enforcing or computationally verifying the differentiability condition is virtually impossible.
 However, in many cases including our applications, the value function is expected to be differentiable almost everywhere.
 Even if $\Phi$ is not differentiable at $(t,\bfx)$ and the optimal control is not unique, $\bfp_{t,\bfx}$ can be recovered from the super differential $\partial^+_{\bfx} \Phi$~\cite[Theorem 7.3.10, 7.4.20]{cannarsa04}.

 Theorem~\ref{thm:p=gradPhi} characterizes optimal controls in a \textit{feedback form} \eqref{eq:opt_control_nec}.
 This means that no further optimization is necessary to find the optimal controls when the value function is known.
 Feedback form representations are valuable in real-world applications.
 If $\nabla \Phi$ can be quickly calculated, optimal controls are readily available at any point in space and time.
 As such, the feedback form avoids recomputing the optimal controls at new points in scenarios when sudden changes to the initial data or the system's state occur.

We can also use Assumption \ref{assum:feedback} to simplify the computation of the trajectories.
Using the \textit{envelope formula} \cite[Sec.~3.1, Theorem 1]{evans2010partial}, we see that
\begin{equation} \label{eq:envelope}
    \begin{split}
    \nabla_{\bfz} \mathcal{H} \big(t,\bfz,\bfp,\bfu^*(t,\bfz,\bfp) \big) &= \nabla_{\bfz} H(t,\bfz,\bfp)\\
    \nabla_{\bfp} \mathcal{H} \big( t,\bfz,\bfp,\bfu^*(t,\bfz,\bfp) \big) &= \nabla_{\bfp} H(t,\bfz,\bfp).
    \end{split}
\end{equation}
Hence, assuming the value function is known, we can express the optimal trajectory as
\begin{equation}\label{eq:HS}
    \begin{cases}
    \partial_s \bfz_{t,\bfx}^*(s)=-\nabla_{\bfp} H \big( s,\bfz_{t,\bfx}^*(s),\nabla_{\bfz} \Phi(s,\bfz_{t,\bfx}^*(s)) \big),\\
    \bfz_{t,\bfx}^*(t)=\bfx,
    \end{cases}
\end{equation}
for $s \in (t,T]$.
These dynamics do not explicitly contain the control, which reduces the problem to the state variables only.
Recall the optimal control can be computed via~\eqref{eq:opt_control_nec}.

The above derivation outlines how to obtain the optimal control and trajectory from the value function under some smoothness assumptions.
Once the value function $\Phi$ is known, this procedure can be applied for any initial data and also adapt the trajectory when the system is perturbed.
 Therefore, if $\Phi$ is computed, we have a \textit{global} solution method.
The key issue that we address next is the computation of $\Phi$.

\subsection{Hamilton-Jacobi-Bellman PDE} \label{sec:HJB}

In the previous section, we reviewed that the solution to the the OC problem~\eqref{eq:OC} for all initial data can be inferred from the value function $\Phi$.
In our approach, we exploit the fact that the value function satisfies the Hamilton-Jacobi-Bellman (HJB) PDE to help train our NN approximation of $\Phi$.

\begin{mytheorem}[Theorems I.5.1, I.6.1 \cite{flemingsoner06}]\label{thm:HJB}
Assume that the value function $\Phi \in C^1([0,T] \times \R^d)$. Then $\Phi$ satisfies the HJB equations (also called the \textit{dynamic programming} equations)
\begin{equation}\label{eq:HJB}
-\partial_s \Phi(s,\bfz) + H\big(s,\bfz,\nabla_{\bfz} \Phi(s,\bfz)\big)=0, \, ~ \Phi(T,\bfz) = G(\bfz)
\end{equation}
for all $(s,\bfz) \in [t,T) \times \R^d$.
Conversely, assume that $\Psi \in  C^1([0,T] \times \R^d)$ is a solution of \eqref{eq:HJB} and $\bfu^*_{t,\bfx}$ is such that
\begin{equation}\label{eq:opt_control_suf}
\begin{split}
    \bfu_{t,\bfx}^*(s) \in {\displaystyle \argmax_{\bfu \in U}} ~\mathcal{H} \big( s,\bfz_{t,\bfx}^*(s), \nabla_{\bfz} \Psi \big(s,\bfz_{t,\bfx}^*(s)\big),\bfu \big)
\end{split}
\end{equation}
for almost all $s\in [t,T]$. Then $\Psi=\Phi$, and $\bfu^*_{t,\bfx}$ is an optimal control.
\end{mytheorem}
The differentiability of $\Phi$ can be relaxed to differentiability almost everywhere in the framework of viscosity solutions~\cite[Chap. II]{flemingsoner06}.

The HJB PDE \eqref{eq:HJB} admits robust existence, uniqueness, and stability theory in the framework of viscosity solutions because \eqref{eq:HJB} is a convex constraint on $\Phi$~\cite{crandall1983viscosity}. Well-established numerical methods, such as ENO/WENO \cite{osher1991high}, benefit from convergence guarantees when solving \eqref{eq:HJB}.
However, these methods rely on grids and therefore are affected by the CoD.
Mitigating this limitation motivates our NN approach.

We note that the PMP is the \textit{method of characteristics} \cite[Sec.~3.2]{evans2010partial} for the HJB equation~\eqref{eq:HJB}.
To be precise, we can compute $\Phi$ along the trajectory $\bfz_{t,\bfx}$ from \eqref{eq:HS} by solving
    \begin{equation*}
    \begin{cases}
    \partial_s \phi_{t,\bfx}(s)=H \big(s,\bfz_{t,\bfx}^*(s),\bfp_{t,\bfx}(s)\big)\\
    \qquad \qquad -\bfp_{t,\bfx}(s)\cdot \nabla_{\bfp} H \big(s,\bfz_{t,\bfx}^*(s),\bfp_{t,\bfx}(s) \big)\\
    \phi_{t,\bfx}(T)=G(\bfz_{t,\bfx}^*(T)).
    \end{cases}
    \end{equation*}
We then have that $\phi_{t,\bfx}(s)=\Phi(s,\bfz^*_{t,\bfx}(s))$.

\section{Neural Network Approach} \label{sec:our_method}

Our approach seeks to minimize~\eqref{eq:Joc} subject to~\eqref{eq:charOC} for initial states sampled from a probability distribution in $\R^d$ whose density we denote by $\rho$.
Hence, it aims at solving the problem for all states along the optimal trajectories originating from those points.
Since the optimal trajectories given by the PMP are characteristics of the HJB equation, our method blends these two approaches.
To enable high-dimensional scalability, our method parameterizes the value function with an NN and computes the controls using~\eqref{eq:opt_control_nec} and~\eqref{eq:HS}.
The NN is trained in an unsupervised fashion by minimizing the sum of the expected cost that results from the trajectories and penalty terms that enforce the HJB equations along the trajectories and at the final-time.

\begin{table*}[t]
  \centering
  \caption{Variables and hyperparameters inherent to the problem itself (shared for NN and baseline) and the hyperparameters tuned for the NN approach. All $\alpha$ values are determined relative to the $\alpha$-less $E$ term in the problem definition. The $\beta$ hyperparameters are tuned relative to the $\alpha$ values.}
  \begin{tabular}{lcccccccccccc}
 	 & \multicolumn{5}{c}{Problem Definition} & \multicolumn{6}{c}{NN-specific Hyperparameters} & \\
 	 \cmidrule(lr){2-6} \cmidrule(lr){7-12}
  	 & $n$ & $d$ & $\alpha_1$ & $\alpha_2$ & $\alpha_3$ & $m$  & $\beta_1$ & $\beta_2$ & $\beta_3$ & $n_t$ & $n_t$ & NN\\
  	 & \# agents  & dim. & on $G$ & on $Q$ & on $W$ & width & on ${\rm HJt}$ & on ${\rm HJfin}$ & on ${\rm HJgrad}$ & training & validation & \# Params\\
  	\midrule
    Corridor & \hphantom{0}2  & \hphantom{00}4  & \hphantom{0}100 & $10^4$  & 300 & \hphantom{0}32 & 0.02 & 0.02 & 0.02 & 20 & 50 & \hphantom{00}1,311\\
    Swap 2~\cite{mylvaganam2017differential} & \hphantom{0}2  & \hphantom{00}4  & \hphantom{0}300  & $10^6$ & $10^5$ & \hphantom{0}16 & 1 & 1 & 3 & 20 & 50 & \hphantom{00,0}415 \\
    Swap 12~\cite{mylvaganam2017differential} & 12 & \hphantom{0}24 & \hphantom{0}300  & -   & $10^5$ & \hphantom{0}32 & 5 & 2 & 5 & 20 & 50 & \hphantom{00}2,196\\
    Swarm~\cite{honig2018trajectory} & 50 & 150 &  \hphantom{0}900  & $10^7$  & 25000 & 512 & 2 & 1 & 3 & 26 & 80 & 342,654  \\
    Quadcopter~\cite{lin2018splitting}     & \hphantom{0}1  & \hphantom{0}12 & 5000 & -   & - & 128 & 0.1 & 0 & 0 & 26 & 50 & \hphantom{0}18,576\\
    \bottomrule
  \end{tabular}
  \label{tab:hyperparameters}
\end{table*} 

\addtolength{\tabcolsep}{-1pt}  
\begin{table}
  \centering
  \caption{NN Statistics. Training times are approximate from running on a shared NVIDIA Quadro RTX 8000 GPU. Deployment times are from running on a single 2.6 GHz Intel(R) Xeon(R) CPU E5-4627 v3 core (Sec.~\ref{sec:computation}).}
  \begin{tabular}{lccccccc}
  & \multicolumn{4}{c}{Training} & \multicolumn{2}{c}{Deploy Time (ms)} \\
  \cmidrule(lr){2-5}  \cmidrule(lr){6-7}
   &  \multirow{2}{*}{\# Iters} & Batch &\multirow{2}{*}{\large$\frac{\text{ms}}{\text{Iter}}$} & Time & NN & Baseline \\
  &  & Size & & (min) & Step & Estimate\\
  	\midrule
    Corridor                                  &   1800 & 1024 & 320 & 10 & 4.4 & 2899 \\
    Swap 2~\cite{mylvaganam2017differential}  &   4000 & 1024 & 560 & 37 & 4.5 & 2571\\
    Swap 12~\cite{mylvaganam2017differential} &  4000 & 2048 & 260 & 17  & 3.6 & 1730\\
    Swarm~\cite{honig2018trajectory}          &  6000 & 1024 & 570 & 57  & 9.6 & 4026 \\
    Quadcopter~\cite{lin2018splitting}        &  6000 & 1024 & 720 & 72  & 5.2 & 3110 \\
    \bottomrule
  \end{tabular}
  \label{tab:stat}
\end{table}
\addtolength{\tabcolsep}{1pt} 

\subsection{Main Formulation}

We consider the semi-global version of the control problem and seek an approximately optimal control for initial states $\bfx \sim \rho$. 
We do so by approximating the value function $\Phi(\cdot)$ with an NN with parameters $\bfth$, which we denote by $\Phi(\,\cdot\,;\bfth)$. 
Thus, we can write the controls in feedback form and the loss in terms of the parameters.
In particular, we solve
\begin{equation} \label{eq:full_opt}
		\begin{split}
		\min_{\bfth} \; \E_{\bfx \sim \rho} \;\; \big\{ \ell_{\bfx}(T)   + G(\bfz_{0,\bfx}(T))
		+ \beta_1 c_{\rm HJt,\bfx}(T) \\ \quad + \beta_2 c_{\rm HJfin,\bfx} + \beta_3 c_{\rm HJgrad,\bfx} \big\},
		\end{split}
	\end{equation}
	subject to
	\begin{equation}\label{eq:ODE_cons}
		\renewcommand*{\arraystretch}{1.1} 
		\begin{split}
		\partial_s\begin{pmatrix}
		\bfz_{0,\bfx}(s)\\
		\ell_{\bfx}(s)\\
		c_{\rm HJt,\bfx}(s) \\
		\end{pmatrix} 
		=\begin{pmatrix}
		- \nabla_{\bfp} H(s,\bfz_{0,\bfx}(s),\nabla_{\bfz} \Phi(s, \bfz_{0,\bfx}(s); \bfth ))\\
		 L_{\bfx}(s)\\
		P_{\rm HJt,\bfx}(s)
		\end{pmatrix},
		\end{split}
	\end{equation}
where $\ell_{\bfx}(0)= c_{\rm HJt,\bfx}(0) = 0$ and $s \in [0,T]$. 
Here, $\ell$ accumulates the Lagrangian cost $L$ along the trajectories, the terms $c_{\rm HJt,\bfx}$, $c_{\rm HJfin,\bfx}$, $c_{\rm HJgrad,\bfx}, P_{\rm HJt,\bfx}$ penalize violations of the HJB, and the scalar penalty weights $\beta_1$,$\beta_2$,$\beta_3>0$ are assumed to be fixed. 
The remainder of this section defines and discusses these terms in more detail.

The term $\ell(T)$ corresponds to the time integral in~\eqref{eq:Joc}.
To compute $L$ at a given time, we use~\eqref{eq:H} and~\eqref{eq:envelope} and reformulate the Lagrangian in terms of the NN parameters $\bfth$ as
\begin{equation}
\begin{split}
    &L_{\bfx}(s)= -H \big( s,\bfz_{0,\bfx}(s),\nabla_{\bfz} \Phi(s,\bfz_{0,\bfx}(s) ; \bfth) \big)\\
    &+\nabla_{\bfz} \Phi \big( s, \bfz_{0,\bfx}(s) ; \bfth \big) \cdot \nabla_{\bfp} H \big( s,\bfz_{0,\bfx}(s),\nabla_{\bfz} \Phi(s, \bfz_{0,\bfx}(s) ; \bfth) \big).
\end{split}
\raisetag{27pt}
\end{equation}

We use HJB penalty terms $c_{\rm HJt,\bfx}$, $c_{\rm HJfin,\bfx}$, and $c_{\rm HJgrad,\bfx}$ derived from the HJB PDE~\eqref{eq:HJB} as follows:
\begin{equation} \label{eq:HJt}
\begin{split}
    P_{\rm HJt,\bfx}(s)=& \\
    \big| \partial_s \Phi(s,&\bfz_{0,\bfx}(s) ; \bfth) - H \big( s,\bfz_{0,\bfx}(s),\nabla_{\bfz} \Phi(s, \bfz_{0,\bfx}(s) ; \bfth) \big) \big|\\
    c_{\rm HJfin,\bfx} =& \left| \, \Phi(T,\bfz_{0,\bfx}(T); \bfth) - G(\bfz_{0,\bfx}(T)) \, \right| \\
    c_{\rm HJgrad,\bfx} =& | \, \nabla_{\bfz} \Phi(T, \bfz_{0,\bfx}(T); \bfth) - \nabla_{\bfz} G(\bfz_{0,\bfx}(T)) \, | .
\end{split}
\raisetag{12pt}
\end{equation}
The ${\rm HJ_{t}}$ penalizer arises from the first equation in~\eqref{eq:HJB}, whereas ${\rm HJ_{fin}}$ and ${\rm HJ_{grad}}$ are direct results of the final-time condition in~\eqref{eq:HJB} and its gradient, respectively.
Penalizers prove helpful in problems similar to~\eqref{eq:full_opt}~\cite{ruthotto2020machine,lin2020apac,onken2020otflow,finlay2020train}. 
These penalizers improve the training convergence (Sec.~\ref{sec:penalizers}) without altering the solution of \eqref{eq:full_opt}. 
The $P_{\rm HJt,\bfx}$ penalizer is accumulated along the trajectory similar to $L$.
The scalar terms $\beta_1$,$\beta_2$,$\beta_3$ weight the importance of each HJB penalizer and are hyperparameters of the NN (Sec.~\ref{sec:hyperparams}, Sec.~\ref{sec:beta_tuning}).

\subsection{Value Function Approximation}

	To enable scalability to high dimensions, we approximate the value function $\Phi$ with an NN. While our formulation supports a wide range of NNs, we design a specific model that enables efficient computation.

	We parameterize the value function as
	\begin{equation}
	\label{eq:NNArchitecture}
		\begin{split}
		\Phi(\bfs ; \bfth) = \bfw^\top N(\bfs;\bfth_N) + \frac{1}{2} \bfs^\top (\bfA^\top 
		\bfA)\bfs + \bfb^\top \bfs + c, \\ \quad \text{where} \quad\bfth = (\bfw, \bfth_N, \bfA, \bfb, c).
		\end{split}
	\end{equation}
	Here, $\bfs{=}(\bfx,t) \in \R^{d+1}$ are the inputs corresponding to space-time, $N(\bfs;\bfth_N) \colon \R^{d+1} \to \R^m$ is an NN,
	and $\bfth$ contains the trainable weights: $\bfw\,\,{\in}\,\,\R^m$, $\bfth_N\,\,{\in}\,\,\R^p$, $\bfA\,\,{\in}\,\,\R^{\gamma \times (d+1)}$, $\bfb\,\,{\in}\,\,\R^{d+1}$, $c\,{\in}\,\R$, where rank $\gamma{=}\min (10,d)$ limits the number of parameters in $\bfA^\top \bfA$. 
	Here, $\bfA$, $\bfb$, and $c$ model quadratic potentials, i.e., linear dynamics; $N$ models nonlinear dynamics.
	
	In our experiments, for $N$, we use a simple two-layer residual neural network (ResNet)~\cite{he2016deep}	
	\begin{equation} \label{eq:ResNet}
	    \begin{split}
	    \bfa_0 & = \sigma(\bfK_0 \bfs + \bfb_0) \\ 
	    N(\bfs;\bfth_N) & = \bfa_0 + \sigma(\bfK_1 \bfa_0 + \bfb_1),\\ 
	    \end{split}
	\end{equation}	
	for $\bfth_N{=}(\bfK_0,\bfK_1, \bfb_0,\bfb_1)$ where $\bfK_0 \in \R^{m \times (d+1)}$, $\bfK_1 \in \R^{m \times m}$, and $\bfb_0, \bfb_1 \in \R^{m}$.
	We use the element-wise nonlinearity $\sigma(\bfx)=\log(\exp(\bfx) + \exp(-\bfx))$, which is the antiderivative of the hyperbolic tangent, i.e., $\sigma'(\bfx)=\tanh(\bfx)$~\cite{ruthotto2020machine,onken2020otflow}.

\subsection{Numerical Implementation}

	We solve the ODE-constrained optimization problem~\eqref{eq:full_opt} using the discretize-then-optimize approach~\cite{gholami2019anode,onken2020do}, in which we define a discretization of the ODE, then optimize on that discretization.
	The forward pass of the model uses a Runge-Kutta (RK) 4 integrator with $n_t$ time steps to eliminate the constraints~\eqref{eq:ODE_cons}. 
	The objective function is then computed, and automatic differentiation~\cite{nocedal2006numerical} calculates the gradient of the objective function with respect to $\bfth$. We use the ADAM optimizer~\cite{kingma2014adam}, a stochastic subgradient method with momentum, to update the parameters $\bfth$. We iterate this process a selected maximum number of times. For the learning rate (step size) provided to ADAM, we follow a piece-wise constant decay schedule.
	For instance, in the experiment in Fig.~\ref{fig:penalizers}, we divide the learning rate by 10 every 800 iterations.
	
	To produce an NN that generalizes to the state-space, we must define initial points in a manner to promote model generalizability. 
	We assume the initial points are drawn from a distribution with density $\rho$. 
	We train the NN on one batch at a time of independent and identically distributed samples from the distribution. After training a number of iterations on that batch, we resample the distribution to define a new batch and train additional iterations on that batch. We repeat this process until we hit the maximum number of iterations. We commonly choose batches of 1024 or 2048 samples which are re-sampled every 25--100 iterations. We found no noticeable empirical difference in solution quality across those ranges. Through this process, the model uses few data points at each iteration, but does not overfit to a specific set of data points.

\begin{figure*}[t]
\centering
\begin{subfigure}{\linewidth}
  \centering
  \includegraphics[width=\linewidth]{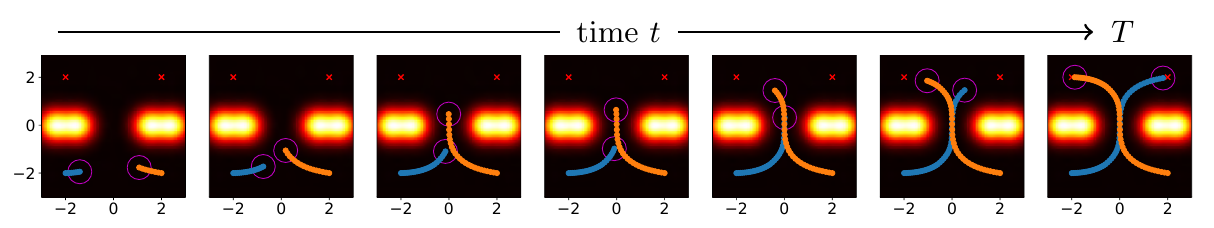}
  \subcaption{The baseline solution for initial state $\bfx_0$.}
  \label{fig:softcorridor_baseline}
\end{subfigure} \\
\begin{subfigure}{\linewidth}
  \centering
  \includegraphics[width=\linewidth]{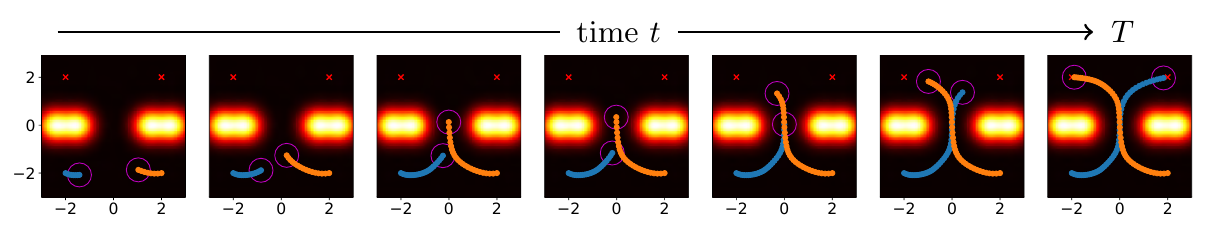}
  \subcaption{For initial state $\bfx_0$ (depicted), the NN learns a similar solution as the baseline. The NN approach is solved for multiple initial states.}
  \label{fig:softcorridor_nn}
\end{subfigure}
\caption{Solutions for the two-agent corridor problem where two agents (orange and blue) pass in between two smooth hills. Taking the terrain into account, the agents seek shortest paths from the initial joint-state $\bfx_0$ to target $\bfy$ (marked with red crosses) while avoiding collision with each other's space bubble (indicated by circles with radius $r$).}
\label{fig:softcorridor}
\end{figure*}

\subsection{Hyperparameters} \label{sec:hyperparams}

	In contrast to the model parameters $\bfth$ learned from the data, NN hyperparameters are values selected \textit{a priori} to training. These include the number of time steps $n_t$, the ResNet width $m$, ResNet depth (the number of layers, tuned to equal 2), and the multipliers $\beta_1,\beta_2,\beta_3$. Additionally, each OC problem has defined $\alpha_1,\alpha_2,\alpha_3$, which both the baseline and NN use; changing these values alters the problem (Table~\ref{tab:hyperparameters}).
	For reproducibility, we include all hyperparameters and settings with a publicly available Python implementation at \url{https://github.com/donken/NeuralOC}. Training on a single NVIDIA Quadro RTX 8000 GPU requires between 10 and 72 minutes for the considered OC problems (Table~\ref{tab:stat}).

\section{Numerical Experiments} \label{sec:numerical_examples}

We solve and analyze five OC problems and compare the NN against a baseline method described in Sec.~\ref{sec:baseline}.
In Sec.~\ref{sec:softcorridor} to \ref{sec:swarm}, we present four centrally controlled multi-agent examples with dimensionality ranging from 4 to 150. In Sec.~\ref{sec:quadcopter}, we consider a quadcopter experiment to demonstrate the NN's ability to solve problems with complicated dynamics.

\begin{table}[t]
  \centering
  \caption{Comparison of solution values for the two-agent corridor problem and single instance $\bfx_0$ shown in Fig.~\ref{fig:softcorridor}.}
	\begin{tabular}{lccc}
	   	\toprule
	 	 Method & $\ell + G$ & $\ell$ & $G$ \\
	   	\midrule
	     Baseline & 61.33 & 61.02 & 0.31 \\ 
	     NN       & 62.19 & 61.98 & 0.21 \\ 
	 	\bottomrule
	\end{tabular} 
	\label{tab:softcorridor}
\end{table}

\begin{figure}
\centering
	\begin{subfigure}{\linewidth}
	  \centering
	  \includegraphics[width=\linewidth]{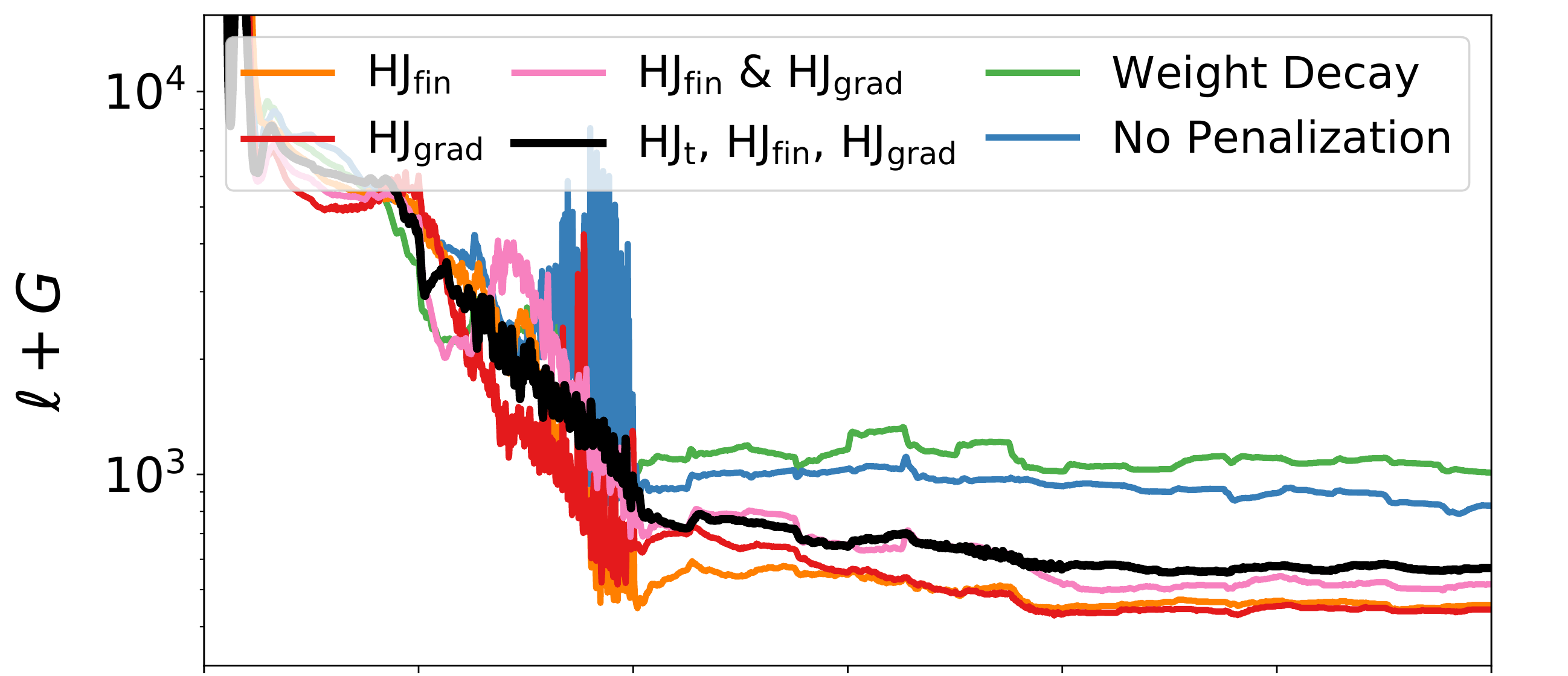} 
	\end{subfigure} \\
	\begin{subfigure}{\linewidth}
	  \centering 
	  \includegraphics[width=\linewidth]{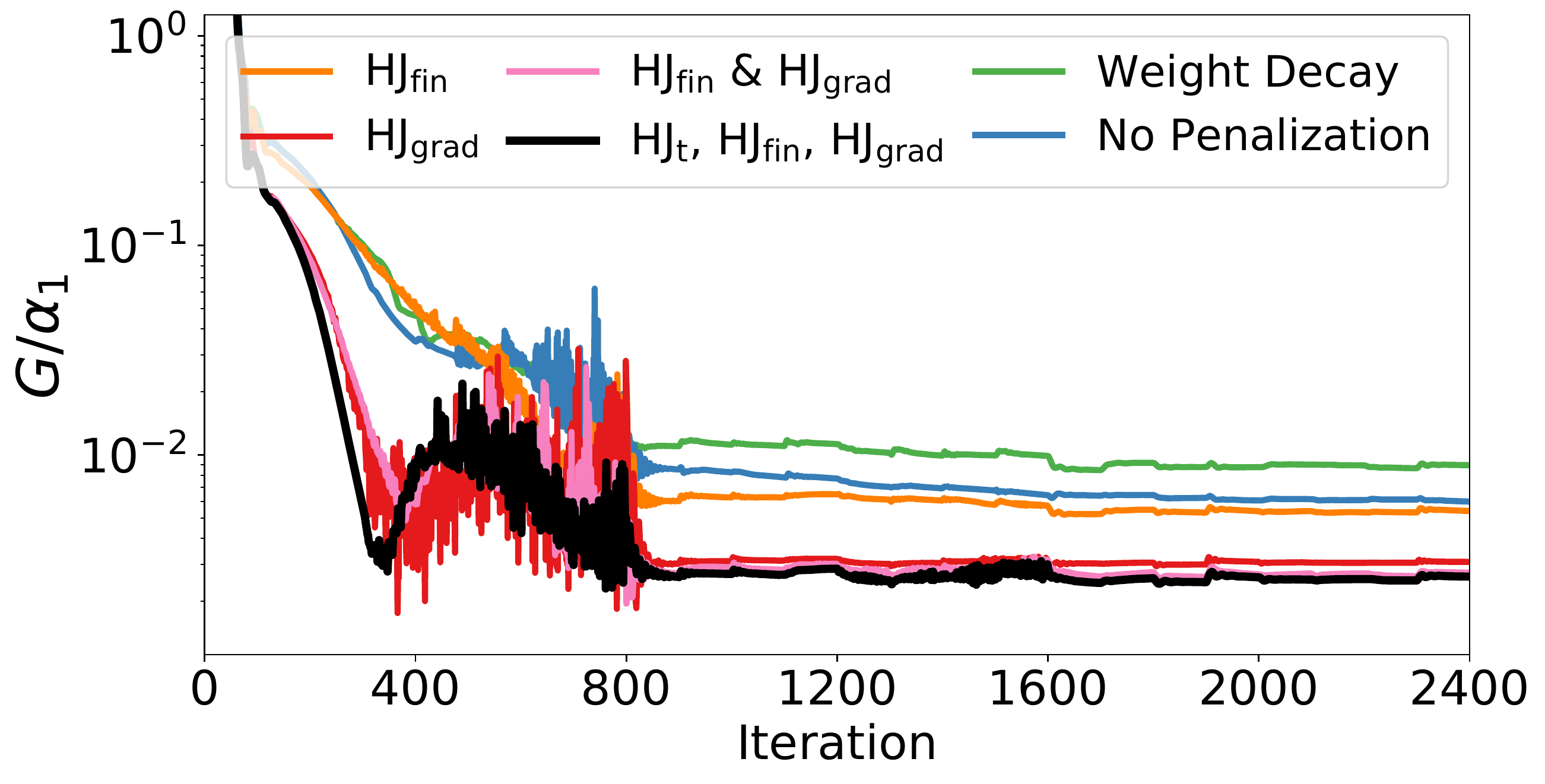}
	\end{subfigure}
	\caption{For the corridor problem (Sec.~\ref{sec:softcorridor}), we train the same model architecture six times using different combinations of the penalty terms. Using all three HJB penalizers leads to quick convergence and a low $G$ value. Each curve is the average of three training instances.}
	\label{fig:penalizers}
\end{figure}

\begin{figure*}[t]
\centering
\begin{subfigure}{0.22\linewidth}
  \centering \captionsetup{width=0.92\linewidth}%
   \includegraphics[width=\linewidth]{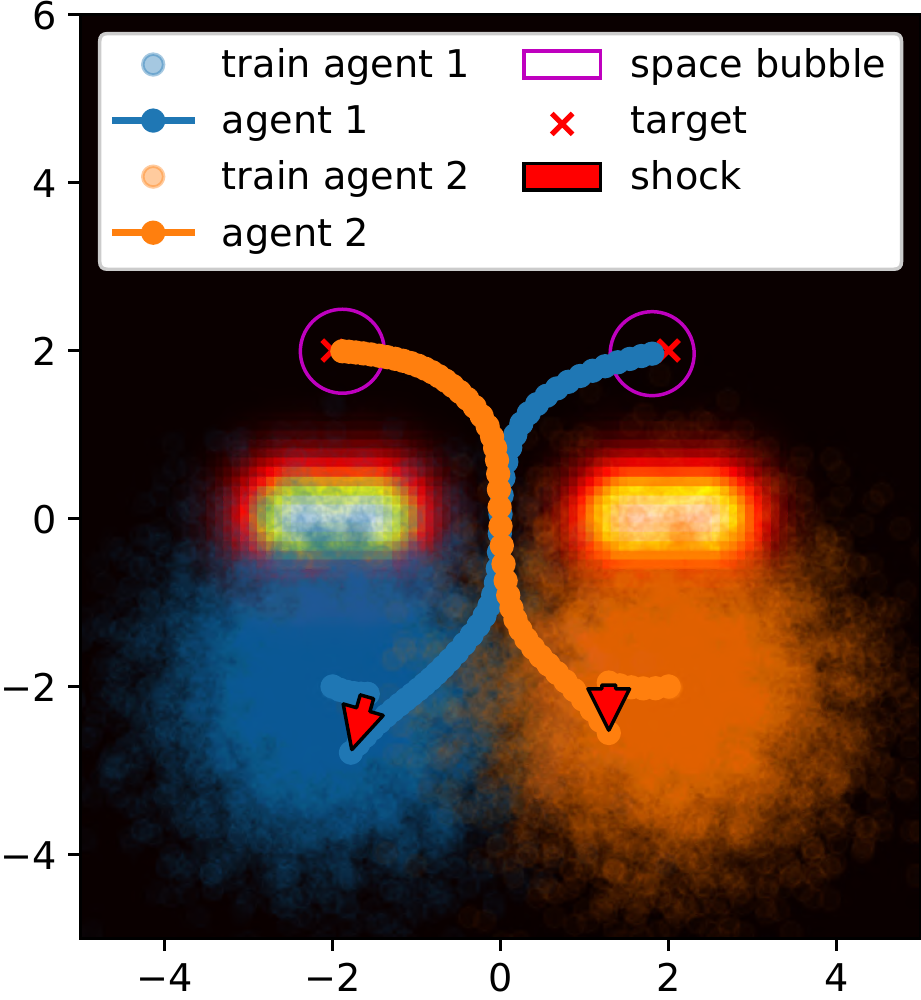}
  \subcaption{Minor shock $\| \bfxi \|=0.94$ within the training space.}
  \label{fig:shock_a}
\end{subfigure}
\begin{subfigure}{0.22\linewidth}
  \centering \captionsetup{width=0.92\linewidth}%
   \includegraphics[width=\linewidth]{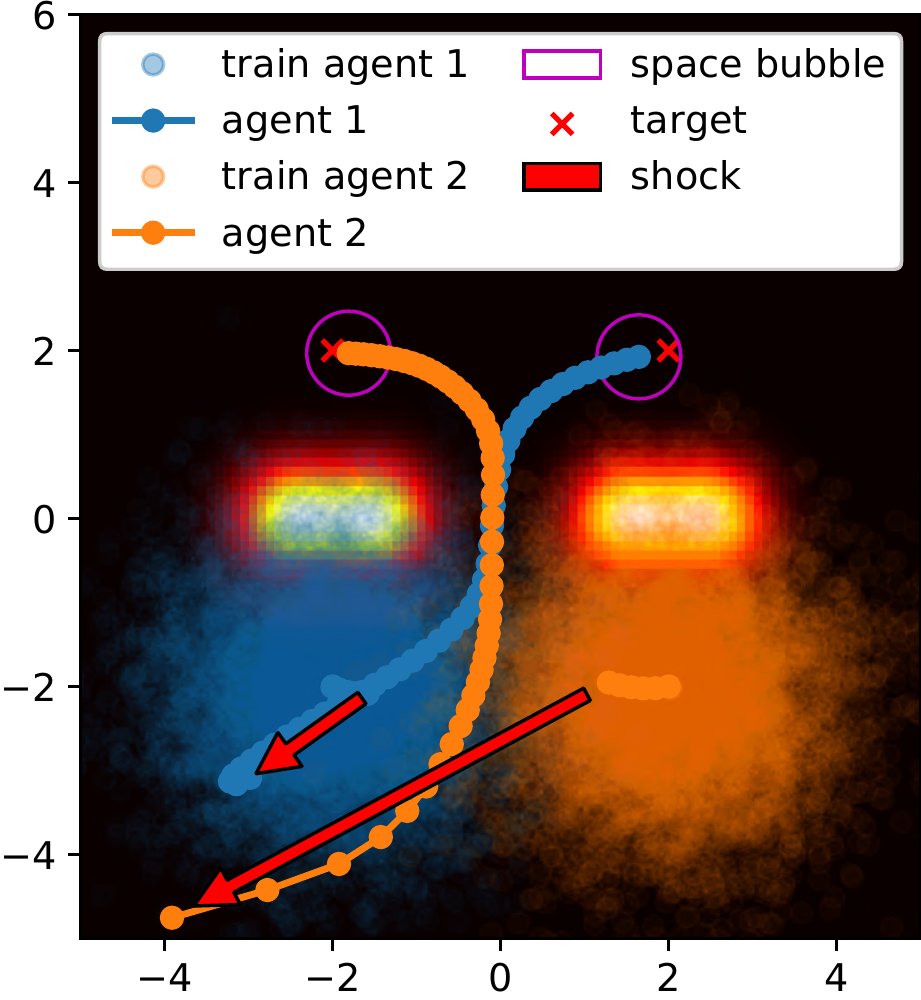}
  \subcaption{Major shock $\| \bfxi \|=6.2$ outside of the training space.}
  \label{fig:shock_b}
  \end{subfigure}
\begin{subfigure}{0.22\linewidth}
	  \captionsetup{width=0.92\linewidth}%
  \centering 
  \vspace{7pt}
  $\hat{\bfx} = \bfz_{0,\bfx_0}(0.1) + \bfxi$, $s= [0.1,1]$\\
  \vspace{3pt}
  \includegraphics[width=\linewidth]{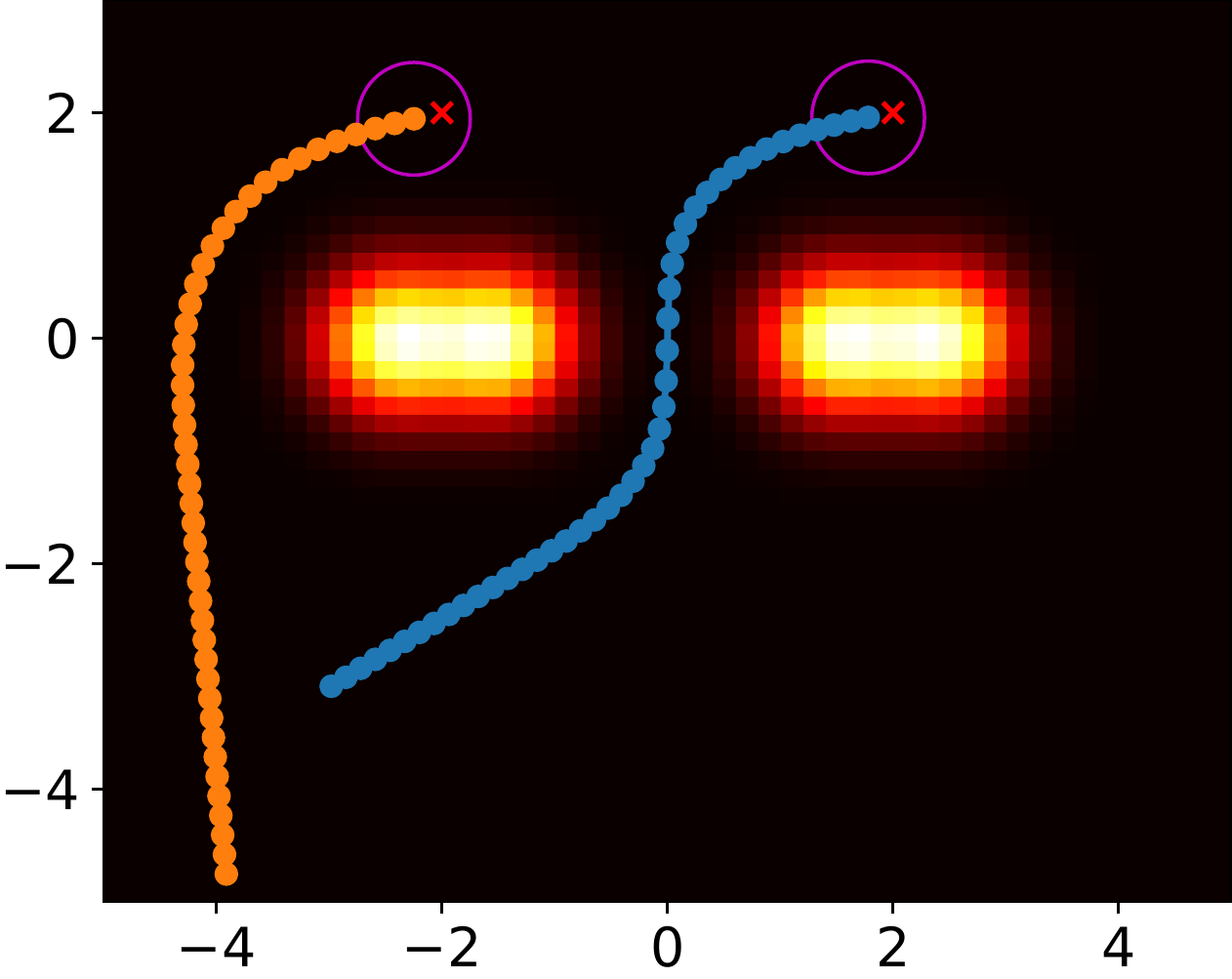}
  \subcaption{Baseline solution for (0.1,$\hat{\bfx}$) after major shock.} 
  \label{fig:shock_c}
  \end{subfigure}  
\begin{subfigure}{0.28\linewidth}
	\captionsetup{width=0.92\linewidth}%
	\vspace{20pt}
  	\centering
	{\small 
	\begin{tabular}{lccc}
	   	\toprule
	 	   & $\ell + G$ & $\ell$ & $G$ \\
	 	 \midrule
	 	  \multicolumn{4}{c}{\textbf{following shock $\| \bfxi \|=0.94$}} \\
	     Baseline & \hphantom{1}59.79 & \hphantom{1}59.46 & 0.33  \\
	     NN       & \hphantom{1}60.54 &  \hphantom{1}60.34  &  0.20  \\
	     \midrule 
	 	  \multicolumn{4}{c}{\textbf{following shock $\| \bfxi \|=6.2$}} \\
	     Baseline & \hphantom{1}71.77 & \hphantom{1}71.22 & 0.55 \\
	     NN       & 151.67 &  150.63 &  1.03\\     
	 	\bottomrule
	\end{tabular}
	} 
	\vspace{16pt}
\subcaption{Solution comparison of the methods on $s{=}[0.1,1]$ following a shock.}
\label{fig:shock_d}
\end{subfigure}
\caption{The NN handles a shock $\bfxi$ at time $s{=}0.1$ (depicted with red arrows) along the trajectory for the depicted corridor problem (Sec.~\ref{sec:softcorridor}). The initial states used during training are depicted as blue and orange point clouds. It can be seen that the major shock causes the system to leave the state-space used during training.}
\label{fig:shock}
\end{figure*}
\begin{figure}[t]
\centering
	\begin{subfigure}{0.46\linewidth}
	  \centering \captionsetup{width=0.9\linewidth}%
	  \includegraphics[width=\linewidth]{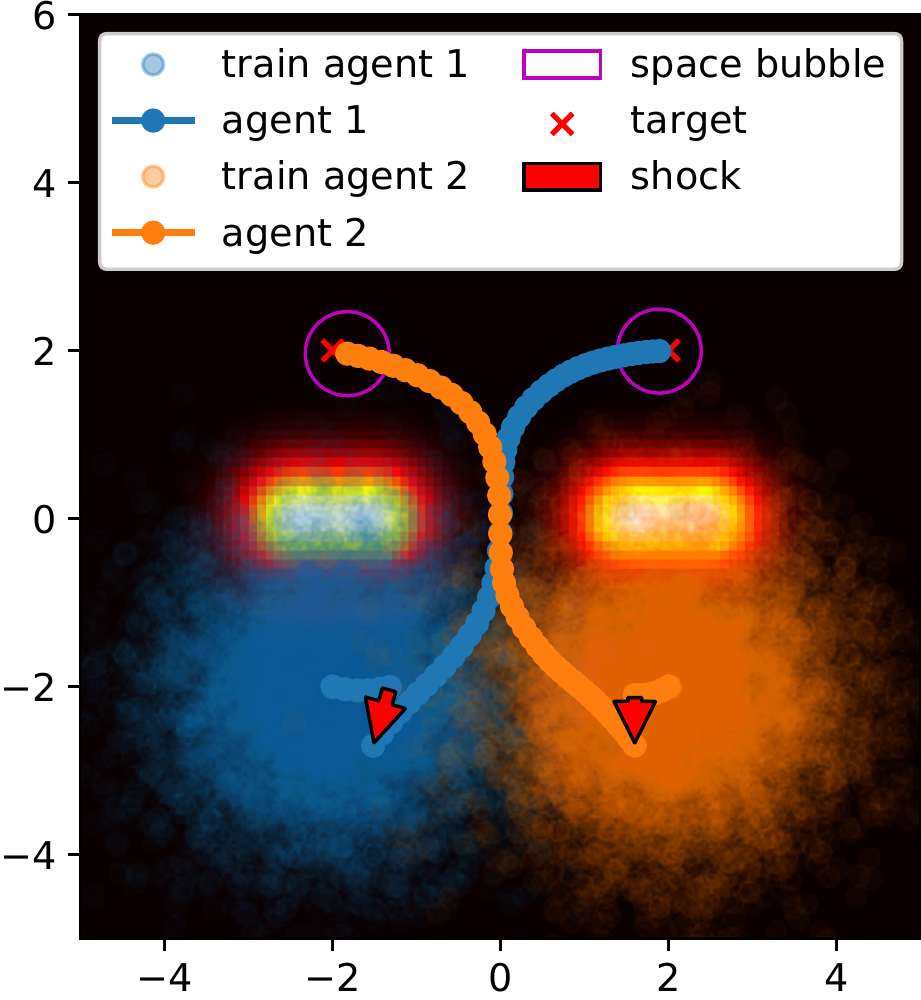}
	  \subcaption{Minor shock $\|\bfxi\|=0.94$ within the training space.}
	\end{subfigure}
	\begin{subfigure}{0.46\linewidth}
	  \centering \captionsetup{width=0.9\linewidth}%
	  \includegraphics[width=\linewidth]{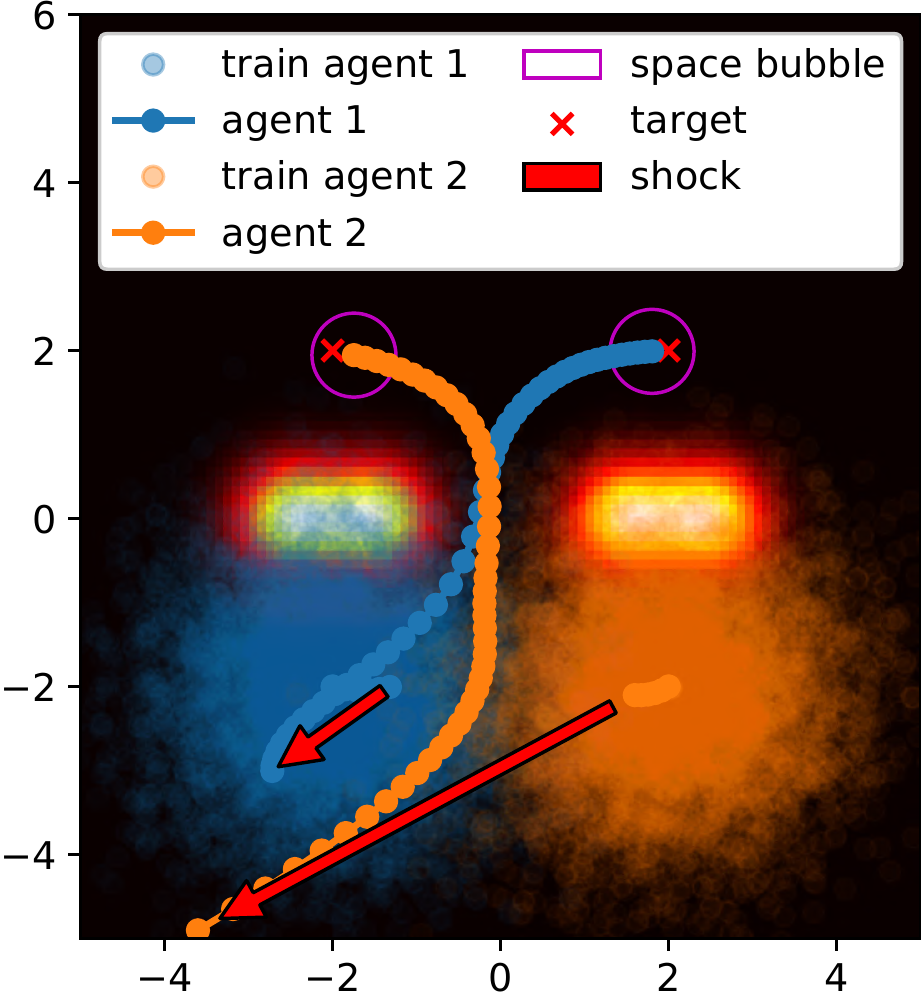}
	  \subcaption{Major shock $\| \bfxi \|=6.2$ outside of the training space.} 
	\end{subfigure}
\caption{We solve the corridor problem with an NN trained without HJB penalizers or weight decay. Comparable to Fig.~\ref{fig:shock}, we see that the penalizers do not alter the solution.}
\label{fig:wt_decay}
\end{figure}

\subsection{Baseline: Optimization for a Single Initial State} \label{sec:baseline}

For comparison with the NN approach, we provide a local solution method that solves the OC problem for a fixed initial state $\bfx_0$. We consider the baseline approach's solution as the ground truth optimal solution and compute the \textit{suboptimality} of the semi-global NN approach's solution evaluated for the initial state relative to the baseline's solution.

For the baseline, we obtain an optimization problem by applying forward Euler to the state equation and a midpoint rule to the integrals, i.e., the ``direct transcription method''~\cite{betts98survey},
\begin{equation} \label{eq:baseline}
	\begin{split}
    \min_{ \{\bfu^{(k)}\} } \quad &
    \, G \left(\bfz_{n_t} \right) + h \sum_{k = 0}^{n_t-1}  L\left(s_k,\bfz_k,\bfu_k\right)\\
    \text{s.t.}  \quad  &\bfz_{k+1} = \bfz_k + h \, f\big(s_k, \bfz_k ,\bfu_k \big), \, \, \bfz_0 = \bfx_0,
    \end{split}
\end{equation}
where $h{=}T/n_t$. Here, we use $\bfz_k$ to denote $\bfz_{0,\bfx_0}(s_k)$, where time point $s_k = hk$.
We use $T{=}1$ and $n_t{=}50$ and solve~\eqref{eq:baseline} using ADAM with initialization of the controls set as straight paths from $\bfx_0$ to $\bfy$ with small added Gaussian noise. 

We arrived at these training decisions empirically. First, 
when solving~\eqref{eq:baseline} in our experiments, ADAM finds slightly more optimal solutions ($1{-}2\%$ more optimal) in practice than L-BFGS. Second, the initialization of the controls substantially influences the solution. As a particular example, the baseline solution depicted in Fig.~\ref{fig:shock_c} learns to send agent 2 around the left side of the left obstacle, resulting in the lowest value of the objective function. If initialized with controls that pass through the right of that obstacle or through the corridor, the baseline struggles to learn this optimal trajectory. 
As a response, we initialize the controls uniformly that lead to a straight path from $\bfx_0$ to $\bfy$. Third, we add random Gaussian noise to the initialization because doing so empirically helps avoid local minima and overall achieves better results.

\subsection{Two-Agent Corridor Example} \label{sec:softcorridor}

We design a $d{=}4$-dimensional problem in which two agents attempt to reach fixed targets on the other side of two hills (Fig.~\ref{fig:softcorridor}). We design the hills in such a manner that one agent must pass through the corridor between the two hills while the other agent waits. For this example, the hills use a smooth terrain, and we assess the resilience of the control to shocks.

\subsubsection{Set-up}\label{subsec:corr_setup}

Suppose two homogeneous agents with safety radius $r{=}0.5$ start at $x^{(1)}{=}[-2, -2]^{\top}$ and $x^{(2)}{=}[2,-2]^{\top}$ with respective targets $y^{(1)}{=}[2, 2]^{\top}$ and $y^{(2)}{=}[-2, 2]^{\top}$. Thus, the initial and target joint-states are
$\bfx_{0}{=}[-2,-2,2,-2]^{\top}$ and $\bfy{=}[2,2,-2,2]^{\top}$.
We sample from $\rho$, which is a Gaussian centered at $\bfx_0$ with an identity covariance. These sampled initial positions form the training set $\bfX$.

The running costs depend on the spatio-temporal cost function $Q_i$. 
Throughout, obstacles are defined using the Gaussian density function with mean $\bfmu \in \R^q$ and covariance $\bfSigma \in \R^{q \times q}$
\begin{equation*}
    \eta(z^{(i)} \, ; \, \bfmu, \bfSigma) = \frac{\exp \left( -\hf (z^{(i)} - \bfmu)\bfSigma^{-1} (z^{(i)} - \bfmu) \right)}{\sqrt{(2\pi)^d \det{\bfSigma}}}.
\end{equation*}

In this experiment, we define obstacles as
\begin{equation*}
	\begin{split}
	Q_i \Big( z^{(i)} \Big) = \eta \left( z^{(i)} \, ; \, \begin{bmatrix} -2.5 \\ 0 \end{bmatrix}, 0.2\bfI \right) + \eta \left( z^{(i)} \, ; \,  \begin{bmatrix} 2.5 \\ 0 \end{bmatrix}, 0.2\bfI \right) \\ + \eta \left(z^{(i)} \, ; \, \begin{bmatrix} -1.5 \\ 0 \end{bmatrix}, 0.2\bfI \right) + \eta \left(z^{(i)} \, ; \, \begin{bmatrix} 1.5 \\ 0.0 \end{bmatrix}, 0.2\bfI \right).
	\end{split}
\end{equation*}
The energy terms are given by
\begin{equation} \label{eq:E_i}
	E_i \Big( u^{(i)} \Big) = \hf \big\| u^{(i)} \big\|^2,
\end{equation}
and the dynamics are given by $f(s,\bfz,\bfu) = \bfu$.

We compute the Hamiltonian~\eqref{eq:H} as
\begin{equation}
\begin{split}
	&H(s,\bfz,\bfp) = \sup_{\bfu \in U} \Big\{ -\bfp^{\top} \bfu - L\big(s, \bfz, \bfu \big) \Big\}  \\
    &= \sup_{\bfu \in U} \Big\{ -\bfp^{\top} \bfu -  E \big(\bfu \big) - \alpha_2 Q \big( \bfz \big) - \alpha_3 W \big( \bfz \big)  \Big\}.
\end{split}
\raisetag{45pt}
\end{equation}
We then can solve for the first-order necessary condition
\begin{equation} \label{eq:p_equals_u}
\begin{split}
	0 &= -\bfp -  \nabla_{\bfu} E \big( \bfu \big) \\
	\Rightarrow \quad \bfp &= -\nabla_{\bfu} \left( \sum_{i=1}^n \hf \big\| u^{(i)} \big\|^2 \right) = - \bfu
\end{split}
\end{equation}

Using the closed-form solution for the controls~\eqref{eq:p_equals_u}, we rewrite the Hamiltonian as
\begin{equation}
\begin{split}
	H(s,\bfz,\bfp)
	&= \| \bfp \|^2 -  \hf \| \bfp \|^2 - \alpha_2 Q \big( \bfz \big)  - \alpha_3 W \big( \bfz \big)  \\
	&= \hf \| \bfp \|^2 - \alpha_2 Q \big( \bfz \big) - \alpha_3 W \big( \bfz \big),  \\
\end{split}
\end{equation}
where the characteristics are given by
\begin{equation} \label{eq:corridor_f}
 	\partial_s \bfz_{t,\bfx}(s) = - \nabla_{\bfp} H\big(s,\bfz_{t,\bfx}(s),\bfp_{t,\bfx}(s)\big) =  -\bfp_{t,\bfx}(s).
\end{equation}

\begin{figure}
\centering	
	\begin{subfigure}{\columnwidth}
 	\centering 
 		\includegraphics[width=0.77\columnwidth]{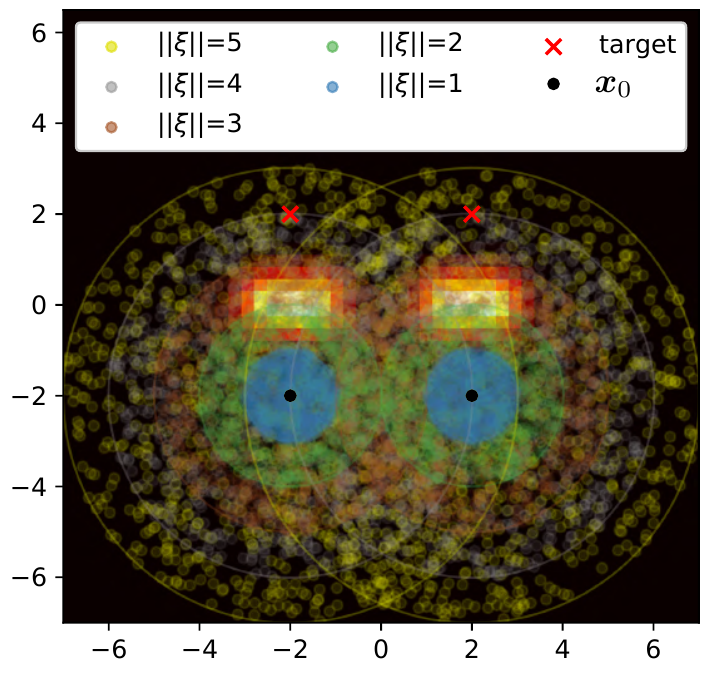}
	  \subcaption{The initial points $\bfx_0+\bfxi$ for the corridor problem sampled from the hyperspheres of radius $\|\bfxi\|$. }
	\end{subfigure}\\	
	\begin{subfigure}{\columnwidth}
	  \centering
	  \includegraphics[width=0.9\columnwidth]{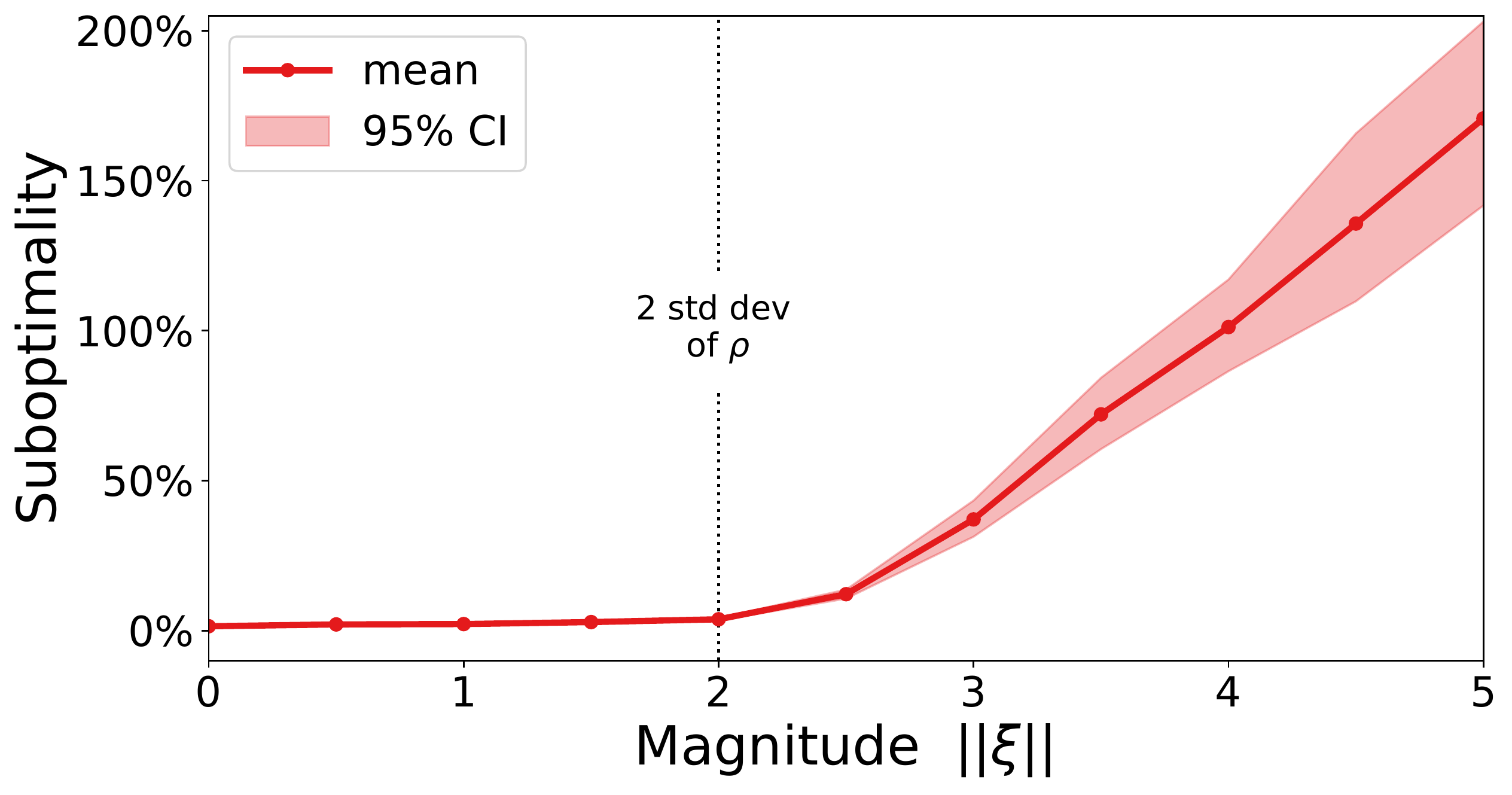}
	  \subcaption{The mean suboptimality of the NN's solution $\ell +G$, where the baseline solution for each initial point is considered optimal. }
	\end{subfigure} \\
	\begin{subfigure}{\columnwidth}
	  \centering 
	  \includegraphics[width=0.9\columnwidth]{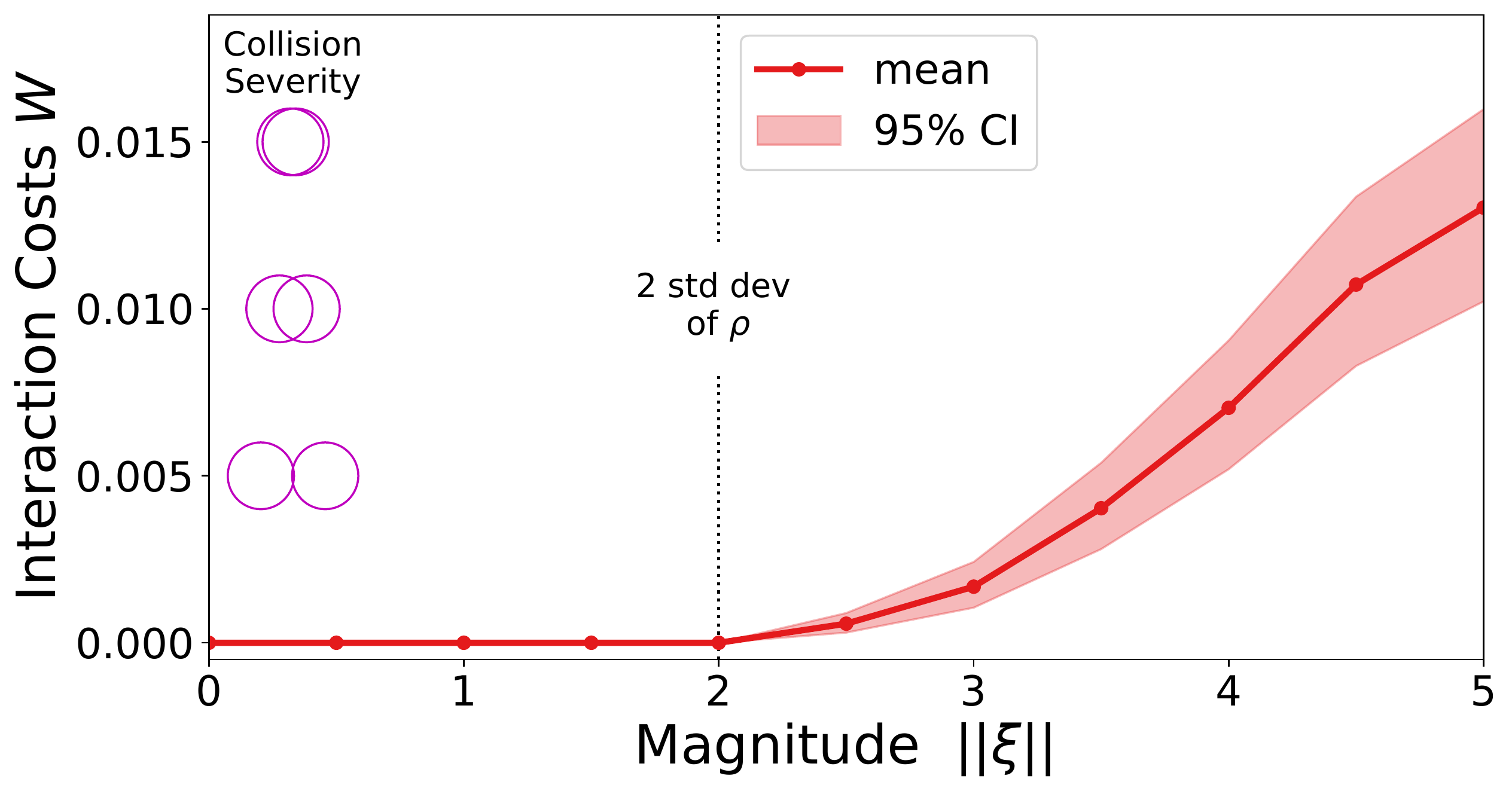}
	 \subcaption{NN interaction costs with comparable example collision severity of two circular agents.}
	\end{subfigure}\\
	\begin{subfigure}{\columnwidth}
	  \centering 
	  \includegraphics[width=0.9\columnwidth]{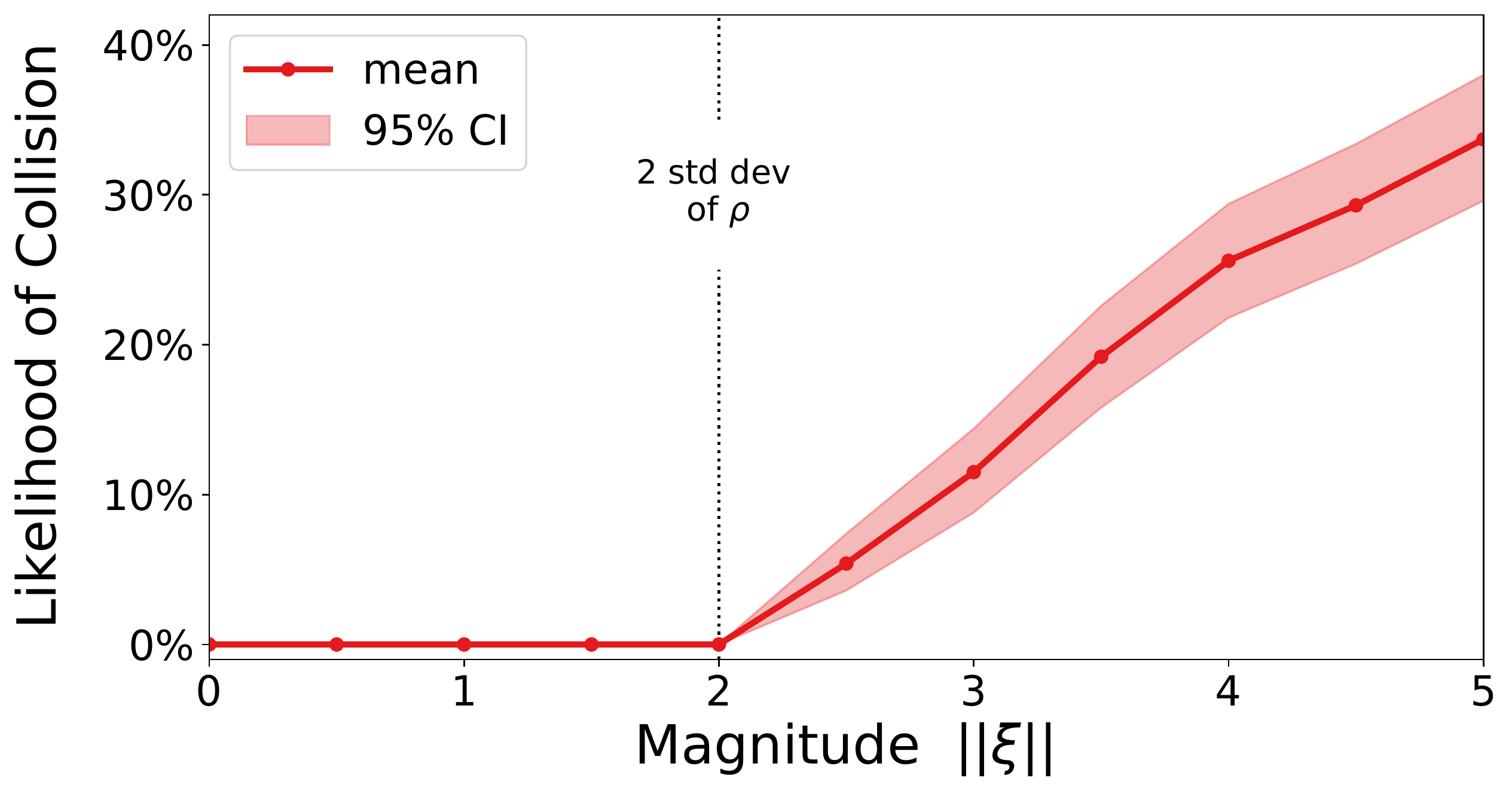}
	  \subcaption{For initial points at each magnitude, we present the percentage of those resulting in a collision of any severity when run with the NN.}
	\end{subfigure}		
	\caption{We compare one NN model with 10,001 baseline models for 1000 initial points $(0,\bfx_0+\bfxi)$ at each magnitude $\|\bfxi\|$. Confidence intervals are computed via 10,000 sub-samplings of size 500 from each set of 1000 points.}
	\label{fig:globalNN}
\end{figure}

\subsubsection{Results}\label{subsec:corridor_results}

	The baseline and the NN learn to wait for one agent to pass through the corridor first, followed by the second agent (Fig.~\ref{fig:softcorridor}). 
	The NN performs marginally worse in $L$ values (Table~\ref{tab:softcorridor}), which can be seen in the early stages of the trajectories of agent 1 (Fig.~\ref{fig:softcorridor_nn}).
	The NN achieves a slightly better $G$ value than the baseline. 
	Although we solve the NN by optimizing the expectation value of a set of points in the region, the NN achieves a near-optimal solution for $\bfx_0$.

\subsubsection{Effect of the HJB Penalizers} \label{sec:penalizers}

We experimentally assess the effectiveness of the penalizers $c_{\rm HJt}$, $c_{\rm HJfin}$, $c_{\rm HJgrad}$ in~\eqref{eq:full_opt}. To this end, we define six models (various combinations of the three HJB penalizers and one model with weight decay) and train three instances of each on the corridor problem. Using the HJB penalizers results in a quicker model convergence on a hold-out validation set (Fig.~\ref{fig:penalizers}).

$\mathbf{HJ_t}\colon$
We enforce the PDE~\eqref{eq:HJB} describing the time derivative of $\Phi$ along the trajectories. Including this penalizer improves regularity and reduces the necessary number of time steps when solving the dynamics~\cite{yang2019,ruthotto2020machine,lin2020apac,onken2020otflow}.

$\mathbf{HJ_{\rm fin}}\colon$ We enforce the final-time condition of the PDE~\eqref{eq:HJB}. The inclusion of this penalizer helps the network achieve the target~\cite{ruthotto2020machine}. Experimentally, using ${\rm HJ_{fin}}$ correlates with a slightly lower $G$ value (Fig.~\ref{fig:penalizers}). 

$\mathbf{HJ_{\rm grad}}\colon$ We enforce the transversality condition $\nabla_{\bfz} \Phi(T,\bfz(T)){=}\nabla_{\bfz} G(\bfz(T))~\forall \bfz$, a consequence of the final-time HJB condition~\eqref{eq:HJB}. Numerically, all conditions are enforced on a finite sample set. Therefore, higher-order regularization may help the generalization; i.e., achieving a better match of $\Phi(T,\cdot)$ and $G$ for samples not used during training (the hold-out validation set). We observe the latter experimentally;
${\rm HJ_{grad}}$ impacts $G$ more than ${\rm HJ_{fin}}$ (Fig.~\ref{fig:penalizers}). Nakamura-Zimmerer \textit{et al.}~\cite{nakamura2019adaptive} similarly enforce $\nabla \Phi$ values.

\begin{figure*}
  \centering 
	\includegraphics[width=\linewidth]{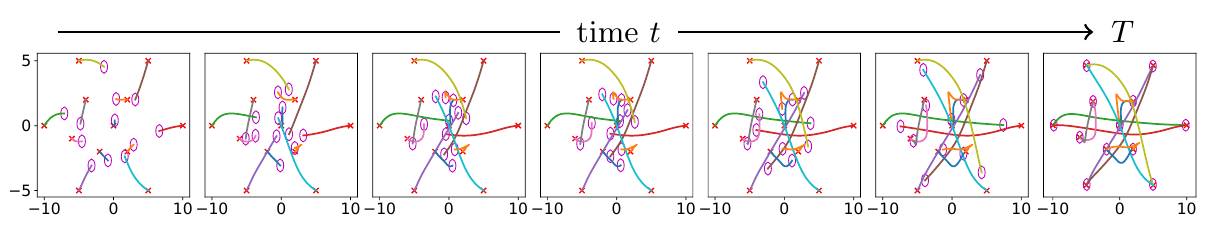}
	\caption{Numerical results of the 12-agent swap experiment (Sec.~\ref{sec:12AgentSwap}). The agents' targets are indicated by red crosses, and the space bubble or safety region around each agent is depicted with a circle. The agents aim to pairwise exchange their positions while avoiding each other and minimizing the length of their trajectories.}
	\label{fig:swap12}
\end{figure*}

\subsubsection{Shocks}\label{subsec:shocks}

	We use this experiment to demonstrate how our approach is robust to shocks to the system's state (Fig.~\ref{fig:shock}). Consider solving the control problem for $s\in[0,T]$ as always. Then for $T=1$, we consider a shock $\bfxi$ (implemented as a random shift) to the system at time $s=0.1$.
	Our method is designed to handle minor shocks that stay within the space of trajectories of the initial distribution about $\bfx_0$. Our model computes a trajectory to $\bfy$ for many initial points. Therefore, for point $\widetilde{\bfx} \in \bfX$, the model provides dynamics $f(s,\bfz_{\widetilde{\bfx}}(s), \bfu_{\widetilde{\bfx}}(s))$ before the shock. After the shock, the state picks up the trajectory of some other point $\hat{\bfx} \in \bfX$ and follows that trajectory to $\bfy$ (Fig.~\ref{fig:shock_a}).
	In this scenario, the total trajectory contains two portions: before the shock and after the shock. That is, 	
	\begin{equation*} \label{eq:shock}
		\begin{split}
			&\bfz_{0,\widetilde{\bfx}}(0.1) = \int_0^{0.1} f  \big( s,\bfz_{0,\widetilde{\bfx}}(s),\bfu_{0,\widetilde{\bfx}}(s) \big) \, \du s, \quad \text{and} \\
			&\bfz_{0,\hat{\bfx}}(1) = \int_{0.1}^{1} f \big( s,\bfz_{0,\hat{\bfx}}(s),\bfu_{0,\hat{\bfx}}(s) \big) \, \du s , \quad  \text{where} \\ 
			& \bfz_{0,\hat{\bfx}}(0.1)=\bfz_{0,\widetilde{\bfx}}(0.1)+\bfxi ,
		\end{split}
	\end{equation*} 
	respectively.
	We view a minor shock then as moving from one trajectory to another (Fig.~\ref{fig:shock_a}). 
	The NN and baseline achieve similar results for the problem along $s{=}[0.1,1]$ (Fig.~\ref{fig:shock_d}).

\begin{figure}
  \centering
    \includegraphics[width=\columnwidth]{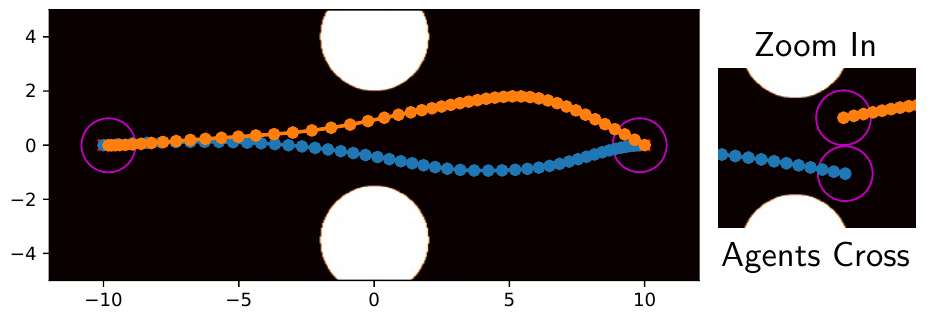}
  \caption{Numerical results of swap experiment with hard-boundary obstacles (Sec.~\ref{sec:TwoAgentSwap}). The agents seek to exchange their positions while keeping a safe distance (indicated by circle) and avoiding the obstacles (white circles). The close-up on the right shows agents at the time of minimal distance.}
  \label{fig:swap2}
\end{figure}

	Interestingly, our model extends outside the training region (Fig.~\ref{fig:shock_b}).
	Although the vast majority of NNs cannot extrapolate, our NN still solves the control problem after a major shock, demonstrating some extrapolation capabilities. We note that the NN solves the original problem for $\bfx_0$ to near optimality. After a large shock, the NN still drives the agents to their targets, although suboptimally. In our example, we compare the NN's solution (Fig.~\ref{fig:shock_b}) with the baseline for $s{=}[0.1,1]$ (Fig.~\ref{fig:shock_c}). The NN learns a solution where agent 2 passes through the corridor followed by agent 1. After the major shock, the NN still applies these dynamics (Fig.~\ref{fig:shock_b}) while the baseline finds a more optimal solution (Fig.~\ref{fig:shock_c}). The NN is roughly 100\% less optimal in this example (Fig.~\ref{fig:shock_d}). 
	
	We attribute the shock robustness to the NN's semi-global nature.
	Experimentally, the shock robustness of our model (Fig.~\ref{fig:shock}) does not noticeably differ from a model trained without penalization (Fig.~\ref{fig:wt_decay}).
	Since the NN is trained offline prior to deployment, it handles shocks in real-time. In contrast, methods that solve for a single trajectory---e.g., the baseline---must pause to recompute following a shock.

\subsubsection{Semi-Global Capabilities of NN model}
	
	For thorough analysis of the NN, we assess one NN's performance for many different initial conditions $(0,\bfx_0+\bfxi)$. We sample 1000 random $\bfxi$ for each magnitude $\| \bfxi \|{=}0.5,1.0,\dots,5.0$. For each $(0,\bfx_0+\bfxi)$, we train a baseline model and compute the suboptimality of the trained NN (Fig.~\ref{fig:globalNN}). This experiment equivalently compares the NN and baseline on samples from concentric hyperspheres. Since a shock can be phrased as picking up a trajectory from an initial condition, testing the NN's semi-global capabilities and shock-robustness are synonymous.
	
	We observe that the NN suboptimality grows as $\| \bfxi \|$ increases (Fig.~\ref{fig:globalNN}). Specifically, for the corridor experiment, the NN performs near optimality within $\| \bfxi \| \leq 2$. Since the NN was trained on $\rho$ which was a Gaussian about $\bfx_0$ with covariance $\bfI$. The bound $\| \bfxi \| \leq 2$ then equates to being within two standard deviations of $\bfx_0$.

\subsection{Multi-Agent Swap Examples} \label{sec:swap}
	
	We present experiments inspired by~\cite{mylvaganam2017differential}, where agents swap positions while avoiding each other. All agents are two-dimensional, and the formulation mostly matches that presented in the corridor example (Sec.~\ref{sec:softcorridor}).
	Specifically, we only alter $\bfx_0$, $\bfy$, and $Q$ for the swap experiments.

\subsubsection{2-Agent Swap} \label{sec:TwoAgentSwap}

	We begin with two agents that swap positions with each other while passing through a corridor with hard edges. To enforce these hard edges, we enforce a space bubble around obstacles similar to how we implement multi-agent interactions~\eqref{eq:Wij}. Therefore, we train with this space bubble but evaluate and plot the results without it. 
	The actual obstacles (two circles with radius 2) are formulated as follows.
	Let $\Omega_{\mathrm{obs}} = \{z \;\; | \;\;   \| z - \bfmu_1 \|< 2  \quad \text{or} \quad \| z - \bfmu_2 \|<  2 \}$, then
	\begin{equation*}
		Q_{i} \Big( z^{(i)} \Big) =
		\begin{cases}
	        1, & \text{if }  z^{(i)} \in \Omega_{\mathrm{obs}},
	        \\
	        0, &\text{otherwise},
		\end{cases}
	\end{equation*}
	where $\bfmu_1 = \begin{bmatrix} 0 \\ 4 \end{bmatrix}$ and $\bfmu_2 = \begin{bmatrix} 0 \\ -3.5 \end{bmatrix}$.
	However, for training, we encode this as 
	\begin{equation*}
		Q_{i,\rm trn} \Big( z^{(i)} \Big) =
		\begin{cases}
		 \sum_{j=1}^2 \eta \left( z^{(i)} \, ; \, \bfmu_j,\bfI \right) , & \text{if } z^{(i)} \in \Omega_{\mathrm{obs,trn}},
		 \\
		 0, &\text{otherwise},
		\end{cases}
	\end{equation*}
	where $\Omega_{\mathrm{obs,trn}} = \{z \; | \;\;   \| z - \bfmu_1 \|< 2.2  \quad \text{or} \quad \| z - \bfmu_2 \|<  2.2 \}$. By training with Gaussian repulsion---which has gradient information within the obstacles---we incentivize the model to learn trajectories avoiding the obstacles. Additionally, $\Omega_{\mathrm{obs,trn}}$ contains an obstacle radial bound ten percent more than in $\Omega_{\mathrm{obs}}$ because we found this additional training buffer alleviates collisions during validation. We use the same obstacle definitions for the baseline and NN approaches.   
	
	For initial and target states, we choose $\bfx_{0}{=}[10,0,-10,0]^{\top}$ and $\bfy{=}[-10,0,10,0]^{\top}$.
	These values are a scaled down version of those in~\cite{mylvaganam2017differential} for ease of visualization.
	For the two-agent problem, the agents successfully switch positions while avoiding each other (Fig.~\ref{fig:swap2}). In validation, the obstacle $Q$ and interaction costs $W$ are exactly 0, so we can confirm that the agents avoid collisions. Qualitatively, our method learns trajectories with shorter arclength than those in~\cite{mylvaganam2017differential}.

\subsubsection{12-Agent Swap} \label{sec:12AgentSwap}

	We also replicate the 12-agent case in \cite{mylvaganam2017differential}. For this experiment, six pairs of agents swap positions. Since there are no obstacles, $Q{=}0$.
	In our setup, the problem is slightly adjusted as our semi-global approach solves for a fixed $\bfy$ but with initial conditions in $\rho$, instead of just $\bfx_{0}$. We display the solution for the single initial case $\bfx_{0}$ (Fig.~\ref{fig:swap12}).

\begin{figure}
\centering
\begin{subfigure}{\columnwidth}
	\includegraphics[width=\columnwidth]{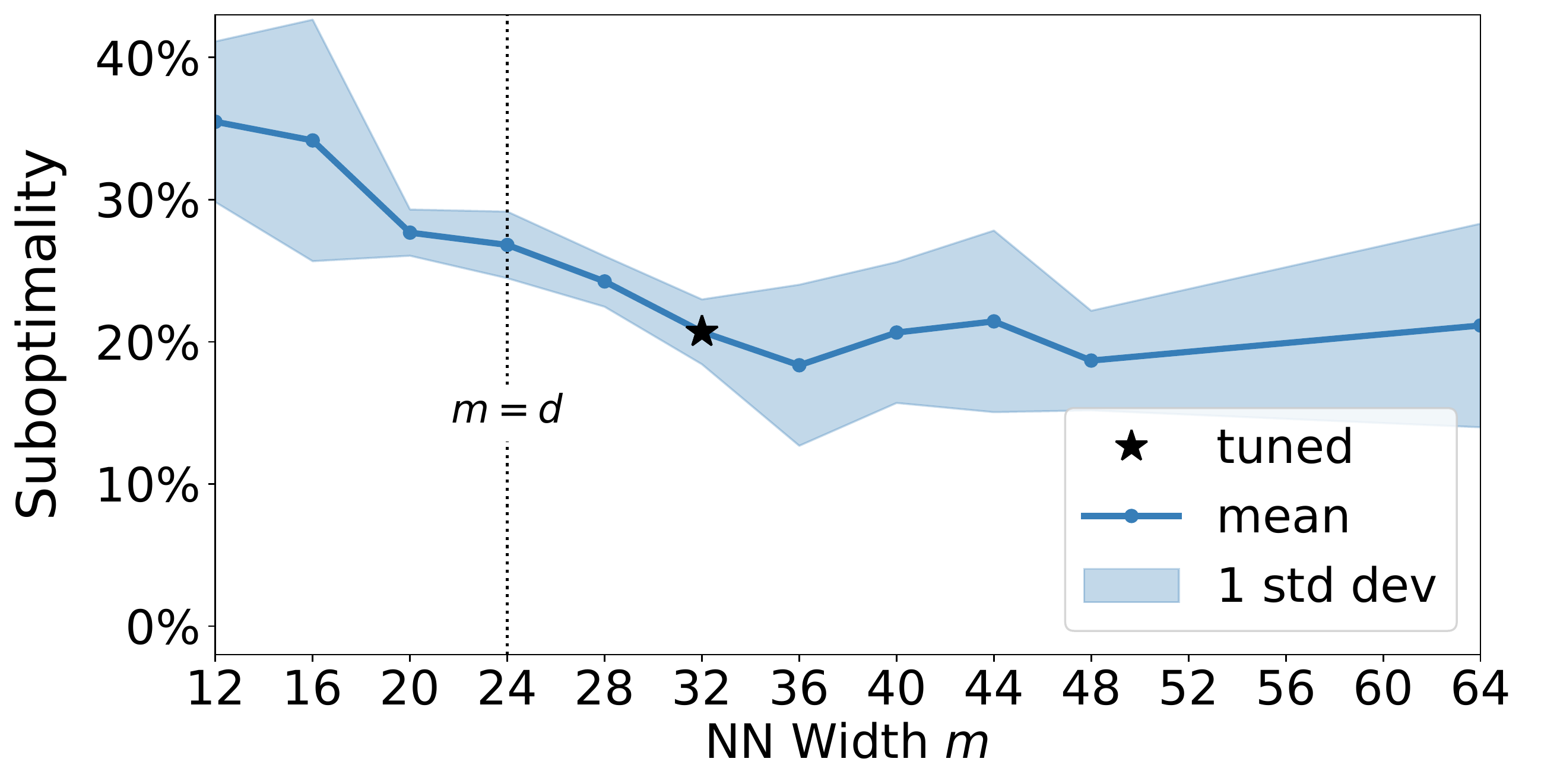} 
	\subcaption{Tuning ResNet width $m$ while keeping all other settings fixed.}
	\label{fig:ablation_m}
\end{subfigure}	
\begin{subfigure}{\columnwidth}
	\centering
	\includegraphics[width=\columnwidth]{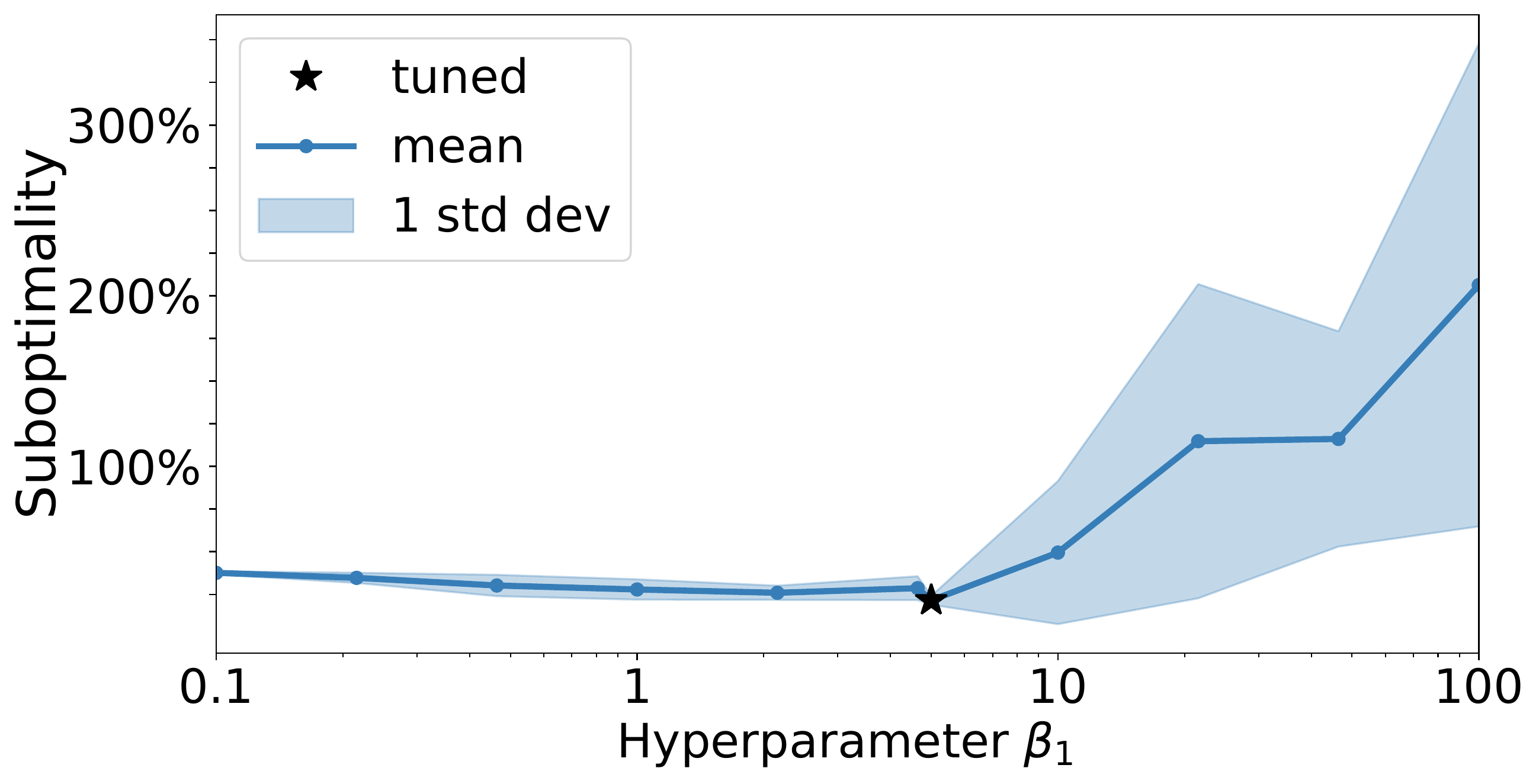} 
	\subcaption{Tuning the scalar hyperparameter $\beta_1$ on the ${\rm HJ_{t}}$ term while keeping all other settings fixed. The x-axis is log-scaled.}
	\label{fig:ablation_beta1}
\end{subfigure}
\begin{subfigure}{\columnwidth}
	\centering
	\includegraphics[width=\columnwidth]{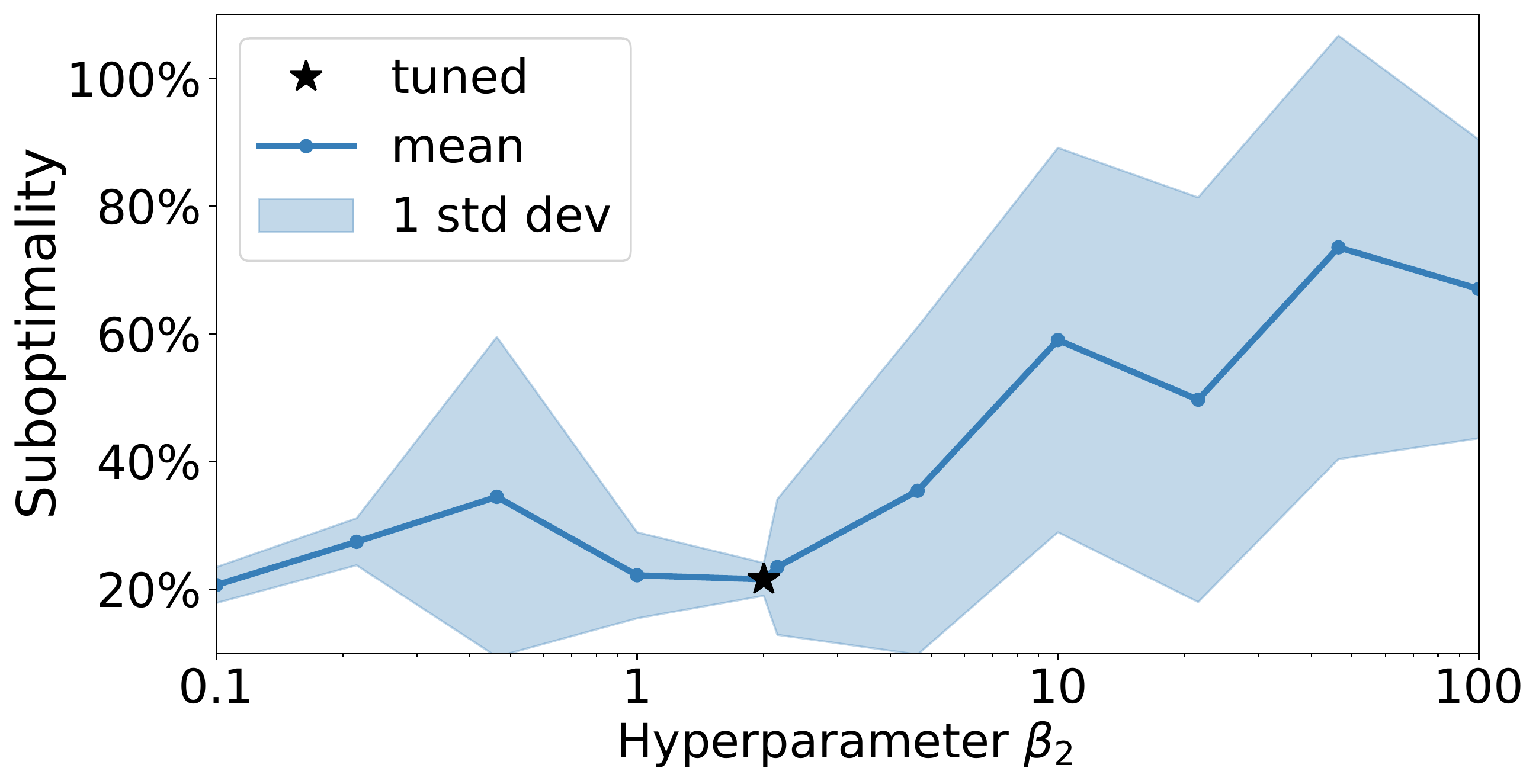} 
	\subcaption{Tuning the scalar hyperparameter $\beta_2$ on the ${\rm HJ_{fin}}$ term while keeping all other settings fixed. The x-axis is log-scaled.}
	\label{fig:ablation_beta2}
\end{subfigure}
\begin{subfigure}{\columnwidth}
	\includegraphics[width=\columnwidth]{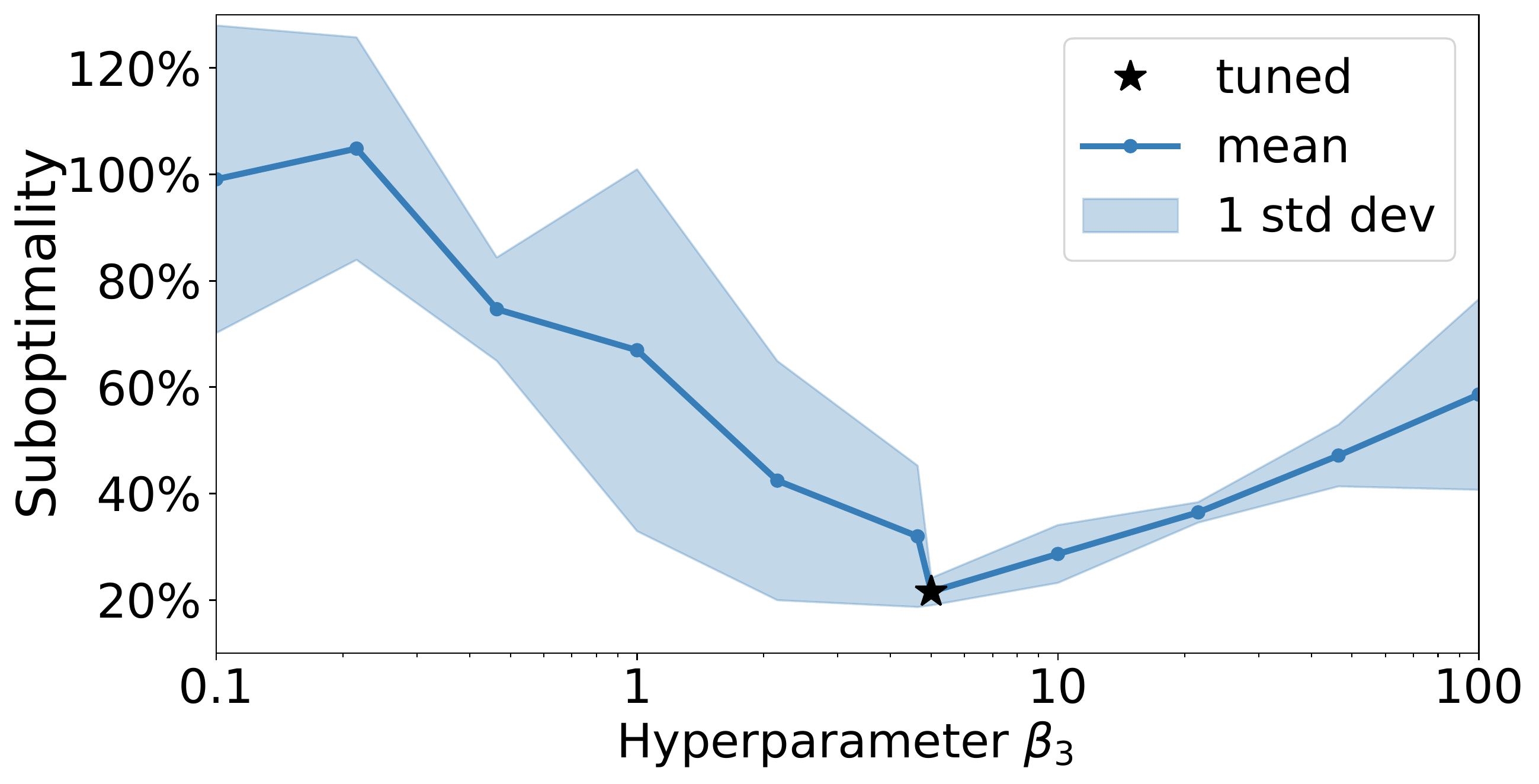} 
	\subcaption{Tuning the scalar hyperparameter $\beta_3$ on the ${\rm HJ_{grad}}$ term while keeping all other settings fixed. The x-axis is log-scaled.}
	\label{fig:ablation_beta3}
\end{subfigure}
\caption{Tuned hyperparameters for the 12-agent swap experiment (Sec.~\ref{sec:12AgentSwap}) where suboptimality is computed relative to the baseline method (Sec.~\ref{sec:baseline}). Each plotted point and bounds are the mean and standard deviation of three model instances.}	
\label{fig:ablation}
\end{figure}

\subsubsection{Impact of ResNet Width}

	We demonstrate the influence of the ResNet width $m$ by observing the results of models with varied width for the 12-agent swap experiment (Fig.~\ref{fig:ablation_m}). We select several $m$ values in the range 12--64 and train three model instances for each while controlling for the rest of the architecture. 
	We then compute the suboptimality of the NN solution relative to the baseline (Sec.~\ref{sec:baseline}) for objective function~\eqref{eq:Joc} of a single initial point $\bfx_0$ (Fig.~\ref{fig:ablation_m}).
	We observe that, for the 12-agent swap experiment, the underlying manifold exists somewhere near dimension $32$ as values $m\geq32$ are relatively stagnant. We note that smaller values of $m$ perform poorly. When $m < d$, we essentially ask the model to condense the input to a lower dimensional manifold. For the 12-agent swap problem, the $d{=}24$ dimensions, though coupled, present no obvious method for reduction to a lower basis. Therefore, we observe poor model performance for $m < d$. 
	
	Based on the experiment (Fig.~\ref{fig:ablation_m}), we use a width of $m=32$ to balance between a small model and performance. We prefer smaller models as a model with few parameters is easier to evaluate. However, due to the GPU parallelization, different model widths in our experiment (Fig.~\ref{fig:ablation_m}) have negligible influence on time per training iteration.

\subsubsection{Hamilton-Jacobi-Bellman Penalty Hyperparameters} \label{sec:beta_tuning}

	In general, we tune hyperparameters relative to each other and set optimizer settings based on architecture design and hyperparameter choices. Thus, in a nuanced response to the findings of Fig.~\ref{fig:penalizers}, we find that, by training longer with an adjusted learning rate scheme, one can achieve a similar NN solution without any HJB penalizers (\textit{cf.} Fig.~\ref{fig:shock},\ref{fig:wt_decay}). This holds because the HJB penalizers do not mathematically alter the problem. 
	
	We design experiments to demonstrate the sensitivity of the NN solution with respect to the hyperparameters $\beta_1,\beta_2,\beta_3$ (Fig.~\ref{fig:ablation}). We train NNs to solve the 12-agent swap experiment (Sec.~\ref{sec:swap}). We check the sensitivity of the NN solution with respect to changing one $\beta$ hyperparameter while keeping all other tuned $\beta$s and hyperparameters fixed (Table~\ref{tab:hyperparameters}). 
	
	$\beta_1\colon$ We observe best performance when $\beta_1 \in (1,5)$ (Fig.~\ref{fig:ablation_beta1}). Since $\beta_1$ weights the ${\rm HJ_{t}}$ term, setting $\beta_1$ too high leads to model training that underprioritizes reaching the target which can result in very suboptimal solutions. 
	
	$\beta_2\colon$ We observe best performance when  $\beta_2 \in (1,2)$ (Fig.~\ref{fig:ablation_beta2}). Since $\beta_2$ weights the ${\rm HJ_{fin}}$ term, setting $\beta_2$ too high leads to NN training that overprioritizes fitting the $\Phi$ value at time $T$. Specifically, the training overprioritizes fitting $\Phi$ rather than $\nabla_{\bfz} \Phi$, which more directly relates to the dynamics. 	
	
	$\beta_3\colon$ We observe best performance when $\beta_3 \in (4,10)$ (Fig.~\ref{fig:ablation_beta3}). Since $\beta_3$ weights the ${\rm HJ_{grad}}$ term, setting $\beta_3$ too high leads to model training that overprioritizes the model reaching the target with less leeway in altering the trajectory for a more optimal $L$. Alternatively, setting $\beta_3$ too small leads to an increase in suboptimality as the model is less likely to satisfactorily reach the target.

\begin{figure}
	\centering
	  \includegraphics[width=0.8\columnwidth]{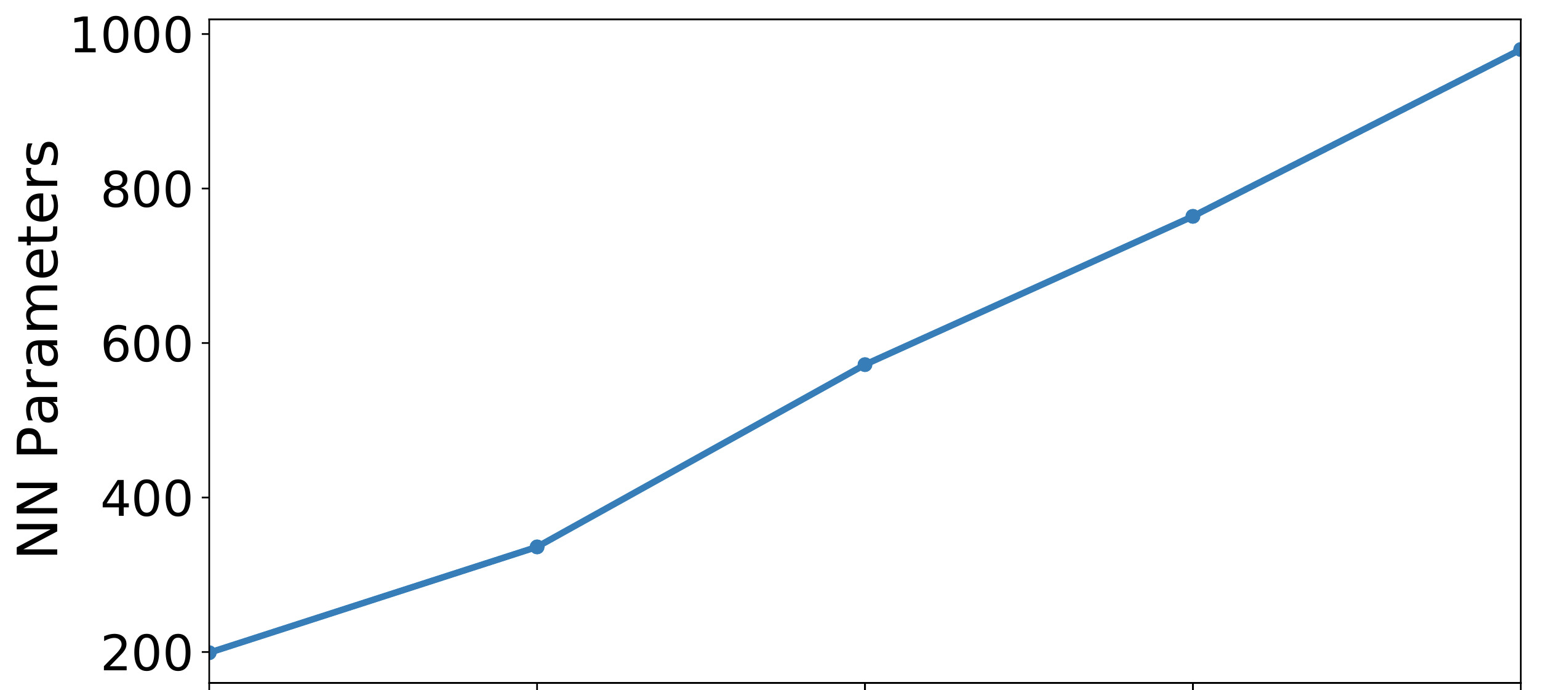} \\
	  \includegraphics[width=0.8\columnwidth]{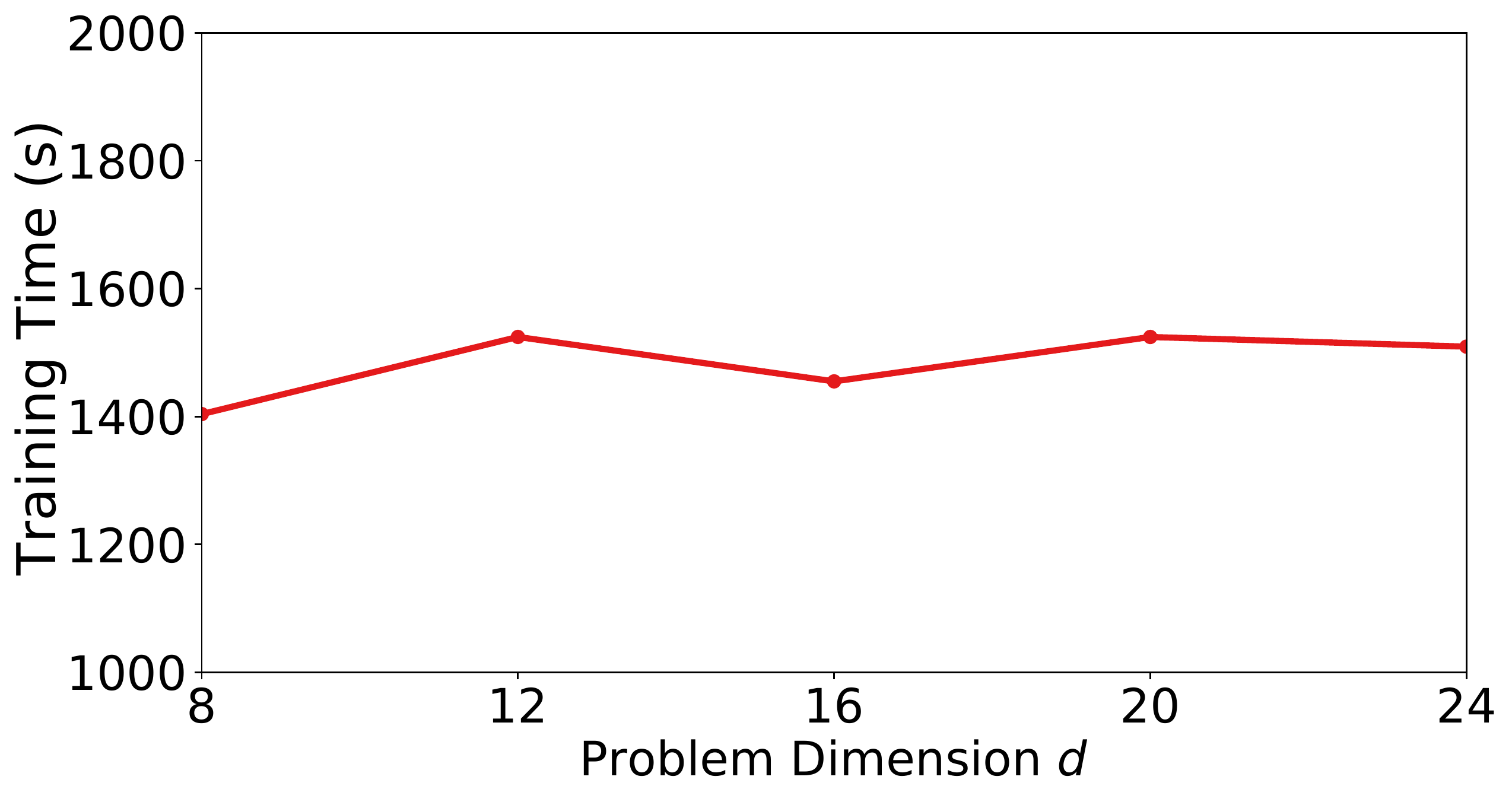} \\
	\caption{The NN's number of parameters scales linearly with the problem dimension as the computational cost remains mostly unchanged, mitigating CoD. For each problem (subproblem of the 12-agent swap experiment), we train the smallest NN that achieves at least 10\% suboptimality.}
	\label{fig:cod}
\end{figure}

\subsubsection{Mitigating the CoD} \label{subsec:cod}

	We expand the 12-agent swap experiment to demonstrate how the NN approach mitigates the CoD (Fig.~\ref{fig:cod}). We design four additional similar problems by removing agents from the original 12-agent version. Thus, we arrive at problems containing 2, 3, 4, 5, and 6 pairs of agents that swap positions. We then determine the smallest NN that is at most 10\% suboptimal. We only tune the width $m$, which dictates the number of NN parameters, and keep all other settings, including the number of training samples and iterations, fixed. The resulting NNs follow a linear growth of number of parameters relative to the problem dimension $d$ (Fig.~\ref{fig:cod}). Due to the parallelization of the GPU, the training times of these NNs remain comparable.

\subsection{Swarm Trajectory Planning Example} \label{sec:swarm}

	We demonstrate the high-dimensional capabilities of our NN approach by solving a 150-dimensional swarm trajectory planning problem in the spirit of~\cite{honig2018trajectory}.
	The swarm problem contains 50 three-dimensional agents that fly from initial to target positions while avoiding each other and obstacles. We construct $Q$ to model two rectangular prism obstacles $[-2,2] \times [-0.5,0.5] \times [0,7]$ and $[2,4] \times [-1,1] \times [0,4]$. 
	We train with Gaussian repulsion inside the obstacles similar to the swap experiment (Sec.~\ref{sec:swap}) and use the same dynamics~\eqref{eq:corridor_f}.
	Due to the complexity of the collision avoidance, we find it beneficial to switch the weights on the HJB penalizers during training---recall that the penalizers do not alter the solution (Sec.~\ref{sec:penalizers}). For the first portion of training, we choose $\beta_1{=}2$, $\beta_2{=}1$, and $\beta_3{=}3$ (Table~\ref{tab:hyperparameters}); for the rest of training, we use $\beta_1{=}\beta_2{=}\beta_3{=}0$. This set-up focuses the model on solving the control problem in the first portion of training as the final-time penalizers help the agents reach their destinations. We then reduce the weights of the penalizers for optimal fine-tuning.
	
\begin{figure}
	\centering
	  \includegraphics[width=0.835\columnwidth]{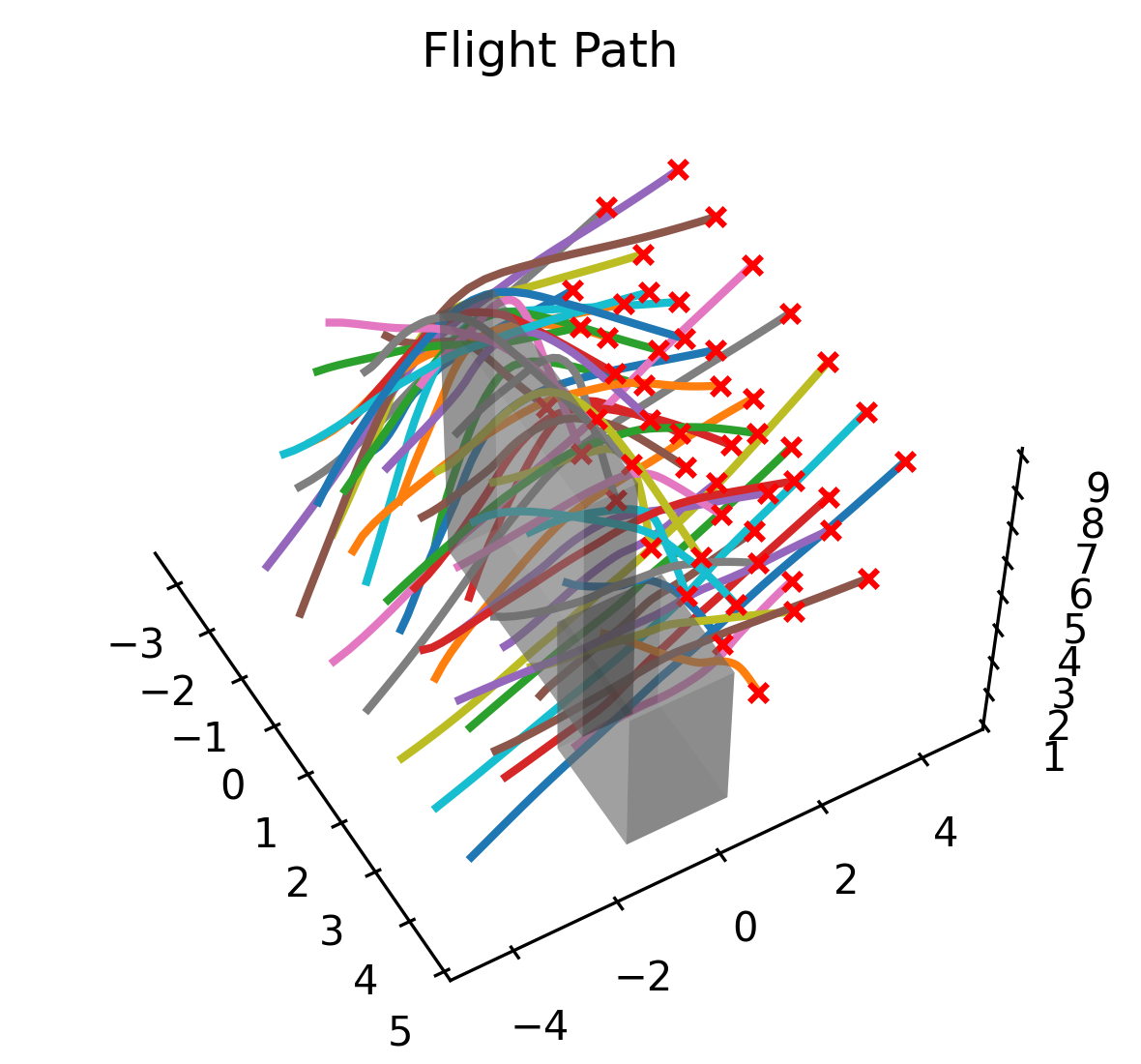}
	\caption{The NN solution for the swarm with 50 agents in $\R^3$ (Sec.~\ref{sec:swarm}). The agents avoid the prism obstacles and each other as they travel from one side of the obstacles to the other.}
	\label{fig:quadtraj}
\end{figure}

 	In validation, we observe that the values for terrain $Q$ and interaction $W$ are exactly 0. Thus, the NN learns to guide all agents around the obstacles and avoid collisions (Fig.~\ref{fig:quadtraj}).

\subsection{Quadcopter Trajectory Planning Example} \label{sec:quadcopter}

In this experiment from~\cite{lin2018splitting}, a quadcopter, i.e., a multirotor helicopter, utilizes its four rotors to propel itself across space from an initial state in the vicinity of $\bfx_{0}$ to target state $\bfy$. 
We choose values $\bfx_{0} {=} [-1.5,-1.5,-1.5,0,\dots,0]^{\top} \in \R^{12}$ and $\bfy {=} [2,2,2,0,\dots,0]^{\top} \in \R^{12}$. 
Denoting gravity as $g$, the acceleration of a quadcopter with mass $m$ is given by%
\begin{equation}
    \begin{cases}
      \ddot{x} = \frac{u}{m} \big( \sin(\psi) \sin(\varphi) + \cos(\psi)\sin(\theta) \cos(\varphi) \big)
      \\
      \ddot{y} = \frac{u}{m} \big(-\cos(\psi) \sin(\varphi) + \sin(\psi)\sin(\theta)\cos(\varphi) \big)
      \\
      \ddot{z} = \frac{u}{m} \cos(\theta)\cos(\varphi) - g
      \\
      \ddot{\psi} = \tau_\psi
      \\
      \ddot{\theta} = \tau_{\theta}
      \\
      \ddot{\varphi} = \tau_{\varphi}
    \end{cases},    
\end{equation}
where $(x,y,z)$ is the spatial position of the quadcopter, $(\psi,\theta,\varphi)$ is the angular orientation with corresponding torques $\tau_{\psi}$, $\tau_{\theta}$, $\tau_{\varphi}$, and $u$ is the main thrust directed out of the bottom of the aircraft~\cite{carrillo2013modeling}.
 The dynamics can be written as the following first-order system
 \begin{equation}
     \label{eq:nablapH}
     \dot{\bfz} = f(s,\bfz,\bfu) \implies \begin{cases}
       \dot{x} = v_x
       \\
       \dot{y} = v_y
       \\
       \dot{z} = v_z
       \\
       \dot{\psi} = v_\psi
       \\
       \dot{\theta} = v_\theta
       \\
       \dot{\varphi} = v_\varphi
       \\
       \dot{v}_x = \frac{u}{m} f_7(\psi, \theta, \varphi)
       \\
       \dot{v}_y = \frac{u}{m} f_8(\psi, \theta, \varphi) 
       \\
       \dot{v}_z = \frac{u}{m} f_9(\theta, \varphi) - g
       \\
       \dot{v}_\psi = \tau_\psi
       \\
       \dot{v}_\theta = \tau_\theta
       \\
       \dot{v}_\varphi = \tau_\varphi
     \end{cases}  ,     
 \end{equation}
 where 
 \begin{equation}
     \begin{split}
         \begin{cases}
             f_7(\psi, \theta, \varphi) &=  \sin(\psi) \sin(\varphi) + \cos(\psi) \sin(\theta) \cos(\varphi)
             \\
             f_8(\psi, \theta, \varphi)  &= -\cos(\psi) \sin(\varphi) + \sin(\psi) \sin(\theta) \cos(\varphi)
             \\
             f_9(\theta, \varphi)  &= \cos(\theta)\cos(\varphi)
         \end{cases}.
     \end{split}
 \end{equation}%
Here, $\bfz = [x \; y \; z \; \psi \; \theta \; \varphi \; v_x \; v_y \; v_z \; v_\psi \; v_\theta \; v_\varphi]^\top \in \R^{12}$ is the state with velocities $v$, and $\bfu = [u \; \tau_\psi \; \tau_\theta \; \tau_\varphi]^\top \in \R^{4}$ is the control. 
 For the energy term, we consider 
 \begin{equation}
 	\begin{split}
     E(\bfu(s)) &= 2 + \| \bfu(s) \|^2  \\ 
     &= 2 + u^2(s) + \tau_\psi^2(s) + \tau_\theta^2(s) + \tau_\varphi^2(s).
     \end{split}
 \end{equation}
 For this problem, we have no obstacles nor other agents, so $L(s,\bfz,\bfu)=E(\bfu)$.

\begin{figure}
	\centering
	\begin{subfigure}{.46\columnwidth}
	\centering
	\includegraphics[width=\linewidth]{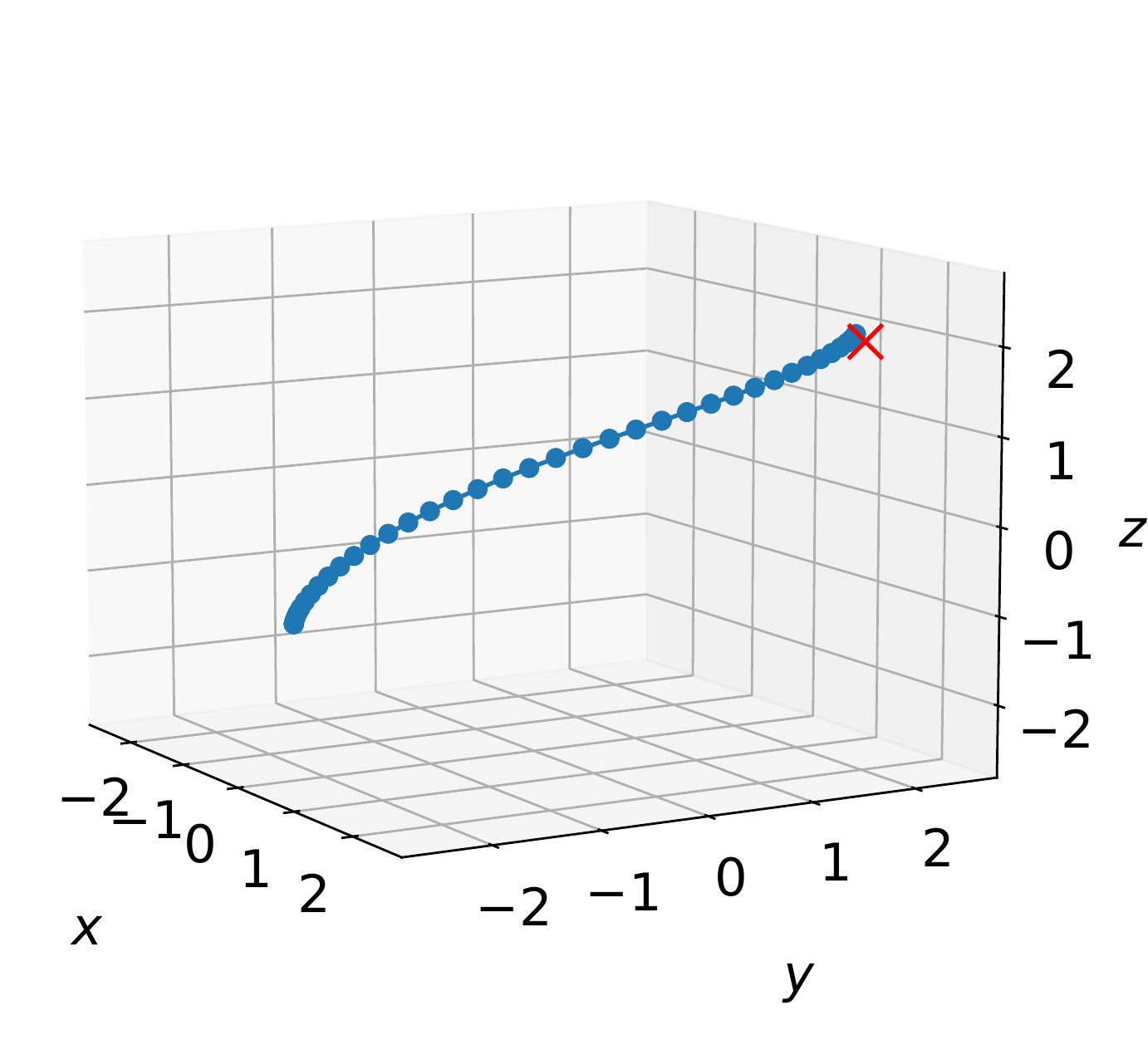}
	\captionsetup{width=0.95\linewidth}%
	\subcaption{Baseline trajectory solved using the four controls on $n_t=50$.}
	\end{subfigure}
	\begin{subfigure}{.46\columnwidth}
	\centering
	\includegraphics[width=\linewidth]{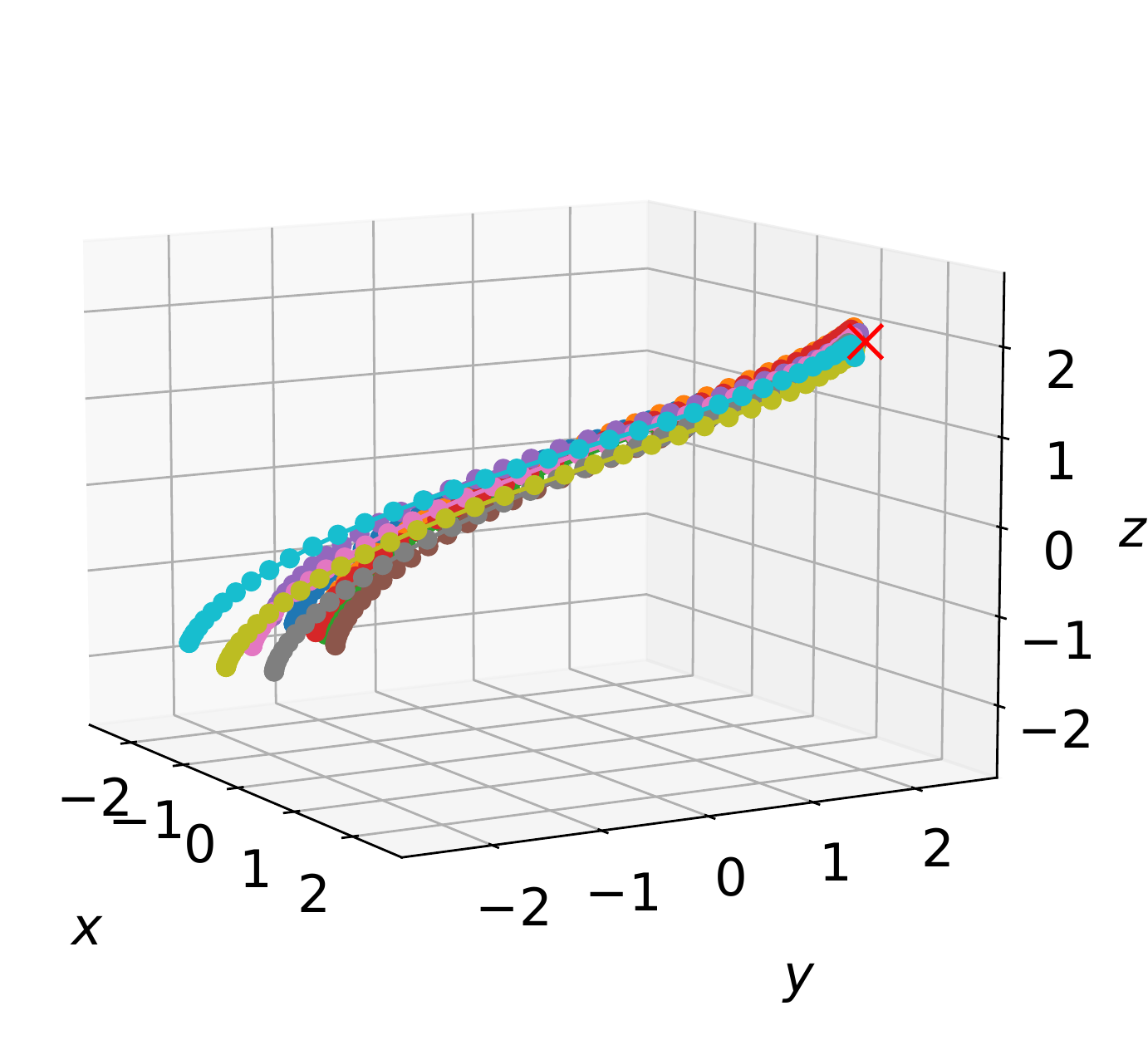}
	\captionsetup{width=0.95\linewidth}%
	\subcaption{NN trajectories, demonstrating the NN's usability for many initial conditions.}
	\end{subfigure}
	\\
	\begin{subfigure}{.46\columnwidth}
	\centering
	\includegraphics[width=\linewidth]{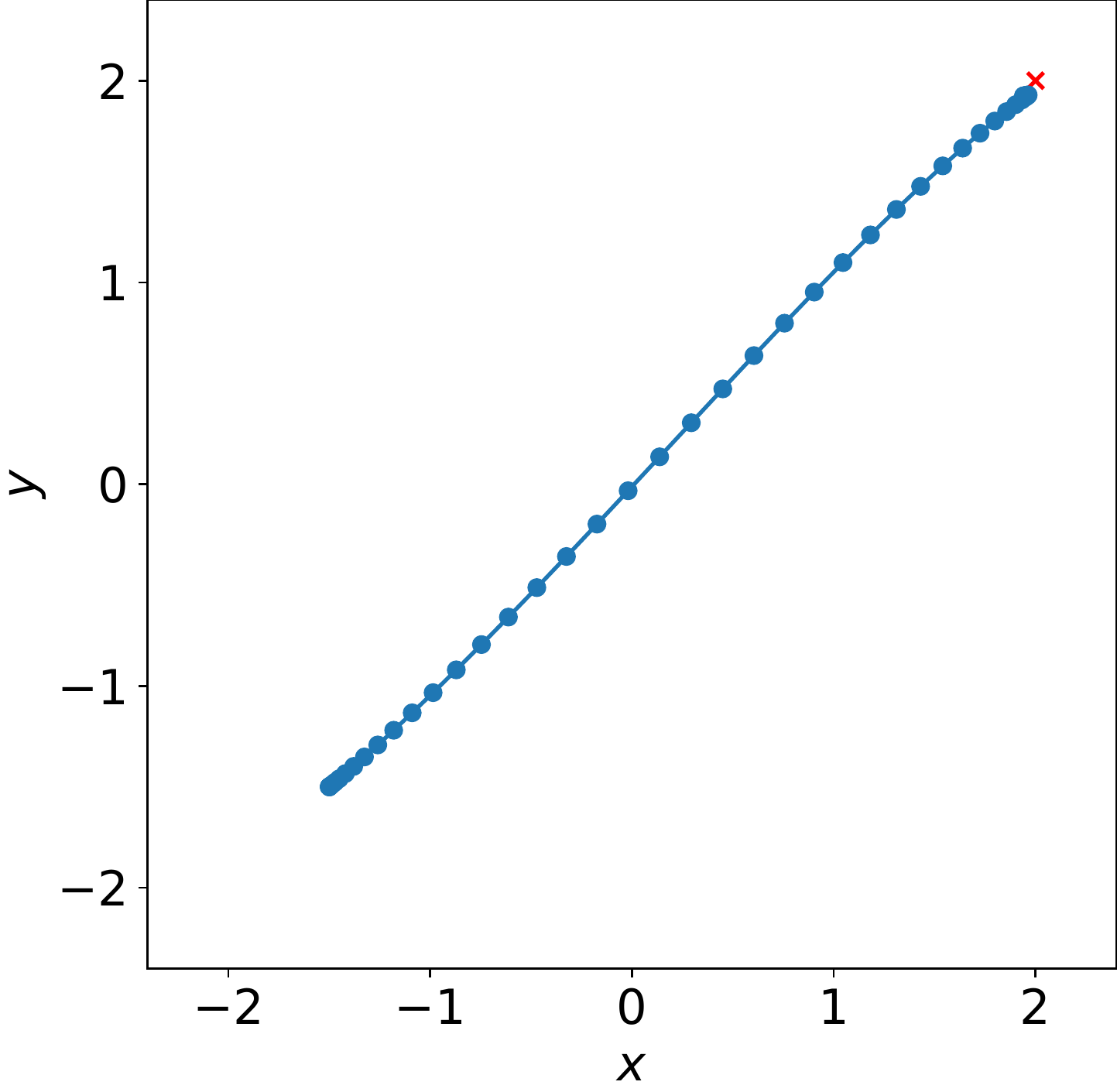}
	\captionsetup{width=0.95\linewidth}%
	\subcaption{Baseline trajectory from bird view.}
	\end{subfigure}
	\begin{subfigure}{.46\columnwidth}
	\centering
	\includegraphics[width=\linewidth]{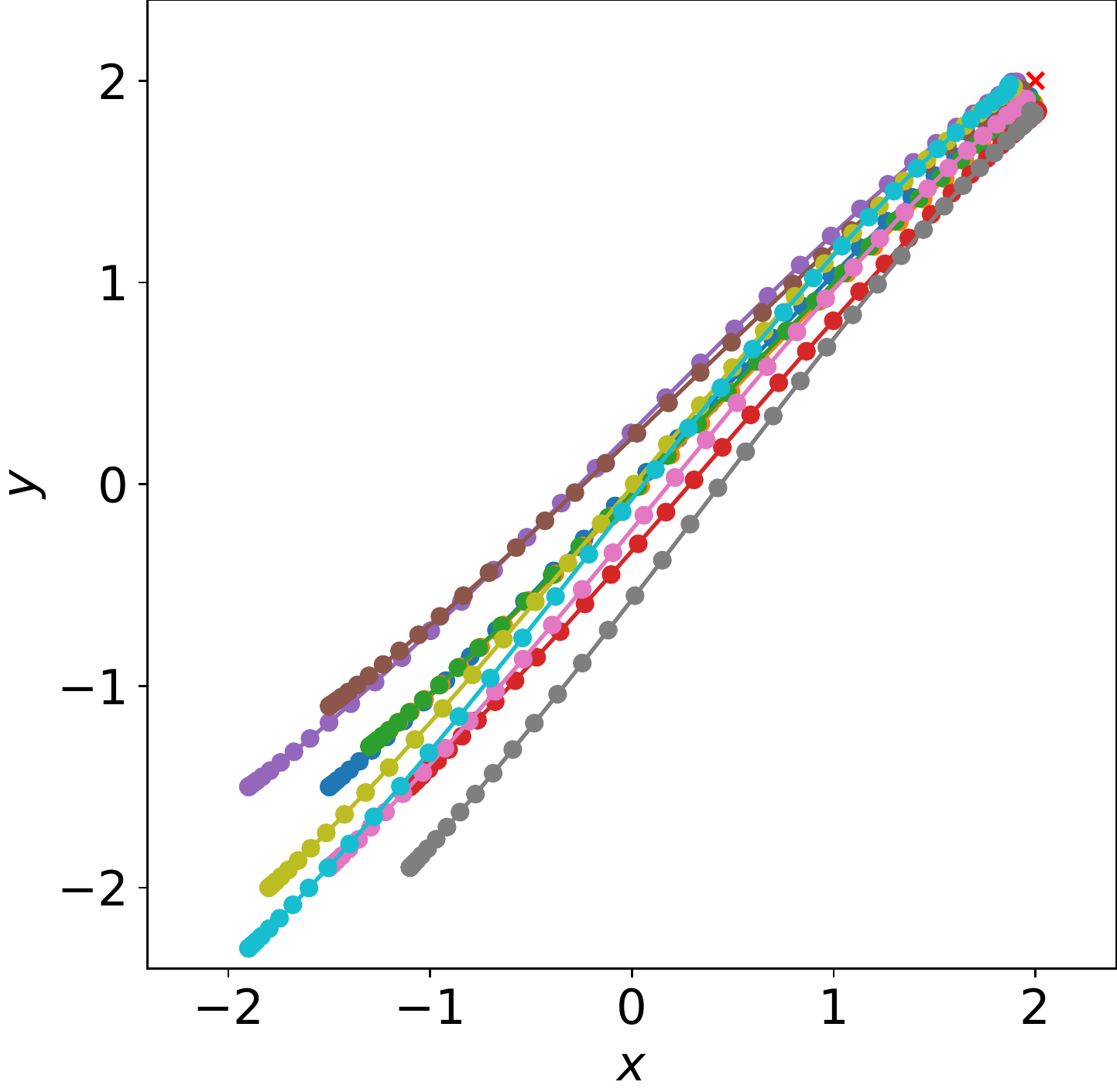}
	\captionsetup{width=0.95\linewidth}%
	\subcaption{NN trajectories from bird view.}
	\end{subfigure}
	\\
	\begin{subfigure}{0.86\columnwidth}
	\centering
	\includegraphics[width=\linewidth]{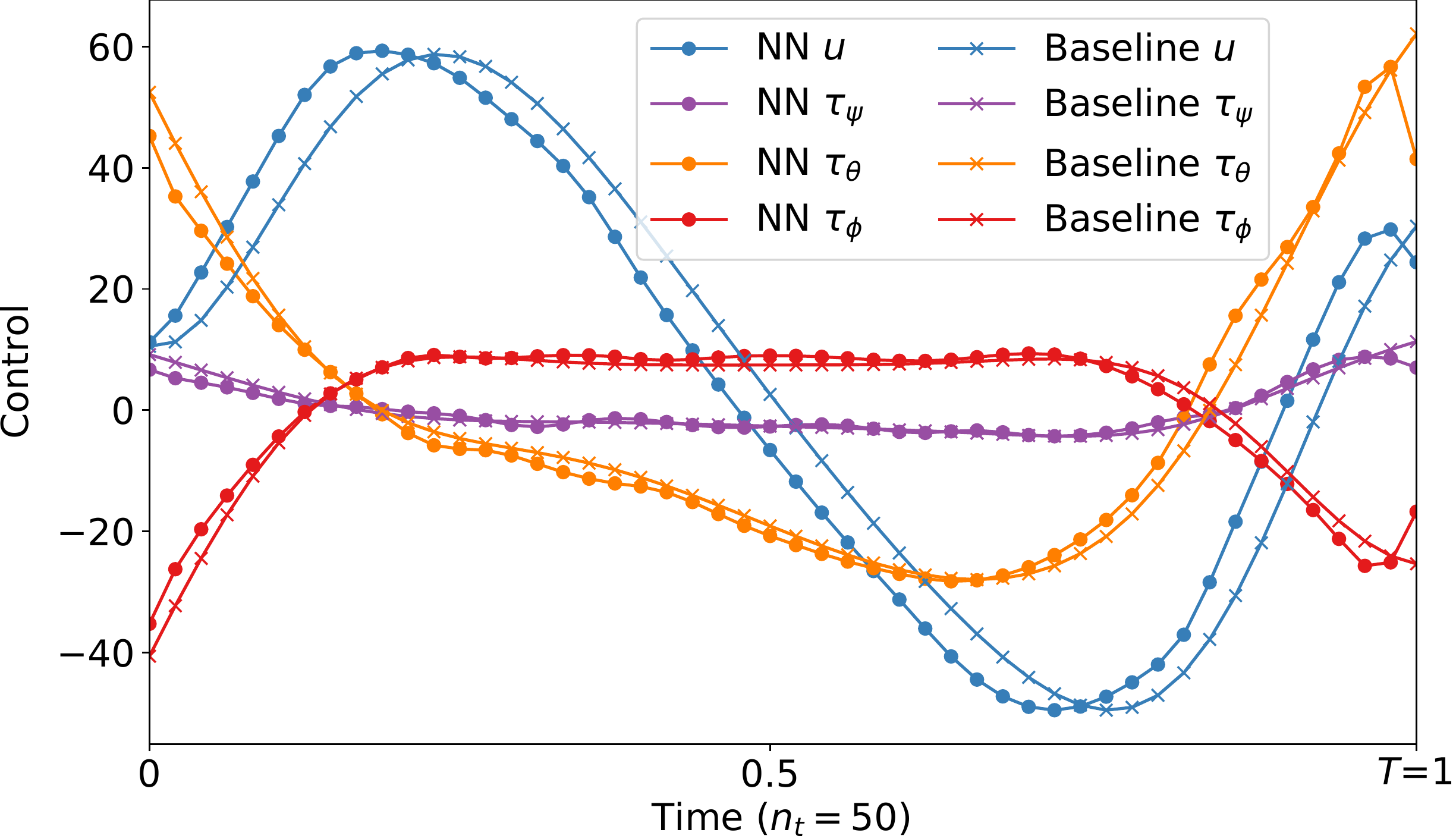}
	\subcaption{Comparison of controls.}
	\label{fig:quad_controls}
	\end{subfigure}
	\\
	\begin{subfigure}{\columnwidth}
	\centering
	\begin{tabular}{lccc}
	   	\toprule
	 	   & $\ell + G$ & $\ell$ & $G$ \\
	 	 \midrule
	     Baseline & 2,182.7 & 2,111.2 & 71.47  \\ 
	     NN       & 2,184.9 & 2,122.0 & 62.90  \\   
	 	\bottomrule
	\end{tabular} 
	\subcaption{Comparison of loss values for single initial point $\bfx_0$.}
	\label{fig:quad_table}
	\end{subfigure}
	\caption{Quadcopter problem results and comparison.}
	\label{fig:quadcopter}
\end{figure}

 We consider the Hamiltonian in~\eqref{eq:H} where $\bfp = [p_1 \; p_2 \;\ldots \; p_{12}]^\top \in \R^{12}$.
 Noting the optimality conditions of \eqref{eq:H} for the quadcopter problem are obtained by
 \begin{equation}
 	\begin{split}
 		-\nabla_{\bfu} E(\bfu) - \bfp^{\top} \, \nabla_{\bfu} f &= \bf0 
 		\\
 		\Rightarrow -2 \begin{bmatrix} u \\ \tau_{\psi} \\ \tau_{\theta} \\  \tau_{\phi} \\ \end{bmatrix} - \begin{bmatrix} p_7 \\ p_8 \\ p_9 \\ p_{10} \\ p_{11} \\ p_{12} \end{bmatrix}^{\top}  
 		\begin{bmatrix} f_7/m & 0 & 0 & 0 \\ f_8/m & 0 & 0 & 0 \\ f_9/m & 0 & 0 & 0 \\ 0  & 1 & 0 & 0\\  0 & 0 & 1 & 0\\ 0 & 0 & 0 & 1 \\ \end{bmatrix} &= \bf0 
 		\\
 		\Rightarrow -2 \begin{bmatrix} u \\ \tau_{\psi} \\ \tau_{\theta} \\  \tau_{\phi} \\ \end{bmatrix} - \begin{bmatrix}  \frac{1}{m}(f_7 p_7 + f_8 p_8 + f_9 p_9) \\ p_{10} \\p_{11} \\ p_{12} \end{bmatrix} &= \bf0,
 	\end{split}
 \end{equation} 
 we can derive an expression for the controls as 
 \begin{equation}
     \label{eq:controlsInP}
     \begin{split}
         u = \frac{-1}{2m}(f_7 p_7 + f_8 p_8 + f_9 p_9), \;\; \\
         \tau_\psi = \frac{-p_{10}}{2}, \;\;
         \tau_\theta = \frac{-p_{11}}{2}, \;\;
         \tau_\varphi = \frac{-p_{12}}{2}.
     \end{split}
 \end{equation}
 We therefore can compute the Hamiltonian
 \begin{equation}
 \begin{split}
 	&H(s,\bfz,\bfp) = -L(\bfu)  - [v_x \; v_y \; v_z] \begin{bmatrix} p_1 \\ p_2 \\ p_3  \end{bmatrix} - [v_\psi \; v_\theta \; v_\varphi] \begin{bmatrix} p_4 \\ p_5 \\ p_6  \end{bmatrix} \\
 	&+ \frac{1}{2m^{2}} \big( p_7 f_7 + p_8 f_8 + p_9 f_9 \big)^{2} + p_9 g+ \frac{1}{2} ( p_{10}^2 + p_{11}^2 + p_{12}^2).
 	\end{split}
 \end{equation} 
 Finally, using~\eqref{eq:gradPhi} and~\eqref{eq:controlsInP}, we compute the controls $\bfu$ using the NN (Fig.~\ref{fig:quad_controls}) with 
 \begin{equation}
     \label{eq:controlsInPhi}
     \begin{split}
         u = \frac{-1}{2m}\left(f_7 \frac{\partial \Phi}{\partial v_x} + f_8 \frac{\partial \Phi}{\partial v_y} + f_9 \frac{\partial \Phi}{\partial v_z} \right), \;\; \\ 
         \tau_\psi = -\frac{1}{2}\frac{\partial \Phi}{\partial v_\psi}, \;\;
         \tau_\theta = -\frac{1}{2}\frac{\partial \Phi}{\partial v_\theta}, \;\;
         \tau_\varphi = -\frac{1}{2}\frac{\partial \Phi}{\partial v_\varphi}.
     \end{split}
 \end{equation}

\begin{figure}
	\centering
	\begin{subfigure}{0.86\columnwidth}
	\centering
	\includegraphics[width=\linewidth]{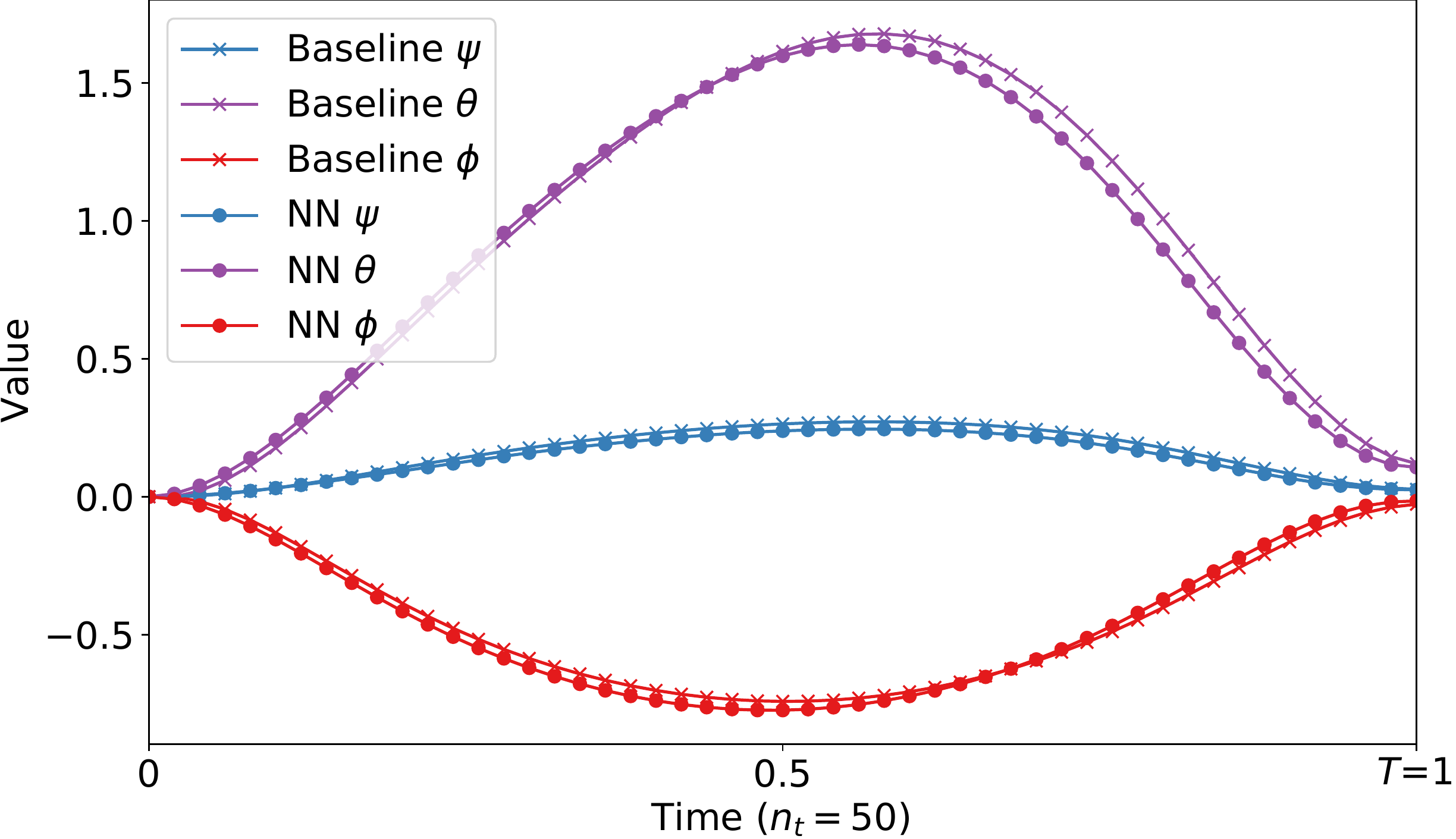}
	\subcaption{Comparison of angular values.}
	\end{subfigure}
	\\
	\begin{subfigure}{0.86\columnwidth}
	\centering
	\includegraphics[width=\linewidth]{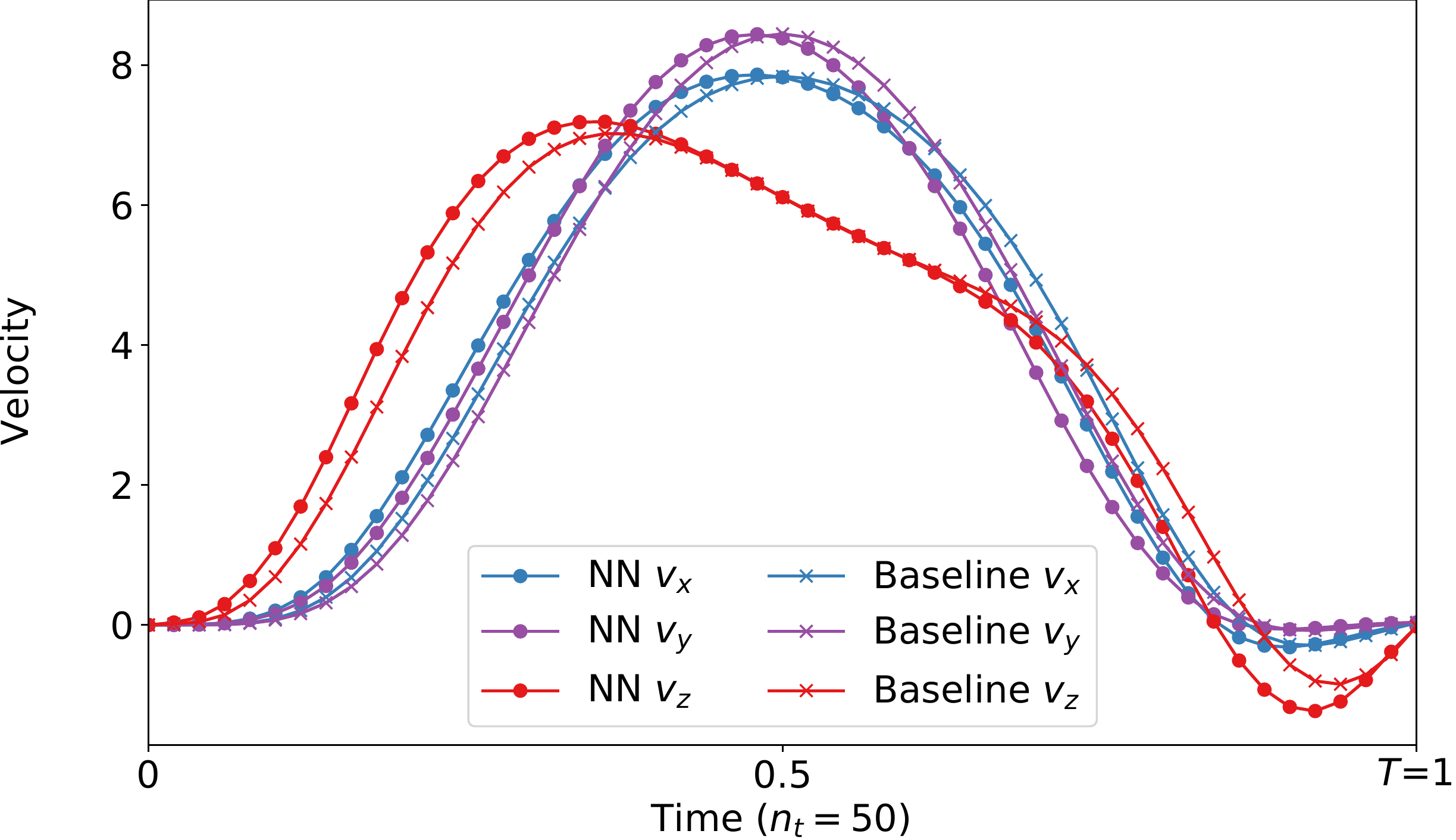}
	\subcaption{Comparison of spatial position velocities.}
	\end{subfigure}
	\\
	\begin{subfigure}{0.86\columnwidth}
	\centering
	\includegraphics[width=\linewidth]{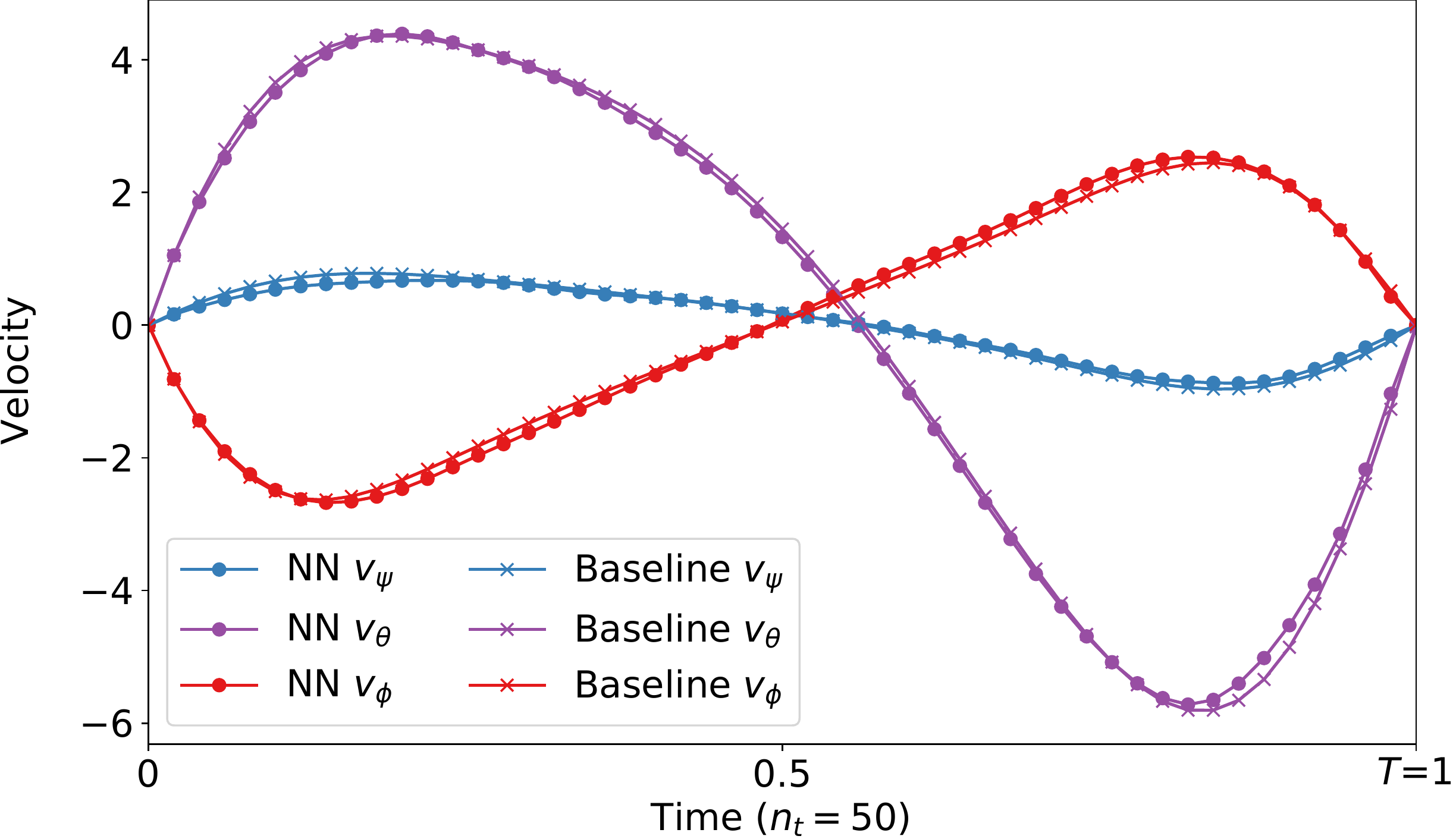}
	\subcaption{Comparison of angular velocities.}
	\end{subfigure}
	\caption{Quadcopter comparison of the additional states as a supplement to Fig.~\ref{fig:quadcopter}.}
	\label{fig:quadcopter_states}
\end{figure}

The quadcopter contains highly coupled 12-dimensional dynamics, which lead to time-consuming model training despite its dimension and lack of obstacles and interactions (Table~\ref{tab:stat}). The HJB terminal conditions seem to offer little impact as no obstacle or interaction costs interfere with the terminal cost.

The NN approach learns similar controls (Fig.~\ref{fig:quad_controls}) and states (Fig.~\ref{fig:quadcopter_states}) as the baseline method. Both methods learn a similar flight path though the NN approach learns for many initial conditions (Fig.~\ref{fig:quadcopter}). As with the corridor problem, the NN learned a solution with better terminal cost, but less optimal $\ell$ than the baseline (Fig.~\ref{fig:quad_table}).

\subsection{Computational Cost} \label{sec:computation}

The offline-online paradigm of our NN approach is specifically designed for efficient deployment in real-time applications. 
To demonstrate this, we compare the computation of the control at time $s$, updating one time step to time $s+1$ on a single CPU core for both methods.
For the baseline, we compute the cost of 100 function and gradient evaluations with $n_t = 20$ for all experiments. NLP algorithms typically require sampling many initial conditions to solve these non-convex problems.
Thus, we believe 100 function-gradient evaluations is a conservative estimate of the cost to generate a trajectory with the baseline method. 
For the NN, we compute the approximate cost of one RK4 step; this is computed by dividing the total cost of the trajectory by the number time steps $n_t$.
Naturally, the set-up of the real-time scenario requires online control generation for a space-time $\bfs$ that may not lie on the pre-computed trajectory. A local solution method, the baseline approach must recompute the entire trajectory from $\bfs$ to target space-time $(\bfy,T)$. Evaluating the NN model to obtain a control is 400x-600x faster than solving a control problem with a baseline method (Table~\ref{tab:stat}).
To compute the full NN trajectory for comparison against the baseline, one simply multiplies the deployment timing of one NN step (Table~\ref{tab:stat}) by the number of time steps $n_t$ (Table~\ref{tab:hyperparameters}).

A thorough timing comparison may involve reducing the NN approach to an RK1 scheme to match the baseline. This would reduce the NN time cost in Table~\ref{tab:stat} by more than a factor of four. Additionally, the baseline timing used a fixed number of time steps $n_t{=}20$. 
In practice, when started later in the time-horizon, the baseline may use $n_t < 20$ and therefore perform faster.
Conversely, timing the baseline with a finer discretization ($n_t > 20$) results in higher time cost.

\section{Discussion} \label{sec:discussion}

Our experiments demonstrate the effectiveness of our NN approach for solving several high-dimensional control problems arising in multi-agent collision avoidance. 
Problems with more complex dynamics and Lagrangians in the finite time-horizon setting are also within reach, as long as the underlying Hamiltonian can be computed efficiently (see Assumption 1). Future work will also involve experimentation with problems that render $\mathcal{H}$ non-concave in $\bfu$ (see Eq.~\ref{eq:H}) and extending our framework to infinite time-horizon control problems such as the ones in~\cite{jiang2015global,michailidis2017adaptive,lopez2019solutions}. Future work will include extensions to situations where the terminal state is unknown or uncertain.
While the focus of this work is on numerical simulations and validations, our positive results motivate the application of our technique to real-world systems.

Our approach does not require first solving for sample trajectories to generate training data and thus differs from the supervised training approaches in~\cite{nakamura2019adaptive,kang2017mitigating}. 
A more similar approach to ours is deep reinforcement learning (RL)~\cite{bertsekas2019reinforcement}.
However, while RL approaches learn from observations of the dynamics and reward functions (and are thus more general), we assume known dynamics and rewards that satisfy Assumption~\ref{assum:feedback}. 
We argue that these assumptions are not too restrictive as demonstrated by our multi-agent examples with non-convex interactions and many examples in the literature, e.g.,~\cite{sanchez2018real}.
We benefit from using the model because the solution's properties---e.g., the approximated value function must satisfy the HJB equations and optimal actions can be obtained from a feedback form---inform the training process.
While this makes our approach less general than RL (applicable in model-free fashion), we expect that this prior knowledge contributes to our approach's effectiveness.
As part of future work, we intend to compare our method to RL approaches in terms of sample-efficiency, network choice, and robustness to hyperparameters.

Among the many RL approaches that have been applied to control problems, the perhaps closest to our method is the actor-critic framework~\cite{konda1999actor1,konda1999actor2}, which employs two neural networks to approximate the policy (actor) and the value function (critic), respectively. 
Notably, the weights of the two networks are not shared, and thus, we should not expect them to generally satisfy the feedback form~\eqref{eq:opt_control_nec}.
In our approach, we parameterize the value function and compute the optimal policy directly using the feedback form.
In addition to requiring only one network, this potentially simplifies the training process, which we plan to investigate in future work.

In the CoD experiment, we observe linear scaling of the NN's parameters for problems of dimensions 8 to 24 (Fig.~\ref{fig:cod}). 
Recall that the number of parameters in a grid-based method scales exponentially with the dimension, leading to prohibitive computational complexity and memory costs.
Since the NN formulation leverages the GPU parallelization and we use the same number of training samples and iterations regardless of dimension, we observe little noticeable change in the time cost across dimensions 8 to 24 (Fig.~\ref{fig:cod}).
Factors that influence the training time stem more from the sequential nature of solving the ODE constraints~\eqref{eq:ODE_cons}. In multi-agent problems, the memory scales quadratically with the number of agents due to the interaction costs $W$. Eventually, for a large enough dimension $d$, the memory costs of the model may exceed the GPU RAM, and implementation changes become necessary.

In our experiments, we show how the semi-global nature of the NN optimally solves the problem within the relevant state-space (Fig.~\ref{fig:globalNN}). As with most machine learning approaches, our method may fail to generalize, i.e., extrapolate beyond the selected training space.  Specifically, the NN often solves the control problem outside the training space, but has potential to do so suboptimally (Fig.~\ref{fig:shock_b}) or cause collisions (Fig.~\ref{fig:globalNN}).

The ability of the NN to avoid collisions and the time needed to train the model  depend most crucially on the number of time steps $n_t$ selected (Sec.~\ref{sec:hyperparams}). 
Large $n_t$ leads to high computation and training time while reducing error; meanwhile, too small $n_t$ leads to overfitting to a refinement of the time discretization of the trajectories. A coarsely discretized approximation of the ODE constraints can result in the model unrealistically jumping over obstacles or other agents. Thus, we use large $n_t$ for the hold-out validation set (Table~\ref{tab:hyperparameters}) to check for overfitting and that the agent movement is sensible.

\section{Conclusion} \label{sec:conclusions}

We formulate and demonstrate an NN approach for solving high-dimensional OC problems arising in multi-agent optimal control that consists of an offline and an online phase. 
In the offline phase, we compute an NN approximation of the control problem's value function in the relevant subset of the space-time domain.
Our learning problem combines the high-dimensional scalability from the PMP and the global nature from the HJB approach.
In the online phase, the NN approximation is used to compute approximately optimal controls using the feedback form in milliseconds.
Our numerical experiments show the effectiveness of our approach for multi-agent problems with state dimension up to 150.
Our experiments show that the obtained controls are nearly optimal relative to a baseline and that the network size and computational costs grow only moderately with the dimension of the problem.
Moreover, our approach is robust to shocks and can handle complicated interaction and obstacle terms.

\appendices

\bibliographystyle{IEEEtran}
\bibliography{main}%

\begin{thebibliography}{10}
\providecommand{\url}[1]{#1}
\csname url@samestyle\endcsname
\providecommand{\newblock}{\relax}
\providecommand{\bibinfo}[2]{#2}
\providecommand{\BIBentrySTDinterwordspacing}{\spaceskip=0pt\relax}
\providecommand{\BIBentryALTinterwordstretchfactor}{4}
\providecommand{\BIBentryALTinterwordspacing}{\spaceskip=\fontdimen2\font plus
\BIBentryALTinterwordstretchfactor\fontdimen3\font minus
  \fontdimen4\font\relax}
\providecommand{\BIBforeignlanguage}[2]{{%
\expandafter\ifx\csname l@#1\endcsname\relax
\typeout{** WARNING: IEEEtran.bst: No hyphenation pattern has been}%
\typeout{** loaded for the language `#1'. Using the pattern for}%
\typeout{** the default language instead.}%
\else
\language=\csname l@#1\endcsname
\fi
#2}}
\providecommand{\BIBdecl}{\relax}
\BIBdecl

\bibitem{honig2018trajectory}
W.~H{\"o}nig, J.~A. Preiss, T.~S. Kumar, G.~S. Sukhatme, and N.~Ayanian,
  ``Trajectory planning for quadrotor swarms,'' \emph{IEEE Transactions on
  Robotics}, vol.~34, no.~4, pp. 856--869, 2018.

\bibitem{kim2020real}
S.~J. Kim and G.~J. Lim, ``A real-time rerouting method for drone flights under
  uncertain flight time,'' \emph{Journal of Intelligent \& Robotic Systems},
  vol. 100, pp. 1355--1368, 2020.

\bibitem{elsayed2020uncertainty}
M.~ElSayed and M.~Mohamed, ``The uncertainty of autonomous unmanned aerial
  vehicles’ energy consumption,'' in \emph{IEEE Transportation
  Electrification Conference \& Expo (ITEC)}, 2020, pp. 8--13.

\bibitem{florence2018nanomap}
P.~R. {Florence}, J.~{Carter}, J.~{Ware}, and R.~{Tedrake}, ``{NanoMap}: Fast,
  uncertainty-aware proximity queries with lazy search over local {3D} data,''
  in \emph{IEEE International Conference on Robotics and Automation (ICRA)},
  2018, pp. 7631--7638.

\bibitem{yang2021opt}
J.~{Yang}, C.~{Liu}, M.~{Coombes}, Y.~{Yan}, and W.~H. {Chen}, ``Optimal path
  following for small fixed-wing uavs under wind disturbances,'' \emph{IEEE
  Transactions on Control Systems Technology (TCST)}, vol.~29, no.~3, pp.
  996--1008, 2021.

\bibitem{ross06realtime}
I.~M. Ross and F.~Fahroo, ``Issues in the real-time computation of optimal
  control,'' \emph{Mathematical and Computer Modelling}, vol.~43, no.~9, pp.
  1172--1188, 2006.

\bibitem{tang2018learning}
G.~Tang, W.~Sun, and K.~Hauser, ``Learning trajectories for real-time optimal
  control of quadrotors,'' in \emph{IEEE/RSJ International Conference on
  Intelligent Robots and Systems (IROS)}, 2018, pp. 3620--3625.

\bibitem{sanchez2018real}
C.~S{\'a}nchez-S{\'a}nchez and D.~Izzo, ``Real-time optimal control via deep
  neural networks: study on landing problems,'' \emph{Journal of Guidance,
  Control, and Dynamics}, vol.~41, no.~5, pp. 1122--1135, 2018.

\bibitem{mo21fasttrack}
M.~Chen, S.~Herbert, H.~Hu, Y.~Pu, J.~Fernandez~Fisac, S.~Bansal, S.~Han, and
  C.~J. Tomlin, ``{FaSTrack}: {A} modular framework for real-time motion
  planning and guaranteed safe tracking,'' \emph{IEEE Transactions on Automatic
  Control (TAC)}, 2021.

\bibitem{bansal2021deepreach}
S.~Bansal and C.~J. Tomlin, ``{DeepReach}: A deep learning approach to
  high-dimensional reachability,'' in \emph{IEEE International Conference on
  Robotics and Automation (ICRA)}, 2021, pp. 1817--1824.

\bibitem{betts98survey}
J.~T. Betts, ``Survey of numerical methods for trajectory optimization,''
  \emph{Journal of Guidance, Control, and Dynamics}, vol.~21, no.~2, pp.
  193--207, 1998.

\bibitem{pontryagin62}
L.~S. Pontryagin, V.~G. Boltyanskii, R.~V. Gamkrelidze, and E.~F. Mishchenko,
  \emph{The Mathematical Theory of Optimal Processes}, ser. Translated by K. N.
  Trirogoff; edited by L. W. Neustadt.\hskip 1em plus 0.5em minus 0.4em\relax
  Interscience Publishers John Wiley \& Sons, Inc.\, New York-London, 1962.

\bibitem{bellman57}
R.~Bellman, \emph{Dynamic Programming}.\hskip 1em plus 0.5em minus 0.4em\relax
  Princeton University Press, Princeton, N. J., 1957.

\bibitem{osher1991high}
S.~Osher and C.-W. Shu, ``High-order essentially nonoscillatory schemes for
  {H}amilton--{J}acobi equations,'' \emph{SIAM Journal on Numerical Analysis},
  vol.~28, no.~4, pp. 907--922, 1991.

\bibitem{ruthotto2020machine}
L.~Ruthotto, S.~J. Osher, W.~Li, L.~Nurbekyan, and S.~W. Fung, ``A machine
  learning framework for solving high-dimensional mean field game and mean
  field control problems,'' \emph{Proceedings of the National Academy of
  Sciences}, vol. 117, no.~17, pp. 9183--9193, 2020.

\bibitem{lin2020apac}
A.~T. Lin, S.~W. Fung, W.~Li, L.~Nurbekyan, and S.~J. Osher, ``Alternating the
  population and control neural networks to solve high-dimensional stochastic
  mean-field games,'' \emph{Proceedings of the National Academy of Sciences},
  vol. 118, no.~31, 2021.

\bibitem{onken2020otflow}
D.~Onken, S.~W. Fung, X.~Li, and L.~Ruthotto, ``{OT-Flow}: Fast and accurate
  continuous normalizing flows via optimal transport,'' \emph{AAAI Conference
  on Artificial Intelligence}, vol.~35, no.~10, pp. 9223--9232, 2021.

\bibitem{onken2020neural}
D.~Onken, L.~Nurbekyan, X.~Li, S.~W. Fung, S.~Osher, and L.~Ruthotto, ``A
  neural network approach applied to multi-agent optimal control,'' in
  \emph{European Control Conference (ECC)}, 2021, pp. 1036--1041.

\bibitem{mylvaganam2017differential}
T.~Mylvaganam, M.~Sassano, and A.~Astolfi, ``A differential game approach to
  multi-agent collision avoidance,'' \emph{IEEE Transactions on Automatic
  Control (TAC)}, vol.~62, no.~8, pp. 4229--4235, 2017.

\bibitem{lin2018splitting}
A.~T. Lin, Y.~T. Chow, and S.~J. Osher, ``A splitting method for overcoming the
  curse of dimensionality in {H}amilton--{J}acobi equations arising from
  nonlinear optimal control and differential games with applications to
  trajectory generation,'' \emph{Communications in Mathematical Sciences},
  vol.~16, no.~7, pp. 1933--1973, 2018.

\bibitem{darbonosher16}
J.~Darbon and S.~Osher, ``Algorithms for overcoming the curse of dimensionality
  for certain {H}amilton--{J}acobi equations arising in control theory and
  elsewhere,'' \emph{Research in the Mathematical Sciences}, vol.~3, no.~1,
  2016.

\bibitem{Kirchner18}
M.~R. {Kirchner}, R.~{Mar}, G.~{Hewer}, J.~{Darbon}, S.~{Osher}, and Y.~T.
  {Chow}, ``Time-optimal collaborative guidance using the generalized {H}opf
  formula,'' \emph{IEEE Control Systems Letters}, vol.~2, no.~2, pp. 201--206,
  2018.

\bibitem{Kirchner2018APM}
M.~R. Kirchner, G.~Hewer, J.~Darbon, and S.~Osher, ``A primal-dual method for
  optimal control and trajectory generation in high-dimensional systems,'' in
  \emph{IEEE Conference on Control Technology and Applications (CCTA)}, 2018,
  pp. 1583--1590.

\bibitem{Chow2018AlgorithmFO}
Y.~Chow, J.~Darbon, S.~Osher, and W.~Yin, ``Algorithm for overcoming the curse
  of dimensionality for certain non-convex {H}amilton-{J}acobi equations,
  projections and differential games,'' in \emph{Annals of Mathematical
  Sciences and Applications}, vol.~3, no.~2.\hskip 1em plus 0.5em minus
  0.4em\relax International Press of Boston, 2018, pp. 369--403.

\bibitem{CHOW2019376}
Y.~T. Chow, J.~Darbon, S.~Osher, and W.~Yin, ``Algorithm for overcoming the
  curse of dimensionality for state-dependent {H}amilton-{J}acobi equations,''
  \emph{Journal of Computational Physics}, vol. 387, pp. 376--409, 2019.

\bibitem{claudel2010lax1}
C.~G. Claudel and A.~M. Bayen, ``Lax--{H}opf based incorporation of internal
  boundary conditions into {H}amilton--{J}acobi equation. {P}art {I}:
  {T}heory,'' \emph{IEEE Transactions on Automatic Control (TAC)}, vol.~55,
  no.~5, pp. 1142--1157, 2010.

\bibitem{claudel2010lax2}
------, ``Lax--{H}opf based incorporation of internal boundary conditions into
  {H}amilton-{J}acobi equation. {P}art {II}: {C}omputational methods,''
  \emph{IEEE Transactions on Automatic Control (TAC)}, vol.~55, no.~5, pp.
  1158--1174, 2010.

\bibitem{kang2017mitigating}
W.~Kang and L.~C. Wilcox, ``Mitigating the curse of dimensionality: Sparse grid
  characteristics method for optimal feedback control and {HJB} equations,''
  \emph{Computational Optimization and Applications}, vol.~68, no.~2, pp.
  289--315, 2017.

\bibitem{nakamura2019adaptive}
T.~Nakamura-Zimmerer, Q.~Gong, and W.~Kang, ``Adaptive deep learning for high
  dimensional {H}amilton-{J}acobi-{B}ellman equations,'' \emph{SIAM Journal on
  Scientific Computing}, vol.~43, no.~2, pp. A1221--A1247, 2021.

\bibitem{sirignano2018}
J.~Sirignano and K.~Spiliopoulos, ``{DGM}: {A} deep learning algorithm for
  solving partial differential equations,'' \emph{Journal of Computational
  Physics}, vol. 375, p. 1339–1364, Dec 2018.

\bibitem{kunisch2020semiglobal}
K.~Kunisch and D.~Walter, ``Semiglobal optimal feedback stabilization of
  autonomous systems via deep neural network approximation,'' \emph{ESAIM:
  Control, Optimisation and Calculus of Variations}, vol.~27, 2021.

\bibitem{E_2017}
W.~E, J.~Han, and A.~Jentzen, ``Deep learning-based numerical methods for
  high-dimensional parabolic partial differential equations and backward
  stochastic differential equations,'' \emph{Communications in Mathematics and
  Statistics}, vol.~5, no.~4, pp. 349–--380, Nov 2017.

\bibitem{Han_2018}
J.~Han, A.~Jentzen, and W.~E, ``Solving high-dimensional partial differential
  equations using deep learning,'' \emph{Proceedings of the National Academy of
  Sciences}, vol. 115, no.~34, pp. 8505--8510, Aug 2018.

\bibitem{nusken2020solving}
N.~N{\"u}sken and L.~Richter, ``Solving high-dimensional
  {H}amilton-{J}acobi-{B}ellman {PDE}s using neural networks: Perspectives from
  the theory of controlled diffusions and measures on path space,''
  \emph{Partial Differential Equations and Applications}, vol.~2, no.~4, pp.
  1--48, 2021.

\bibitem{moon2020generalized}
J.~Moon, ``Generalized risk-sensitive optimal control and
  {H}amilton-{J}acobi-{B}ellman equation,'' \emph{IEEE Transactions on
  Automatic Control (TAC)}, 2020.

\bibitem{han2016deep}
J.~Han and W.~E, ``Deep learning approximation for stochastic control
  problems,'' \emph{arXiv:1611.07422}, 2016.

\bibitem{stern2019multi}
R.~Stern, N.~R. Sturtevant, A.~Felner, S.~Koenig, H.~Ma, T.~T. Walker, J.~Li,
  D.~Atzmon, L.~Cohen, T.~S. Kumar \emph{et~al.}, ``Multi-agent pathfinding:
  Definitions, variants, and benchmarks,'' \emph{Twelfth Annual Symposium on
  Combinatorial Search}, 2019.

\bibitem{jing2019multiagent}
G.~Jing and L.~Wang, ``Multiagent flocking with angle-based formation shape
  control,'' \emph{IEEE Transactions on Automatic Control (TAC)}, vol.~65,
  no.~2, pp. 817--823, 2019.

\bibitem{zhao2018affine}
S.~Zhao, ``Affine formation maneuver control of multiagent systems,''
  \emph{IEEE Transactions on Automatic Control (TAC)}, vol.~63, no.~12, pp.
  4140--4155, 2018.

\bibitem{sharon2015conflict}
G.~Sharon, R.~Stern, A.~Felner, and N.~R. Sturtevant, ``Conflict-based search
  for optimal multi-agent pathfinding,'' \emph{Artificial Intelligence}, vol.
  219, pp. 40--66, 2015.

\bibitem{wagner2011m}
G.~Wagner and H.~Choset, ``M*: A complete multirobot path planning algorithm
  with performance bounds,'' in \emph{{IEEE/RSJ} International Conference on
  Intelligent Robots and Systems}, 2011, pp. 3260--3267.

\bibitem{erdmann1987multiple}
M.~Erdmann and T.~Lozano-Perez, ``On multiple moving objects,''
  \emph{Algorithmica}, vol.~2, no. 1--4, 1987.

\bibitem{richards2002aircraft}
A.~Richards and J.~P. How, ``Aircraft trajectory planning with collision
  avoidance using mixed integer linear programming,'' in \emph{American Control
  Conference}, vol.~3.\hskip 1em plus 0.5em minus 0.4em\relax IEEE, 2002, pp.
  1936--1941.

\bibitem{blackmore2006optimal}
L.~Blackmore and B.~Williams, ``Optimal manipulator path planning with
  obstacles using disjunctive programming,'' in \emph{American Control
  Conference}.\hskip 1em plus 0.5em minus 0.4em\relax IEEE, 2006.

\bibitem{patel2011trajectory}
R.~B. Patel and P.~J. Goulart, ``Trajectory generation for aircraft avoidance
  maneuvers using online optimization,'' \emph{Journal of Guidance, Control,
  and Dynamics}, vol.~34, no.~1, pp. 218--230, 2011.

\bibitem{zhang2021opt}
X.~{Zhang}, A.~{Liniger}, and F.~{Borrelli}, ``Optimization-based collision
  avoidance,'' \emph{IEEE Transactions on Control Systems Technology (TCST)},
  vol.~29, no.~3, pp. 972--983, 2021.

\bibitem{standley2011complete}
T.~Standley and R.~Korf, ``Complete algorithms for cooperative pathfinding
  problems,'' in \emph{International Joint Conference on Artificial
  Intelligence (IJCAI)}, 2011, pp. 668--673.

\bibitem{wagner2015subdimensional}
G.~Wagner and H.~Choset, ``Subdimensional expansion for multirobot path
  planning,'' \emph{Artificial Intelligence}, vol. 219, pp. 1--24, 2015.

\bibitem{boyarski2015icbs}
E.~Boyarski, A.~Felner, R.~Stern, G.~Sharon, O.~Betzalel, D.~Tolpin, and
  E.~Shimony, ``{ICBS}: The improved conflict-based search algorithm for
  multi-agent pathfinding,'' in \emph{Eighth Annual Symposium on Combinatorial
  Search}.\hskip 1em plus 0.5em minus 0.4em\relax Citeseer, 2015.

\bibitem{Cohen16Improved}
L.~Cohen, T.~Uras, T.~K.~S. Kumar, H.~Xu, N.~Ayanian, and S.~Koenig, ``Improved
  solvers for bounded-suboptimal multi-agent path finding,'' in
  \emph{International Joint Conference on Artificial Intelligence (IJCAI)},
  2016, pp. 3067--3074.

\bibitem{riviere2020glas}
B.~Rivi{\`e}re, W.~H{\"o}nig, Y.~Yue, and S.-J. Chung, ``{GLAS}:
  Global-to-local safe autonomy synthesis for multi-robot motion planning with
  end-to-end learning,'' \emph{IEEE Robotics and Automation Letters}, vol.~5,
  no.~3, pp. 4249--4256, 2020.

\bibitem{shi2020neuralswarm}
G.~Shi, W.~H{\"o}nig, Y.~Yue, and S.-J. Chung, ``{Neural-Swarm}: Decentralized
  close-proximity multirotor control using learned interactions,'' \emph{IEEE
  International Conference on Robotics and Automation (ICRA)}, pp. 3241--3247,
  2020.

\bibitem{flemingsoner06}
W.~H. Fleming and H.~M. Soner, \emph{Controlled {M}arkov Processes and
  Viscosity Solutions}, 2nd~ed., ser. Stochastic Modelling and Applied
  Probability.\hskip 1em plus 0.5em minus 0.4em\relax Springer, New York, 2006,
  vol.~25.

\bibitem{kang2019algorithms}
W.~Kang, Q.~Gong, and T.~Nakamura-Zimmerer, ``Algorithms of data development
  for deep learning and feedback design,'' \emph{arXiv:1912.00492}, 2019.

\bibitem{cannarsa04}
P.~Cannarsa and C.~Sinestrari, \emph{Semiconcave Functions, {H}amilton-{J}acobi
  Equations, and Optimal Control}, ser. Progress in Nonlinear Differential
  Equations and their Applications.\hskip 1em plus 0.5em minus 0.4em\relax
  Boston, MA: Birkh\"{a}user Boston, Inc., 2004, vol.~58.

\bibitem{evans2010partial}
L.~C. Evans, \emph{Partial Differential Equations}.\hskip 1em plus 0.5em minus
  0.4em\relax American Mathematical Society, 2010, vol.~19.

\bibitem{crandall1983viscosity}
M.~G. Crandall and P.-L. Lions, ``Viscosity solutions of {H}amilton-{J}acobi
  equations,'' \emph{Trans. Amer. Math. Soc.}, vol. 277, no.~1, pp. 1--42,
  1983.

\bibitem{finlay2020train}
C.~Finlay, J.-H. Jacobsen, L.~Nurbekyan, and A.~M. Oberman, ``How to train your
  neural {ODE}: the world of {J}acobian and kinetic regularization,'' in
  \emph{International Conference on Machine Learning (ICML)}, 2020, pp.
  3154--3164.

\bibitem{he2016deep}
K.~He, X.~Zhang, S.~Ren, and J.~Sun, ``Deep residual learning for image
  recognition,'' in \emph{IEEE Conference on Computer Vision and Pattern
  Recognition (CVPR)}, 2016, pp. 770--778.

\bibitem{gholami2019anode}
A.~Gholaminejad, K.~Keutzer, and G.~Biros, ``{ANODE}: Unconditionally accurate
  memory-efficient gradients for neural {ODE}s,'' in \emph{International Joint
  Conference on Artificial Intelligence (IJCAI)}, 2019, pp. 730--736.

\bibitem{onken2020do}
D.~Onken and L.~Ruthotto, ``Discretize-optimize vs. optimize-discretize for
  time-series regression and continuous normalizing flows,''
  \emph{arXiv:2005.13420}, 2020.

\bibitem{nocedal2006numerical}
J.~Nocedal and S.~Wright, \emph{Numerical Optimization}.\hskip 1em plus 0.5em
  minus 0.4em\relax Springer Science \& Business Media, 2006.

\bibitem{kingma2014adam}
D.~P. Kingma and J.~Ba, ``Adam: {A} method for stochastic optimization,'' in
  \emph{International Conference on Learning Representations ({ICLR})}, 2015.

\bibitem{yang2019}
L.~Yang and G.~E. Karniadakis, ``Potential flow generator with ${L}_2$ optimal
  transport regularity for generative models,'' \emph{IEEE Transactions on
  Neural Networks and Learning Systems}, 2020.

\bibitem{carrillo2013modeling}
L.~R.~G. Carrillo, A.~E.~D. L{\'o}pez, R.~Lozano, and C.~P{\'e}gard, ``Modeling
  the quad-rotor mini-rotorcraft,'' in \emph{Quad Rotorcraft Control}.\hskip
  1em plus 0.5em minus 0.4em\relax Springer, 2013, pp. 23--34.

\bibitem{jiang2015global}
Y.~Jiang and Z.-P. Jiang, ``Global adaptive dynamic programming for
  continuous-time nonlinear systems,'' \emph{IEEE Transactions on Automatic
  Control (TAC)}, vol.~60, no.~11, pp. 2917--2929, 2015.

\bibitem{michailidis2017adaptive}
I.~Michailidis, S.~Baldi, E.~B. Kosmatopoulos, and P.~A. Ioannou, ``Adaptive
  optimal control for large-scale nonlinear systems,'' \emph{IEEE TAC},
  vol.~62, no.~11, pp. 5567--5577, 2017.

\bibitem{lopez2019solutions}
V.~G. Lopez, F.~L. Lewis, Y.~Wan, E.~N. Sanchez, and L.~Fan, ``Solutions for
  multiagent pursuit-evasion games on communication graphs: Finite-time capture
  and asymptotic behaviors,'' \emph{IEEE Transactions on Automatic Control
  (TAC)}, vol.~65, no.~5, pp. 1911--1923, 2019.

\bibitem{bertsekas2019reinforcement}
D.~P. Bertsekas, \emph{Reinforcement Learning and Optimal Control}.\hskip 1em
  plus 0.5em minus 0.4em\relax Athena Scientific Belmont, MA, 2019.

\bibitem{konda1999actor1}
V.~R. Konda and V.~S. Borkar, ``Actor-critic--type learning algorithms for
  {M}arkov decision processes,'' \emph{SIAM Journal on Control and
  Optimization}, vol.~38, no.~1, pp. 94--123, 1999.

\bibitem{konda1999actor2}
V.~Konda and J.~Tsitsiklis, ``Actor-{C}ritic algorithms,'' \emph{Advances in
  Neural Information Processing Systems (NeurIPS)}, vol.~12, 1999.

\end{thebibliography}

\end{document}